\documentclass{amsart}

\usepackage{graphicx}
\usepackage{bbm,amsmath,amsfonts,amssymb}
\usepackage[mathscr]{eucal}

\usepackage{natbib}

\newtheorem{theorem}{Theorem}
\newtheorem{lemma}{Lemma}
\newtheorem{corollary}{Corollary}
\newtheorem{proposition}{Proposition}
\newtheorem{definition}{Definition}

\newtheorem{remark}{Remark}

\newcounter{figures}[section]

\RequirePackage[colorlinks,allcolors =black]{hyperref}

\newtheorem{assumption*}{Assumption}

\renewcommand{\kappa}{\varkappa}

\newcommand{\rd}{{\rm d}}
\newcommand{\e}{\varepsilon}

\newcommand{\cA}{{\mathcal A}}
\newcommand{\cB}{{\mathcal B}}
\newcommand{\cC}{{\mathcal C}}

\newcommand{\cE}{{\mathcal E}}
\newcommand{\cF}{{\mathcal F}}
\newcommand{\cG}{{\mathcal G}}
\newcommand{\cH}{{\mathcal H}}

\newcommand{\cJ}{{\mathcal J}}
\newcommand{\cK}{{\mathcal K}}

\newcommand{\cN}{{\mathcal N}}
\newcommand{\cO}{{\mathcal O}}

\newcommand{\cR}{{\mathcal R}}
\newcommand{\cS}{{\mathcal S}}

\newcommand{\cU}{{\mathcal U}}

\newcommand{\cW}{{\mathcal W}}

\newcommand{\cZ}{{\mathcal Z}}

\newcommand{\Be}{\boldsymbol{e}}

\newcommand{\blb}{\boldsymbol{\beta}}

\newcommand{\blg}{\boldsymbol{\mu}}

\newcommand{\bbr}{\boldsymbol{r}}
\newcommand{\bB}{\mathbb B}

\newcommand{\bE}{\mathbb E}
\newcommand{\bF}{\mathbb F}

\newcommand{\bK}{{\mathbb K}}
\newcommand{\bL}{{\mathbb L}}

\newcommand{\bN}{{\mathbb N}}

\newcommand{\bP}{{\mathbb P}}

\newcommand{\bR}{{\mathbb R}}

\newcommand{\bZ}{{\mathbb Z}}
\newcommand{\mF}{\mathfrak{F}}
\newcommand{\mS}{\mathfrak{S}}

\newcommand{\mN}{\mathfrak{N}}

\newcommand{\mU}{\mathfrak{U}}

\newcommand{\mh}{\mathfrak{h}}

\newcommand{\mH}{\mathfrak{H}}

\newcommand{\mb}{\mathfrak{b}}
\newcommand{\mz}{\mathfrak{z}}

\newcommand{\epr}{\hfill\hbox{\hskip 4pt
                \vrule width 5pt height 6pt depth 1.5pt}\vspace{0.5cm}\par}

\newcommand{\bsh}{\boldsymbol{\vec{h}}}

\begin{document}

\title[Adaptive upper bounds]
{Adaptive estimation of the $\bL_2$-norm of a probability density and related topics II. Upper bounds via the oracle approach.}

\author{G. Cleanthous}
\address{Department of Mathematics and Statistics,
National University of Ireland, Maynooth}
\email{galatia.cleanthous@mu.ie}
\author{A. G. Georgiadis}
\address{School of Computer Science and Statistics,
Trinity College of Dublin}
\email{georgiaa@tcd.ie}
\author{O. V. Lepski}
\address{Institut de Math\'ematique de Marseille, CNRS,\\
Aix-Marseille  Universit\'e}
\email{oleg.lepski@univ-amu.fr}



\begin{abstract}

This is the second part of the research project initiated in 

\noindent\cite{G+N+O-Part1}. We deal with the problem of
the adaptive estimation of the $\bL_2$--norm
of a probability density
on $\bR^d$, $d\geq 1$, from independent observations.
The unknown density
is assumed to be uniformly bounded by unknown constant and to belong to the union of balls  in the isotropic/anisotropic
Nikolskii's spaces. In \cite{G+N+O-Part1} we have proved that the optimally  adaptive estimators do no exist in the considered problem
and provided with several lower bounds for the adaptive risk. In this part we  show that these bounds are tight and present the adaptive estimator
which is obtained by a data-driven selection from a family of kernel-based estimators.
The proposed  estimation procedure as well as the computation of its
 risk are heavily based on new  concentration inequalities for decoupled $U$-statistics of order two established in Section \ref{sec:subsec-U-stat}.
 It is also worth noting that all our results are derived from the unique oracle inequality which may be of  independent interest.

\end{abstract}

\date{May, 2024}

\maketitle

\section{Introduction}

 Suppose that we observe i.i.d. random vectors $X_i\in\bR^d$, $i=1,\ldots, n,$
 with  common probability density $f$.
 We are interested in estimating
the $\bL_2$-norm of $f$,
 $$
 \|f\|^2_2:=\int_{\bR^d}f^2(x)\rd x,
 $$
from observation $X^{(n)}=(X_1,\ldots,X_n)$.
By estimator of $\|f\|_2$ we mean any $X^{(n)}$-measurable map $\widetilde{N}:\bR^n\to \bR$.

\noindent Accuracy of an estimator $\widetilde{N}$
is measured by the quadratic risk
$$
 \cR_n[\widetilde{N}, \|f\|_2]:=\Big(\bE_f \big[\widetilde{N}-\|f\|_2\big]^2\Big)^{1/2},
$$
 where $\bE_f$ denotes  expectation with respect to the probability measure
$\bP_f$ of the observations $X^{(n)}=(X_1,\ldots,X_n)$.
We adopt minimax approach to measuring estimation accuracy.
Let $\mF$ denote the set of all square-integrable probability densities defined on $\bR^d$.
The maximal risk of an estimator $\widetilde{N}$ on the set
$\bF\subset\mF$ is defined by
$$
\cR_n\big[\widetilde{N}, \bF\big]:=\sup_{f\in\bF}\cR_n[\widetilde{N}, \|f\|_2],
$$
and the minimax risk on $\bF$ is then
$
 \cR_n[\bF]:=\inf_{\widetilde{N}} \cR_n[\widetilde{N}; \bF],
$
where $\inf$ is taken over all possible estimators.  An estimator $\widehat{N}_n$ is called
{\em optimal in order} or {\em rate--optimal} if
\[
 \cR_n[\widehat{N}_n; \bF] \asymp \cR_n[\bF],\;\;\;n\to\infty.
\]
The rate at which $\cR_n[\bF]$ converges to zero as $n$ tends to infinity is referred to as
{\em the minimax rate of convergence}.

\paragraph*{Minimax estimation of $\bL_2$-norm.}
 Let $\bF=\bN_{\vec{r}, d}(\vec{\beta}, \vec{L})\cap\bB_\infty(Q)$, where $\bN_{\vec{r}, d}(\vec{\beta}, \vec{L})$ is an anisotropic Nikolskii's class (see Section \ref{sec:adaptive-results-L_2} for formal definition) and
 \[
\bB_\infty(Q):=\big\{f:\bR^d\to\bR: \|f\|_\infty\leq Q\big\},\; Q>0.
\]
The problem of minimax estimation of $\bL_p$-norms of a density, $1<p<\infty$, on the intersection of the Nikolskii's class with a given ball in $\bL_q$-space, $q\in (p,\infty]$ was studied in \cite{GL20a}, \cite{GL20b}. Their results reduced to the case $p=2$, $q=\infty$ considered in the present paper can be summarized as follows.
Following \cite{G+N+O-Part1} introduce the following notation: $\varsigma=(\vec{\beta},\vec{r})$ and for any $s\in[1,\infty]$ set
\begin{equation*}
\tau_{\varsigma}(s):=1-\tfrac{1}{\omega_\varsigma}+\tfrac{1}{\beta s},\quad\; \tfrac{1}{\beta}:=\sum_{j=1}^d\tfrac{1}{\beta_j},\quad\; \tfrac{1}{\omega_\varsigma}:=\sum_{j=1}^d\tfrac{1}{\beta_jr_j}.
\end{equation*}
Define at last
\begin{gather}
\label{eq:rate-exponent-infty}
\mz^*(\varsigma)=\left\{
\begin{array}{clc}
\frac{1}{\tau_\varsigma(1)},\quad&\tau_\varsigma(2)\geq 1,\; \tau_\varsigma(1)>2;
\\*[2mm]
\tfrac{1}{2-\tau_\varsigma(\infty)},\quad&\tau_\varsigma(2)< 1,\;\tau_\varsigma(\infty)<0;
\\*[2mm]
\frac{1}{2},\quad&\tau_\varsigma(2)\geq 1,\;\tau_\varsigma(1)\leq 2;
\\*[2mm]
\frac{1}{2},\quad&\tau_\varsigma(2)<1,\;\tau_\varsigma(\infty)\geq 0.
\end{array}
\right.
\end{gather}
The  following result corresponding to the case  $p=2,q=\infty$ can be deduced from Theorem 1 in \cite{GL20a}.
For  given $\varsigma\in (0,\infty)^d\times[1,\infty]^d$, $\vec{L}\in(0,\infty)^d$ and $Q>0$
\begin{equation}
\label{eq:th_LB}
 \cR_n\big[\bN_{\vec{r},d}\big(\vec{\beta},\vec{L}\big)\cap \bB_\infty(Q)\big]\asymp n^{-\mz^*(\varsigma)}, \quad n\to\infty.
\end{equation}
The rate optimal estimator of $\bL_p$-norm for an arbitrary integer $p\geq 2$ was constructed in \cite{GL20b} under additional assumptions, which in the case $p=2, q=\infty$, are reduced to $\vec{r}\in[1,2]^d\cup[2,\infty]^d$.

\smallskip

\noindent{\bf Minimax adaptive estimation of $\bL_2$-norm.} The problem of the adaptive estimation $\bL_2$-norm of smooth bounded densities can be formulated as follows:
\begin{quote}
is it possible to construct a single estimator $N^*_n$
 which is  simultaneously minimax on each class
 $\bN_{\vec{r},d}\big(\vec{\beta},\vec{L}\big)\cap \bB_\infty(Q)$, i.e. such that
\begin{equation}
\label{eq:question}
\limsup_{n\to\infty} n^{\mz^*(\varsigma)}\cR_n\big[N^*_n;\; \bN_{\vec{r},d}\big(\vec{\beta},\vec{L}\big)\cap \bB_\infty(Q)\big]<\infty,\quad\forall \varsigma,\vec{L}, Q?
\end{equation}
\end{quote}
In the companion paper \cite{G+N+O-Part1}  we proved that the answer on aforementioned question is negative. Also, we direct the reader to \cite{G+N+O-Part1} for background information on adaptive theory and relevant literature on minimax and minimax adaptive estimation of nonlinear functionals.

Here, without providing the exact definitions of terms and the statements of results from \cite{G+N+O-Part1}, we will recall one of our assertions in an informal manner. Set\footnote{Recall that  the symbol $(r)$ is employed to represent the rate.}
$$
\Theta_{(r)}=\big\{\varsigma\in(0,\infty)^d\times[1,\infty]:\; \text{either}\;\tau_\varsigma(2)\geq 1,\; \tau_\varsigma(1)>2\;\text{or}\;\tau_\varsigma(2)< 1,\;\tau_\varsigma(\infty)<0 \big\}.
$$
Note that in view of (\ref{eq:rate-exponent-infty}) the parameter set $\Theta_{(r)}$ corresponds to nonparametric regime of the rate of convergence.

Then, in order for (\ref{eq:question}) to be upheld, it is necessary to replace the normalization family $\{n^{-\mz^*},\varsigma\in\Theta_{(r)}\}$ by an another one which cannot be "faster" than the family given by
$$
\Psi=\big\{\psi_n(\varsigma),\;\varsigma\in\Theta_{(r)}\big\},\quad \psi_n\big(\varsigma\big)=\Big(\tfrac{\sqrt{\ln(n)}}{n}\Big)^{\mz^*(\varsigma)}.
$$
In Section \ref{sec:adaptive-results-L_2} we prove that  $\Psi$ is an admissible  normalization family over the scale of anisotropic functional classes and that it  the adaptive rate of convergence, i.e. satisfies Definition 2 in
\cite{G+N+O-Part1}, over the collection of isotropic functional classes.

\paragraph{Adaptive upper bounds via oracle approach.}
Let us now discuss, using norm estimation as an example, the application of the oracle approach to adaptive estimation.
Let $\mN_\Lambda=\big\{\widehat{N}_\lambda, \lambda\in \Lambda\big\}$ be a finite\footnote{We consider a finite $\Lambda$ so as not to get into a discussion about measurability.} family of estimators of $\bL_2$-norm  built from the observation $X^{(n)}$.

We want to construct a $\Lambda$-valued random element  $\hat{\lambda}$
 completely determined by the observation  $X^{(n)}$ and to prove that for all  $n$ large enough
\begin{equation}
\label{eq:ell-oracle-inequality}
\cR_n\big[\widehat{N}_{\hat{\lambda}}; \|f\|_2\big]\leq \inf_{\lambda\in \Lambda}\mathrm{O}_{n}(f,\lambda)+r_n,\quad\forall
f\in\cF\subseteq \mF.
\end{equation}
Here $\cF$ is either $\mF$ or its "massive" subset and $r_n$ is the remainder term which may depend on $\cF$ and the  family of estimators $\mN_\Lambda$ only.
The quantity
$\cO_n(\cdot,\cdot)$  is known  explicitly and it is tacitly assumed that $r_n\to 0, n\to\infty$ quite rapidly.
Historically,    (\ref{eq:ell-oracle-inequality}) is called \textit{oracle inequality} if (reduced to the considered problem)
 it would be  proved  with
\begin{equation}
\label{eq:ell-oracle-inequality-new}
\mathrm{O}_{n}(f,\lambda)=\mathrm{O}^{\text{opt}}_{n}(f,\lambda):=C\cR_n\big[\widehat{N}_{\lambda}; \|f\|_2\big].
\end{equation}
Here $C$ is a constant  which may depend on $\cF$ and the  family of estimators $\mN_\Lambda$ only.
 The latter means that the ''oracle''
knowing the true parameter $f$ can construct  the estimator $\widehat{N}_{\lambda(f)}$ which provides the minimal risk over the collection $\mN_\Lambda$  for any $f\in\cF$, that is
\[
\lambda(f):\;\; \cR_n\big[\widehat{N}_{\lambda(f)}; \|f\|_2\big]=\min_{\lambda\in\Lambda}\cR_n\big[\widehat{N}_{\lambda}; \|f\|_2\big].
\]
Since $\lambda(f)$ depends on unknown $f$ the estimator $\widehat{N}_{\lambda(f)}$, called oracle estimator, is not an estimator in usual sense and, therefore, cannot be used. The goal is to construct the estimator $\widehat{N}_{\hat{\lambda}}$  which ''mimics''  the oracle one.
It is worth noting that the oracle inequality (\ref{eq:ell-oracle-inequality}) with right hand side given by (\ref{eq:ell-oracle-inequality-new})  is not always available, and this is the reason why we deal with more general definition given in (\ref{eq:ell-oracle-inequality}). In particular if $\mN_\Lambda$ is the family of minimax estimators constructed in \cite{GL20b} then the oracle inequality with (\ref{eq:ell-oracle-inequality-new}) does not hold.

The important remark   is that being proved the inequality (\ref{eq:ell-oracle-inequality})
provides a simple enough criterion to determine whether the selected estimator  $\widehat{N}_{\hat{\lambda}}$ is
adaptive  with respect to the given scale of functional classes $\big\{\bF_\theta\subset\cF,\theta\in\Theta\big\}$ or not.
Define for any $\theta\in\Theta$
$$
\mathfrak{O}_n(\theta)=\inf_{\lambda\in\Lambda}\sup_{f\in\bF_\theta}\mathrm{O}_{n}\big(f,\lambda\big)+r_n
$$
and consider the normalization family   $\big\{\mathfrak{O}_n(\theta),\theta\in\Theta\big\}$.
First, we remark that it is an admissible family. Indeed, we have in view of (\ref{eq:ell-oracle-inequality}), since $r_n$ is independent of $f$
\begin{align*}
\sup_{f\in\bF_\theta}\cR_n\big[\widehat{N}_{\hat{\lambda}}; \|f\|_2\big]&\leq \sup_{f\in\bF_\theta}\Big[\inf_{\lambda\in \Lambda}\mathrm{O}_{n}(f,\lambda)\Big]+r_n,
\\*[2mm]
&\leq \inf_{\lambda\in \Lambda}\Big[\sup_{f\in\bF_\theta}\mathrm{O}_{n}(f,\lambda)\Big]+r_n=\mathfrak{O}_n(\theta),\quad\forall\theta\in\Theta.
\end{align*}
Hence if $\big\{\mathfrak{O}_n(\theta),\theta\in\Theta\big\}$ satisfies the conditions of Theorems 3 or 4 in \cite{G+N+O-Part1} then we can assert that it is the adaptive rate of convergence and $\widehat{N}_{\hat{\lambda}}$ is an adaptive estimator.

\par
\noindent In Section \ref{sec:selection rule and oracle inequality} we propose a data-driven selection rule from the family $\{\widehat{N}_{\vec{h}}\}$ of kernel-based
estimators parameterized by the collection of multi-bandwidths $\vec{h}$ and   establish for it the oracle inequality (\ref{eq:ell-oracle-inequality}) with $\cF=\bB_\infty\big(\alpha_n^{-1}\big)$, where $\alpha_n\to 0, n\to\infty,$ is a given sequence. Then we derive all our adaptive results from this unique oracle inequality.

\section{Selection rule and oracle inequality.}
\label{sec:selection rule and oracle inequality}

Furthermore, without loss of generality  we will assume that $n=2m$ and  denote by $Y_i=X_{m+i}, i=1,\ldots m$. The asymptotics will be studied under $m\to\infty$.

\vskip0.1cm

\noindent {\bf Why we need data-splitting?}\; Our estimation construction as well as the  results concerning its
minimax risk are heavily based on the concentration inequalities for $U$-statistics of order two.
Unfortunately the best known results namely \cite{gine-zinn} and \cite{patricia2003} are not sufficiently sharp for our purposes.
Our results presented in Section \ref{sec:subsec-U-stat} are sharper, but they are proved for decoupled $U$-statistics only
(see Remark \ref{rem:after-th-main-ustat}  for detailed discussion). The latter necessitates  data-splitting.

\vskip0.1cm

\noindent {\bf More general risk.}\; We will establish our main results (except Theorem \ref{th:main-adap-isotrop} proved for quadratic risk)
for more general type of risks. Let $q\in [1,\infty)$ be given real number.  Accuracy of an estimator $\widetilde{N}$
will be measured by the  $q$-risk that is
\vskip0.2cm
\centerline{$
 \cR^{(q)}_{2m}[\widetilde{N}, \|f\|_2]:=\Big(\bE_f \big|\widetilde{N}-\|f\|_2\big|^q\Big)^{1/q}.
$}
\vskip0.2cm
\noindent The maximal over given set of densities $\bF\subset\mF$ $q$-risk is given as usual by
\vskip0.2cm
\centerline{$
 \displaystyle{\cR^{(q)}_{2m}[\widetilde{N}, \bF]:=\sup_{f\in\bF}\cR^{(q)}_{2m}[\widetilde{N}, \|f\|_2].}
$}
\vskip0.2cm

\subsection{Family of estimators.}
 \label{sec:ideas}
 Let $\cK:\bR\to\bR$ be a symmetric bounded function supported on $(-t/2,t/2)$, $t>0$, satisfying $\int_{\bR}\cK(x)\rd x=1$.
Set for any  $x=(x_1,\ldots,x_d)\in\bR^d$
$$
K(x)=\prod_{j=1}^d\cK(x_j),\quad\; T(x)=2K(x)-\int_{\bR^d}K(v)K(x+v)\rd v
$$

\noindent {\bf Family of kernels.}\;
Set\footnote{Here and later $ab$, $a/b$, $a\wedge b$, $a\vee b$ and $a\gtrless$ for any $a,b\in\bR^d$ are understood coordinately.}
for any $\vec{h}=(h_1,\ldots,h_d)\in (0,\infty)^d$
\begin{equation*}
K_{\vec{h}}(x)=V^{-1}_{\vec{h}}K(x/\vec{h}),\quad T_{\vec{h}}(x)=V^{-1}_{\vec{h}}T(x/\vec{h}),\quad V_{\vec{h}}=\prod_{j=1}^dh_j.
\end{equation*}
The following notation related to the kernel $K$ will be used in the sequel:
$
\varpi_T=1\vee\sup_{p\in[1,\infty]}\|T\|_p.
$
Define also for any $x\in\bR^d$ and $\vec{h}\in \bR^d_+$
\begin{equation}
\label{eq:bais+smoother}
 S_{\vec{h}}(x):= \int_{\bR^d} K_{\vec{h}}(x-y) f(y) \rd y,\;\;\;B_{\vec{h}}(x):= S_{\vec{h}}(x)-f(x).
\end{equation}
Obviously, $B_{\vec{h}}(x)$ is the bias of the kernel density estimator of $f(x)$
associated with kernel~$K$ and multi-bandwidth $\vec{h}$.

\vskip0.1cm

\noindent {\bf Basic observation.}\; The construction of our estimator's collection is based on
a simple observation formulated below
as Lemma \ref{lem:representation}. Its elementary proof is postponed  to Section \ref{sec:Proof of technical lemmas}.

  \begin{lemma}
 \label{lem:representation}
For any    $f\in\mF$
and  $h\in (0, \infty)^d$ one has
  \begin{equation}
   \|f\|_2^2 =  \int_{\bR^d}\int_{\bR^d} T_{\vec{h}}(x-y)f(x)f(y)\rd x\rd y
 + \int_{\bR^d}  B^2_{\vec{h}}(x)\rd x.
\label{eq:N-representation-new}
 \end{equation}
 \end{lemma}

\noindent The natural unbiased estimator for $\int_{\bR^d}\int_{\bR^d} T_{\vec{h}}(x-y)f(x)f(y)\rd x\rd y $ is the decoupled  U--statistics
\begin{equation*}
 \widehat{N}_{\vec{h}}=\tfrac{1}{m^2}\sum_{i,j=1}^mT_{\vec{h}}(Y_i-X_j).
\end{equation*}

\noindent {\bf Set of bandwidths.}\; Let $m\in\bN^*$ be given and set
 $$
 \mH^d=\big\{\vec{h}\in (0,1]^d:\; m^2V_{\vec{h}}\geq \ln(m)\big\}.
 $$
 Introduce also $\cZ_m=\{e^{-1},e^{-2},\ldots,e^{-\lfloor2 \ln(m)\rfloor+1},e^{-2\lfloor\ln(m)\rfloor}\}$ and let
 \begin{equation}
 \label{eq:set-of-band}
 \cH^d_m=\mH^d\cap\cZ_m^d.
 \end{equation}

 \noindent {\bf Set of densities.}\; All our considerations below will be done under condition that $f\in\bB_\infty\big(\alpha_m^{-1}\big)$, where $\alpha_m \downarrow 0, m\to\infty$, be an arbitrary but \textit{a priory} chosen sequence, verifying $\alpha_m\ln(m)\geq 1$ and $\alpha_m\leq 1/3$ for any $m\in\bN^*$.  Also, we will assume that $m\geq 9$.

 \subsection{Selection rule.}
 \label{sec:selection-rule}
In this section we propose a data-driven choice from the estimator's collection $\big\{\widehat{N}_{\vec{h}}, \vec{h}\in\cH^d_m\big\}$.
 This estimation procedure is inspired by the general selection method introduced in \cite{lepski2018} in the context of an abstract statistical experiment.
 Let  $q\in [1,\infty)$ be fixed and introduce the following sets of bandwidths
\begin{align*}
H^d_m=\big\{\vec{h}\in\cH^d_m:\; mV_{\vec{h}}< \ln(m)\big\}&,\quad\; \bar{H}^d_m=\cH^d_m\setminus H^d_m;
\\*[2mm]
\mH^d_*=\big\{\vec{h}\in\mH^d:\; mV_{\vec{h}}\leq \alpha^2_m\big\}&,\quad\;\bar{\mH}^d_*=\mH^d\setminus\mH^d_*.
\end{align*}
Set for any $\vec{h}\in\cH^d_m$ (recall that $\cH^d_m\subset\mH^d$)
\begin{equation}
\label{random-J}
\widehat{J}_{\vec{h}}=\frac{1}{m^2}\sum_{i,j=1}^m|T_{\vec{h}}(Y_i-X_j)|,\quad\; \widehat{W}_{\vec{h}}=\widehat{\cW}_{\vec{h}}\mathrm{1}_{\bar{H}^d_m}\big(\vec{h}\big)+\alpha_m^{-1}\mathrm{1}_{H^d_m}\big(\vec{h}\big).
\end{equation}
Here the random variable $\widehat{\cW}_{\vec{h}}$ is given by
\begin{align*}
\widehat{\cW}_{\vec{h}}&=\tfrac{1}{mV_{\vec{h}}}\Big[256s\ln(m)\bigvee\sup_{\vec{k}\in\bZ^d}
\sum_{j=1}^m\mathrm{1}_{\Pi^*_{\vec{k}}}\big(\tfrac{X_j}{\vec{h}}\big)\Big],
\end{align*}
where $s=12q+2$ and  for any $\vec{k}=(k_1,,\ldots,k_d)\in\bZ^d$
$$
\Pi^*_{\vec{k}}=\big((k_1-1)t,(k_1+2)t\big]\times\cdots\times\big((k_d-1)t,(k_d+2)t\big].
$$
Set $\Upsilon_{\vec{h}}(m)=\big[10\ln(m)\big]\mathrm{1}_{\bar{\mH}^d_*}\big(\vec{h}\big)
+\Big[\tfrac{17\ln(m)}{|\ln(mV_{\vec{h}})|}\Big]\mathrm{1}_{\mH^d_*}\big(\vec{h}\big)$ and define finally
\begin{align*}
\cU_{\vec{h}}&=\Big[\tfrac{16s \varpi_T\widehat{W}_{\vec{h}}\widehat{J}_{\vec{h}} \ln(m)}{m}\Big]^{\frac{1}{2}}+
\Big[\tfrac{16s\varpi_T \widehat{J}_{\vec{h}}\ln(m) }{m^2V_{\vec{h}}}\Big]^{\frac{1}{2}}+\tfrac{147s\varpi_T \widehat{W}_{\vec{h}}\ln(m)}{m}
+\tfrac{s^2\varpi_T\Upsilon_{\vec{h}}(m) \ln(m)}{m^2V_{\vec{h}}}.
\end{align*}
\noindent {\bf Selection rule.}\; Set for any $\vec{h}\in\cH^d_m$
\begin{equation}
\label{eq:sel-rule1}
\widehat{\cR}_{\vec{h}}=
\max_{\vec{w}\in\cH^d_m}\Big(\big|\widehat{N}_{\vec{h}\vee\vec{w}}-\widehat{N}_{\vec{w}}\big|-18\widehat{\cU}_{\vec{w}}\Big)_+,
\quad \widehat{\cU}_{\vec{w}}=\max_{\vec{\mathbf{w}}\in\cH_m^d}\cU_{\vec{w}\vee\vec{\mathbf{w}}},
\end{equation}
and define
\begin{equation}
\label{eq:sel-rule2}
\vec{h}^*=\arg\min_{\vec{h}\in\cH^d_m}\big(\widehat{\cR}_{\vec{h}}+18\widehat{\cU}_{\vec{h}}\big).
\end{equation}
We remark that $\vec{h}^*\in\cH^d_m$ since $\cH^d_m$ is a finite set and, by the same reason, $\vec{h}^*$ is $X^{(n)}$-measurable.
Our final estimator is $|\widehat{N}_{\vec{h}^*}|^{\frac{1}{2}}$.

\subsection{Oracle inequality.}
\noindent Introduce for any $\vec{h}\in\cH^d_m$
\begin{align*}
\mU_{\vec{h}}(f)&=\Big[\tfrac{\kappa_{\vec{h}}(f) \|f\|_2^2 \ln(m)}{m}\Big]^{\frac{1}{2}}+\Big[\tfrac{\|f\|_2^2\ln(m) }{m^2V_{\vec{h}}}\Big]^{\frac{1}{2}}
+\tfrac{\kappa_{\vec{h}}(f) \ln(m)}{m}+
\tfrac{\Upsilon_{\vec{h}}(m)\ln(m)}{m^2V_{\vec{h}}};
\\*[2mm]
\kappa_{\vec{h}}(f)&=\big(\tfrac{256s\ln(m)}{mV_{h}}\vee\|f\|_\infty\big)\mathrm{1}_{\bar{H}^d_m}\big(\vec{h}\big)
+\alpha_m^{-1}\mathrm{1}_{H^d_m}\big(\vec{h}\big).
\end{align*}
Recall that $B_{\vec{h}}(\cdot)$ is defined in (\ref{eq:bais+smoother}) and introduce the quantity that we call the {\it oracle risk}.
\begin{equation*}
\cO_f= \min_{\vec{h}\in\cH^d_m}\Big\{\big\|B_{\vec{h}}\big\|_2^2+\max_{\vec{w}\in \cH_m^d}\Big|\big\|B_{\vec{h}\vee\vec{w}}\big\|^2_2-\big\|B_{\vec{w}}\big\|^2_2\Big|
 +\mU_{\vec{h}}(f)\Big\}
  +[\ln(m)]^{\frac{d}{2q}}m^{-2}.
\end{equation*}

\begin{theorem}
\label{th:oracle-inequality}
Let $q\in[1,\infty)$ be given. There exists $\mathbf{c}>0$ depending on $\cK$ and $q$ only such that for any $m\geq 9$ and $f\in\bB_\infty\big(\alpha_m^{-1}\big)$ one has
$$
\cR^{(q)}_{2m}\big[|\widehat{N}_{\vec{h}^*}|^{\frac{1}{2}},f\big]\leq \mathbf{c}\Big[\tfrac{\cO_f}{\|f\|_2}\bigwedge \sqrt{\cO_f}\Big].
$$
\end{theorem}

\noindent In order to make possible the maximization of the oracle risk over given collection of isotropic/anisotropic Nikolskii classes,
we will need  to bound effectively the quantity
$$
\big\|B_{\vec{h}}\big\|_2^2+\max_{\vec{w}\in \cH_m^d}\Big|\big\|B_{\vec{h}\vee\vec{w}}\big\|^2_2-\big\|B_{\vec{w}}\big\|^2_2\Big|.
$$
For any $a>0$    introduce
\begin{equation*}
\label{eq:def-small-b_h}
b_{a,j}(x)=
\int_{\bR}\cK(u)f\big(x+ua\mathbf{e}_j\big)\rd u-f(x),\quad j=1,\ldots,d,
\end{equation*}
where, recall, $\{\mathbf{e}_1, \ldots, \mathbf{e}_d\}$ denotes the canonical basis in $\bR^d$. Set
$$
\cO^*_f=\min_{\vec{h}\in\cH^d_m}\Big\{\mb^2_{\vec{h}}(f)+\|f\|_2\mb_{\vec{h}}(f) +\mU_{\vec{h}}(f)\Big\}
  +[\ln(m)]^{\frac{d}{2q}}m^{-2},\quad \mb_{\vec{h}}(f)=\max_{l=1,\ldots,d}\sup_{a\leq h_l}\big\|b_{a,l}\big\|_2.
$$

\begin{corollary}
\label{cor2:th-oracle}
Let $q\in[1,\infty)$ be given. There exists $\mathbf{c_1}>0$ depending on $\cK$ and $q$ only such that for any $m\geq 9$ and $f\in\bB_\infty\big(\alpha_m^{-1}\big)$ one has
$$
\cR^{(q)}_{2m}\big[|\widehat{N}_{\vec{h}^*}|^{\frac{1}{2}},f\big]\leq \mathbf{c_1}\Big[\tfrac{\cO^{*}_f}{\|f\|_2}\bigwedge \sqrt{\cO^{*}_f}\Big].
$$
\end{corollary}

\section{Adaptive estimation of $\bL_2$-norm of smooth bounded densities.}
\label{sec:adaptive-results-L_2}

\noindent We begin this section with the definition of the Nikolskii functional class.

\vskip0.1cm

\noindent {\bf Anisotropic Nikolskii class.}\; Let $(\Be_1,\ldots,\Be_d)$ denote the canonical basis of $\bR^d$.
 For function $G:\bR^d\to \bR$ and
real number $u\in \bR$
{\em the first order difference operator}
with step size $u$ in direction of variable
$x_j$ is defined by
$\Delta_{u,j}G (x):=G(x+u\Be_j)-G(x),\;j=1,\ldots,d$.
By induction,
the $k$-th order difference operator with step size $u$ in direction of  $x_j$ is
\begin{equation}
\label{eq:Delta}
 \Delta_{u,j}^kG(x)= \Delta_{u,j} \Delta_{u,j}^{k-1} G(x) = \sum_{l=1}^k (-1)^{l+k}\tbinom{k}{l}\Delta_{ul,j}G(x).
\end{equation}
\begin{definition}
\label{def:nikolskii}
For given  vectors $\vec{\beta}=(\beta_1,\ldots,\beta_d)\in (0,\infty)^d$, $\vec{r}=(r_1,$ $\ldots,r_d)\in [1,\infty]^d$,
 and $\vec{L}=(L_1,\ldots, L_d)\in (0,\infty)^d$ we
say that the function $G:\bR^d\to \bR$ belongs to  anisotropic
Nikolskii's class $\bN_{\vec{r},d}\big(\vec{\beta},\vec{L}\big)$ if
 $\|G\|_{r_j}\leq L_{j}$ for all $j=1,\ldots,d$ and
 there exist natural numbers  $k_j>\beta_j$ such that
\[
 \big\|\Delta_{u,j}^{k_j} G\big\|_{r_j} \leq L_j |u|^{\beta_j},\;\;\;\;
\forall u\in \bR,\;\;\;\forall j=1,\ldots, d.
\]
\end{definition}

\vskip0.1cm

\noindent {\bf Special construction of the kernel $\cK$.}\; We will use the following specific kernel $\cK$ [see, e.g., \cite{lepski-kerk}] in the definition of the estimator's family $\big\{\widehat{N}_{\vec{h}}, \vec{h}\in\cH^d_m\big\}$.

\smallskip

 Let  $b>1$  be an arbitrary large but  \textit{a priory} chosen integer number. Furthermore we will assume that $\vec{\beta}\in (0,b)^d$
 and let $D:[-1/2, 1/2]\to \bR$ be a symmetric, locally integrable bounded function satisfying  $\int_\bR D=1$.
Define
\begin{equation}
\label{eq:cK-function}
 \cK(y)=-\sum_{k=1}^b (-1)^{k}\tbinom{b}{k} \tfrac{1}{k}D\big(\tfrac{y}{k}\big).
\end{equation}
We note that $\cK$ is compactly supported on $[-b/2,b/2]$ and we can easily check that $\int_{\bR}\cK=1$.


\subsection{Main results. Anisotropic case.} Let $\widehat{N}_{\vec{h}^*}$ denote the estimator obtained by the selection rule (\ref{eq:sel-rule1})--(\ref{eq:sel-rule2}) from the family
of  $\big\{\widehat{N}_{\vec{h}}, \vec{h}\in\cH^d_m\big\}$, with function $T$ being constructed with the help of the function $\cK$ given
in (\ref{eq:cK-function}).
Set for any $\varsigma\in (0,b)^d\times[1,\infty]$
 \begin{equation}
 \label{eq:adaptive-rate-anis}
 \psi_m(\varsigma)=(\mu_m)^{\mz^*(\varsigma)},\quad \mu_m=\tfrac{\sqrt{\ln(m)}}{m}\;\mathrm{1}\big(\mz^*(\varsigma)< 1/2\big)+\tfrac{\ln(m)}{m}\mathrm{1}\big(\mz^*(\varsigma)=1/2\big).
  \end{equation}

 \begin{theorem}
 \label{th:adaptive-anisotr}
Let $q\in[1,\infty)$ and $b>1$ be fixed. For any $\varsigma\in(0,b)^d\times\big\{[1,2]^d\cup[2,\infty]^d\big\}$, $\vec{L}\in(0,\infty)^d$ and $Q>0$
$$
\limsup_{m\to\infty}\psi^{-1}_m(\varsigma)\;\cR^{(q)}_{2m}\Big[|\widehat{N}_{\vec{h}^*}|^{\frac{1}{2}}, \bN_{\vec{r},d}(\vec{\beta},\vec{L})\cap\bB_\infty(Q)\Big]<\infty.
$$
 \end{theorem}

\noindent  We were able to construct a completely data-driven estimation procedure whose accuracy looses  only "logarithmic" factor with respect to
the minimax rate of convergence. The latter result is proved under additional assumption $\vec{r}\in [1,2]^d\cup[2,\infty]^d$. It is worth noting that the minimax rate is known only if this assumption holds.
Let us make now some remarks in the case $q=2$. The following set of parameters is referred to as the non-parametric regime.
$$
\mS_{\text{nparam}}=\big\{\varsigma\in(0,b)^d\times\big\{[1,2]^d\cup[2,\infty]^d\big\}:\; \mz^*(\varsigma)< 1/2\big\}.
$$
Combining the results of Theorem \ref{th:adaptive-anisotr} and  Proposition 1 in \cite{G+N+O-Part1} we can assert that
the normalization family $\{\psi_m(\varsigma),\varsigma\in \mS_{\text{nparam}}\}$ is an adaptive rate of convergence on
$$
\mS^{\prime}_{\text{nparam}}=\big\{\varsigma\in\mS_{\text{nparam}}:\; \mz^*(\varsigma)=\tfrac{1}{\tau_\varsigma(1)}\big\}.
$$
Combining the results of Theorem \ref{th:adaptive-anisotr} and  Proposition 2 in \cite{G+N+O-Part1} we can assert that
the normalization family $\{\psi_m(\varsigma),\varsigma\in \mS_{\text{nparam}}\}$ is an adaptive rate of convergence on some part of the set
$$
\mS^{\prime\prime}_{\text{nparam}}=\big\{\varsigma\in\mS_{\text{nparam}}:\; \mz^*(\varsigma)=\tfrac{1}{2-\tau_\varsigma(\infty)}\big\}.
$$

\paragraph{Open problem 1.} \textsf{To prove or to disprove that} $\{\psi_m(\varsigma),\varsigma\in \mS_{\text{nparam}}\}$ \textsf{is an adaptive rate of convergence across entire nonparametric regime.}

\vskip0.1cm

\noindent The following set of parameters is referred to as the parametric regime.
$$
\mS_{\text{param}}=\big\{\varsigma\in(0,b)^d\times\big\{[1,2]^d\cup[2,\infty]^d\big\}:\; \mz^*(\varsigma)=1/2\big\}.
$$
As it follows from Theorem \ref{th:adaptive-anisotr} the accuracy provided by the estimator $|\widehat{N}_{\vec{h}^*}|^{\frac{1}{2}}$ is $\sqrt{\frac{\ln(m)}{m}}$ for any
$\varsigma\in\mS_{\text{param}}$
while minimax rate of convergence is $m^{-1/2}$.

\paragraph{Open problem 2.} \textsf{To prove or to disprove that there is no optimally adaptive estimators when} $\varsigma$ \textsf{run} $\mS_{\text{param}}$  \textsf{and} $\sqrt{\ln(m)}$ \textsf{is unavoidable price to pay for adaptation.}

\vskip0.1cm

\noindent In the next section we will show that optimally adaptive estimators over collection of \textit{isotropic} Nikolskii classes, with parameters belonging
to $\mS_{\text{param}}$, do exist.

\subsection{Main results. Isotropic case.}
\label{sec:isotropic-case}

In this section we study the adaptive estimation over the scale of \textit{isotropic} Nikolskii classes of bounded functions. From now on
$$
\beta_{l}=\blb,\; r_l=\bbr,\; L_l=L,\quad\forall l=1,\ldots, d.
$$
Since $\vec{\beta}=(\blb,\ldots,\blb)$, $\vec{r}=(\bbr,\ldots,\bbr)$ the isotropic Nikolskii class will be denoted by $\bN_{\bbr,d}\big(\blb,L\big)$. Moreover we obviously have for any $s\in[1,\infty]$
$$
\tau_{\varsigma}(s)=1-\frac{d}{\blb \bbr}+\frac{d}{\blb s}, \quad \frac{1}{\omega_\varsigma}=\frac{d}{\blb \bbr}, \quad \frac{1}{\beta}=\frac{d}{\blb}.
$$
The adaptation will be studied over the collection of functional classes
$$
\bN_\vartheta:=\bN_{\bbr,d}\big(\blb,L\big)\cap \bB_\infty(Q), \quad \vartheta=(\blb,\bbr,L,Q).
$$
In this section we adopt the following notations.
For any $b\in (0,\infty]$ set $\Theta^b_{(r)}=(0,b)\times[1,\infty]$, $\Theta_{(c)}=(0,\infty)^2$ and let
$\vartheta=(\theta,\bar{\theta})$, where   $\theta=(\blb,\bbr)\in\Theta^b_{(r)}$ and $\bar{\theta}=(L,Q)\in\Theta_{(c)}$.
Thus, $\vartheta\in\Theta=:\Theta^b_{(r)}\times\Theta_{(c)}$.


Note finally that the exponent figurant  in the minimax rate of convergence  given in (\ref{eq:rate-exponent-infty}), becomes in the "isotropic case"
\begin{gather}
\label{eq:rate-exponent-infty-isotropic}
\mz^*(\theta)=\left\{
\begin{array}{cccc}
\frac{\blb\bbr}{\blb\bbr+d(\bbr-1)},\quad& \bbr\geq 2,\; \blb\bbr< d(\bbr-1);
\\*[2mm]
\frac{\blb\bbr}{\blb\bbr+d},\quad&\bbr< 2,\; \blb\bbr< d;
\\*[2mm]
\frac{1}{2},\quad&\bbr\geq 2,\; \blb\bbr\geq d(\bbr-1);
\\*[2mm]
\frac{1}{2},\quad&\bbr< 2,\; \blb\bbr\geq d,
\end{array}
\right.
\end{gather}
and introduce the following normalization family.
$$
\Psi^b=\Big\{\psi^*_n(\theta):=\big(\mu_m\big)^{\mz^*(\theta)},\;\;\theta\in\Theta^b_{(r)}\Big\},
$$
where $\mu_n$ is given in (\ref{eq:adaptive-rate-anis}) with $\mz^*(\varsigma)$ replaced by $\mz^*(\theta)$.
Define also the following sets of parameters.
\begin{align*}
\Theta^b_{(r),\text{nparam}}&=\big\{\theta\in\Theta^b_{(r)}:\;  \bbr\geq 2,\; \blb\bbr< d(\bbr-1)\big\}\cup\big\{\theta\in\Theta^b_{(r)}:\; \bbr< 2,\; \blb\bbr< d\big\};
\\*[2mm]
\Theta^b_{(r),\text{param}}&=\big\{\theta\in\Theta^b_{(r)}:\;  \bbr\geq 2,\; \blb\bbr\geq  d(\bbr-1)\big\}\cup\big\{\theta\in\Theta^b_{(r)}:\; \bbr< 2,\; \blb\bbr\geq d\big\}.
\end{align*}

\noindent We note that  $\Theta^b_{(r),\text{nparam}}$  and $\Theta^b_{(r),\text{param}}$ are the set of parameters corresponding to nonparametric regime of minimax rate of convergence and parametric one respectively.
 The following result is direct consequence of Theorem \ref{th:adaptive-anisotr} and Theorem 5 from \cite{G+N+O-Part1}.

\begin{theorem}
\label{th:adaptive-isotr-with-log}
For any fixed $b\in (0,\infty]$ the estimator $|\widehat{N}_{\vec{h}^*}|^{\frac{1}{2}}$ is an adaptive estimator and   the normalization family $\Psi^b$ it is the adaptive rate of convergence over scale of functional classes
$\left\{\bN_\vartheta, \vartheta\in \Theta^b_{(r)}\times \Theta_{(c)}\right\}$. The price to pay for adaptation is given by $\sqrt[4]{\ln(m)}$.
\end{theorem}

\begin{remark}
\label{rem:after-th:adaptive-isotr-with-log}
Thus, the estimator $|\widehat{N}_{\vec{h}^*}|^{\frac{1}{2}}$ is optimal across the whole  set of parameters  corresponding to the nonparametric regime. Moreover, the accuracy of this estimator is inferior to the optimal one (minimax rate) by $\ln(m)$-factor across the entire  set corresponding to the parametric regime.
The latter follows from Theorem \ref{th:adaptive-anisotr} since isotropic classes are a particular case of anisotropic ones.
\end{remark}

The first issue that naturally emerges within this context is as follows:: \textit{is whether the loss in estimation accuracy unavoidable on $\Theta^b_{(r),\text{param}}$?} We will see that the answer is negative and surprisingly, the estimator which is  \textit{optimaly adaptive} across the entire  set corresponding to the parametric regime is $|\widehat{N}_{\vec{\mathrm{h}}_m}|^{\frac{1}{2}}$, where $\vec{\mathrm{h}}_m=(m^{-\frac{1}{d}},\ldots,m^{-\frac{1}{d}})$.

From now on for the reader convenience we will consider the quadratic risk, that is $q=2$. However all results below are true for an arbitrary $q\geq 1$.

\begin{theorem}
\label{th:optim-adap-par-zone}
For any $\vartheta\in \Theta^\infty_{(r),\text{param}}\times \Theta_{(c)}$ one has
$$
\limsup_{m\to\infty}\sqrt{m}\;\cR_{2m}\Big[|\widehat{N}_{\vec{\mathrm{h}}_m}|^{\frac{1}{2}}, \bN_\vartheta\Big]<\infty.
$$
\end{theorem}
\noindent A brief and straightforward proof of the theorem will be presented in the Section \ref{sec:proof-th:optim-adap-par-zone}.

\vskip0.1cm

As we see there are two estimators: the first one behaves  optimally across the entire  set of parameters  corresponding nonparametric regime but it is suboptimal on the set corresponding to the parametric regime. The second one is optimally-adaptive on the set corresponding parametric regime but it has very bad behavior across the set corresponding to nonparametric regime.

The problem which we address now is formulated as follows. We seek an estimation procedure  which would behave optimally on whole set of nuisance parameters $\Theta^b_{(r)}\times\Theta_{(c)}$.
Our construction is based on the data-driven choice between  estimators $|\widehat{N}_{\vec{h}^*}|^{\frac{1}{2}}$ and $|\widehat{N}_{\vec{\mathrm{h}}_m}|^{\frac{1}{2}}$ via Lepski method.

\par

\noindent Consider the following selection rule
\begin{equation}
\label{eq:sel-rule-isotr}
\vec{\mathbf{h}}=\left\{
\begin{array}{ll}
\vec{\mathrm{h}}_m,\quad & \big||\widehat{N}_{\vec{h}^*}|^{\frac{1}{2}}-|\widehat{N}_{\vec{\mathrm{h}}_m}|^{\frac{1}{2}}\big|\leq 2\ln(m)m^{-\frac{1}{2}};
\\*[3mm]
\vec{h}^*,\quad & \big||\widehat{N}_{\vec{h}^*}|^{\frac{1}{2}}-|\widehat{N}_{\vec{\mathrm{h}}_m}|^{\frac{1}{2}}\big|> 2\ln(m)m^{-\frac{1}{2}}.
\end{array}
\right.
\end{equation}
Our final estimator is $|\widehat{N}_{\vec{\mathbf{h}}}|^{\frac{1}{2}}$. Introduce the normalization family $\boldsymbol{\Psi^b}$ as follows.
\begin{equation}
 \label{eq:adaptive-rate-iso}
 \boldsymbol{\psi}_m(\theta)=(\blg_m)^{\mz^*(\theta)},\quad \blg_m=\tfrac{\sqrt{\ln(m)}}{m}\;\mathrm{1}\big(\mz^*(\theta)< 1/2\big)+m^{-1}\mathrm{1}\big(\mz^*(\theta)=1/2\big).
  \end{equation}

\begin{theorem}
\label{th:main-adap-isotrop}
Let $b>1$ be fixed. Then, for any $\vartheta\in \Theta^b_{(r)}\times \Theta_{(c)}$ one has
$$
\limsup_{m\to\infty}\boldsymbol{\psi}^{-1}_m(\theta)\;\cR_{2m}\Big[|\widehat{N}_{\vec{\mathbf{h}}}|^{\frac{1}{2}}, \bN_\vartheta\Big]<\infty.
$$
\end{theorem}

\noindent We conclude that the problem of adaptive estimation of $\bL_2$-norm of a probability density over a collection of isotropic Nikolskii classes of bounded functions is completely solved.

\section{Upper functions for decoupled  $U$-statistics.}
 \label{sec:subsec-U-stat}
Basic probabilistic  tool for proving the main results in the paper  is  exponential inequalities  for $U$-statistics of order two. As we already mention the existing results are not sufficiently sharp for our purposes. In this section we present several results in this direction.  The proposed technique of deriving  such type inequalities for decoupled $U$-statistics may be of independent interest and can be applied in other statistical problems.  The proofs of all theorems below are postponed to Section \ref{sec:proof-U-stat}.

\vskip0.1cm

Let   $X^{(m)}=(X_1\ldots, X_m)\in\bR^{dm}$, $m,d\geq 1$,  be  i.i.d. random vectors
 having common density $f$ and  let $Y^{(m)}=(Y_1\ldots, Y_m)\in\bR^{dm}$ be an independent copy of $X^{(m)}$. The probability law of  $(X^{(m)},Y^{(m)})$ will be  denoted by  $\bP_{f}$ and  $\bE_{f}$ will be used for the mathematical expectation  with respect to $\bP_{f}$.

\paragraph{Kernels.} Let  $U:\bR^d\to\bR$ be \textit{symmetric} bounded  function compactly supported on a given cube $(-t,t)^d$, $t>0$. For any  $\vec{h}=(h_1,\ldots,h_d)\in (0,1]^d$ and $x\in\bR^d$ define
$$
U_{\vec{h}}(x)=V^{-1}_{\vec{h}}U\Big(\tfrac{x_1}{h_1},\ldots,\tfrac{x_d}{h_d}\Big),\quad\; V_{\vec{h}}=\prod_{j=1}^dh_j.
$$
Also let
$
\varpi_U=1\vee\sup_{p\in[1,\infty]}\|U\|_p.
$

\subsection{Non-random upper function for decoupled  $U$-statistics.}
 \label{sec:subsec-U-stat-nonrandom}

Our goal is to find non-random upper function for the collection
$\big\{\big|\widehat{U}_{\vec{h}}-\bE_f\big[\widehat{U}_{\vec{h}}\big]\big|,\;\vec{h}\in\cH_{m}^d\big\}$, where
\begin{equation*}
\widehat{U}_{\vec{h}}=\tfrac{1}{m^2}\sum_{i,j=1}^mU_{\vec{h}}\big(Y_i-X_j\big).
\end{equation*}
We will need the following notations. For any $\vec{h}\in\mH^d$ set
\begin{align}
\label{eq:def-J-and-W}
J_{\vec{h}}(f)&=\int_{\bR^d}\int_{\bR^d}|U_{\vec{h}}(y-z)|f(y)f(z)\rd y\rd z;
\\*[2mm]
W_{\vec{h}}(f)&=\sup_{x\in\bR^d}\int_{\bR^d}|U_{\vec{h}}(y-x)|f(y)\rd y.
\end{align}
\noindent Set also $s=12q+2$, $q\geq 1$, and recall that
$$
\Upsilon_{\vec{h}}(m)=\big[10\ln(m)\big]\mathrm{1}_{\bar{\mH}^d_*}\big(\vec{h}\big)
+\Big[\tfrac{17\ln(m)}{|\ln(mV_{\vec{h}})|}\Big]\mathrm{1}_{\mH^d_*}\big(\vec{h}\big).
$$
Introduce the following quantities.
\begin{align}
\label{eq:def-Lambda(u,Q)}
U_{\vec{h}}(f)&=\Big[\tfrac{16s W_{\vec{h}}(f) J_{\vec{h}}(f)\ln(m)}{m}\Big]^{\frac{1}{2}}+\Big[\tfrac{16s\varpi_UJ_{\vec{h}}(f)\ln(m) }{m^2V_{\vec{h}}}\Big]^{\frac{1}{2}}+\tfrac{23sW_{\vec{h}}(f)\ln(m)}{m}
\nonumber\\*[2mm]
&\quad + \tfrac{\Upsilon_{\vec{h}}(m)s^2\varpi_U \ln(m)}{m^2V_{\vec{h}}}.
\nonumber\\*[2mm]
\Lambda_q(U)&=2^{q+1}\varpi_U^q\Big(1 +2^{2q-1}\big[1+2^q\Gamma(q+1)\big]
\nonumber\\*[2mm]
&\quad+2^{2q-2}\Gamma(q+1)\Big[2^{-1}\big(3\sqrt{12q+2}\big)^{q}
+\big(28(6q+1)\big)^{q}\Big]\Big).
\end{align}

\begin{theorem}
\label{th:non-random-upper-function}

For any $m\geq 9$, $q\geq 1$ and   $f\in\bB_\infty\big(\alpha_m^{-1}\big)$    one has
\begin{align*}
\bE_f\bigg\{\sup_{\vec{h}\in\cH^d_m}\Big(\big|\widehat{U}_{\vec{h}}-\bE_{f}[\widehat{U}_{\vec{h}}]\big|-U_{\vec{h}}(f)\Big)_+^q\bigg\}
\leq \Lambda^*_q(U)[2\ln(m)]^dm^{-4q},
\end{align*}
where
$\Lambda^*_q(U)=2^{q-1}\Lambda_q(U)+\Gamma(q+1)\big[8\varpi_{U}\big]^q$.

\end{theorem}

We will say that $\{U_{\vec{h}}(f), \vec{h}\in\cH^d_m\}$ is an {\it upper function} for the corresponding collection of decoupled and centered $U$-statistics.

\begin{remark}
\label{rem:after-th-main-ustat}
As we have already seen there are two different regimes of  adaptive rate: parametric and nonparametric ones. The nonparametric regime corresponds to the case
$mV_{\vec{h}}\to 0$ $m\to \infty$. More precisely $|\ln(mV_{\vec{h}})|\asymp a\ln(m), a>0,$ and we have in this regime when $m\to\infty$
$$
\Upsilon_{\vec{h}}(m)\asymp \text{const}\quad\Rightarrow\quad \tfrac{\Upsilon_{\vec{h}}(m)s^2\varpi_U \ln(m)}{m^2V_{\vec{h}}}\asymp
\tfrac{\ln(m)}{m^2V_{\vec{h}}}.
$$
This is much sharper than the asymptotics provided by the best known results on the concentration inequalities for $U$-statistics,
namely \cite{gine-zinn} and \cite{patricia2003},  in  which the corresponding term  is
proportional to $\tfrac{\ln^2(m)}{m^2V_{\vec{h}}}$.

\end{remark}

\paragraph{Some consequences.} The following simple consequences of Theorem \ref{th:non-random-upper-function} will be steadily exploited in the sequel.
Their   proofs are given in Section \ref{sec:Proof-of-cor}.

We will need some  notations.
Recall that $(-t,t)^{d}$ is the support of the function $U$ and  for any $\vec{k}=(k_1,,\ldots,k_d)\in\bZ^d$
$$
\Pi^*_{\vec{k}}=\otimes_{r=1}^d\big((k_r-1)t,(k_r+2)t\big].
$$
Define
$$
\cW_{\vec{h}}(f)=V_{\vec{h}}^{-1}\bigg[\tfrac{256s\ln(m)}{m}\bigvee\sup_{\vec{k}\in\bZ^d}\int_{\bR^d}\mathrm{1}_{\Pi^*_{\vec{k}}}\big(\tfrac{y}{\vec{h}}\big)f(y)\rd y\bigg]
$$
and
let $W^*_{\vec{h}}(f)=\cW_{\vec{h}}(f)\mathrm{1}_{\bar{H}^d_m}\big(\vec{h}\big)+\alpha_m^{-1}\mathrm{1}_{H^d_m}\big(\vec{h}\big)$. Set
\begin{align}
\label{eq:definition-of-U}
U^*_{\vec{h}}(f)&=\Big[\tfrac{16s \varpi_U W^*_{\vec{h}}(f) J_{\vec{h}}(f)\ln(m)}{m}\Big]^{\frac{1}{2}}+\Big[\tfrac{16s\varpi_UJ_{\vec{h}}(f)\ln(m) }{m^2V_{\vec{h}}}\Big]^{\frac{1}{2}}+\tfrac{147s\varpi_UW^*_{\vec{h}}(f)\ln(m)}{m}
\nonumber\\*[2mm]
&\quad + \tfrac{\Upsilon_{\vec{h}}(m)s^2\varpi_U \ln(m)}{m^2V_{\vec{h}}}.
\end{align}

\begin{corollary}
\label{cor1:th2}
For any $m\geq 9$, $q\geq 1$ and   $f\in\bB_\infty\big(\alpha_m^{-1}\big)$
\begin{equation*}
\bE_f\bigg\{\sup_{\vec{h}\in\cH^d_m}\Big(\big|\widehat{U}_{\vec{h}}-\bE_{f}[\widehat{U}_{\vec{h}}]\big|-U^*_{\vec{h}}(f)\Big)_+^q\bigg\}
\leq \Lambda^*_q(U)[2\ln(m)]^dm^{-4q}.
\end{equation*}
\end{corollary}
\noindent The proof of the corollary is based on simple upper estimate of $ W_{\vec{h}}(f)$.

\noindent Introduce
\begin{align*}
\mU_{\vec{h}}(f)&=\Big[\tfrac{\kappa_{\vec{h}}(f) \|f\|_2^2 \ln(m)}{m}\Big]^{\frac{1}{2}}+\Big[\tfrac{\|f\|_2^2\ln(m) }{m^2V_{\vec{h}}}\Big]^{\frac{1}{2}}
+\tfrac{\kappa_{\vec{h}}(f) \ln(m)}{m}+
\tfrac{\Upsilon_{\vec{h}}(m)\ln(m)}{m^2V_{\vec{h}}}
\\*[2mm]
\kappa_{\vec{h}}(f)&=\big(\tfrac{256s\ln(m)}{mV_{h}}\vee\|f\|_\infty\big)\mathrm{1}_{\bar{H}^d_m}\big(\vec{h}\big)
+\alpha_m^{-1}\mathrm{1}_{H^d_m}\big(\vec{h}\big).
\end{align*}

\begin{corollary}
\label{cor2:th2}
For any $m\geq 9$, $q\geq 1$ and   $f\in\bB_\infty\big(\alpha_m^{-1}\big)$
\begin{equation*}
\bE_f\bigg\{\sup_{\vec{h}\in\cH^d_m}\Big(\big|\widehat{U}_{\vec{h}}-\bE_{f}[\widehat{U}_{\vec{h}}]\big|-\Omega_q(U)\mU_{\vec{h}}(f)\Big)_+^q\bigg\}
\leq \Lambda^*_q(U)[\ln(m)]^dm^{-4q},
\end{equation*}
where $\Omega_q(U)=3^{d+1}(24q+4)^2\varpi_U$.

\end{corollary}
\noindent The proof  is based on simple upper estimates of $ J_{\vec{h}}(f)$ and $ W^*_{\vec{h}}(f)$.

\subsection{Random upper function for decoupled  $U$-statistics.}
\label{sec:subsec-U-stat-random}

The idea behind the construction below consists in the replacing the quantities $W^*_{h}(f)$ and  $J_{h}(f)$   appeared
in the description of the upper function $\{U^*_{\vec{h}}(f), \vec{h}\in\cH^d_m\}$ (depending on $f$ and therefore unknown) by their empirical counterparts.
Similarly to (\ref{random-J}) define for any $\vec{h}\in\cH^d_m$
\begin{equation*}
\widehat{J}_{\vec{h}}=\frac{1}{m^2}\sum_{i,j=1}^m|U_{\vec{h}}(Y_i-X_j)|
\end{equation*}
and let $\widehat{W}_{\vec{h}}$ be defined in (\ref{random-J}) with $s=12q+2$.
Note that $J_{\vec{h}}(f)=\bE_f\{\widehat{J}_{\vec{h}}\}$ and  $\widehat{W}_{\vec{h}}$ seems to be  a natural estimator for $W^*_{\vec{h}}(f)$. Set finally
\begin{align*}
\cU_{\vec{h}}&=\Big[\tfrac{16s \varpi_U\widehat{W}_{\vec{h}}\widehat{J}_{\vec{h}} \ln(m)}{m}\Big]^{\frac{1}{2}}+\Big[\tfrac{16s\varpi_U \widehat{J}_{\vec{h}}\ln(m) }{m^2V_{\vec{h}}}\Big]^{\frac{1}{2}}+\tfrac{147s\varpi_U \widehat{W}_{\vec{h}}\ln(m)}{m}
+\tfrac{s^2\varpi_U \Upsilon_{\vec{h}}(m)\ln(m)}{m^2V_{\vec{h}}}.
\end{align*}

\begin{theorem}
\label{th:random-upper-function}
For any $m\geq 9$, $q\geq 1$ and   $f\in\bB_\infty\big(\alpha_m^{-1}\big)$
\begin{eqnarray*}
&&\bE_f\bigg\{\sup_{\vec{h}\in \cH_m^d}\big(\cU_{\vec{h}}-2U^*_{\vec{h}}(f)\big)^q_+\bigg\}\leq B_1[2\ln(m)]^{d/2}m^{-2q};
\\*[2mm]
&&\bE_f\bigg\{\sup_{\vec{h}\in \cH_m^d}\big(U^*_{\vec{h}}(f)-9\cU_{\vec{h}}(f)\big)^q_+\bigg\}\leq B_2[2\ln(m)]^{d/2}m^{-2q}.
\end{eqnarray*}

\end{theorem}
Here  $B_1$ and $B_2$ are the constants which  depend on $q$ and $U$ only.
Their explicit expressions can be extracted from the corresponding proofs but it lies beyond the scope of the paper.

\section{Proof of main results.}

\subsection{Proof of Theorem \ref{th:oracle-inequality}.} Set for any $\bB_\infty\big(\alpha_m^{-1}\big)$
\begin{align*}
 \cR^q(f)&=\bE_{f}\Big\{\big(\widehat{N}_{\vec{h}^*}-\|f\|_2^2\big)^q\Big\}.
 \end{align*}
 Our first goal is to prove the following result. There exists $\mathbf{c}>0$ depending only on $\cK$  such for any $f\in\bB_\infty\big(\alpha_m^{-1}\big)$
 \begin{align}
 \label{eq1:proof-or-ineq}
 \mathbf{c}^{-1}\cR(f)&\leq \min_{\vec{h}\in\cH^d_m}\Big\{\big\|B_{\vec{h}}\big\|_2^2+2\max_{\vec{w}\in \cH_m^d}\Big|\big\|B_{\vec{h}\vee\vec{w}}\big\|^2_2-\big\|B_{\vec{w}}\big\|^2_2\Big|
 +\mU_{\vec{h}}(f)\Big\}
 \nonumber\\*[2mm]
  &\quad +[\ln(m)]^{\frac{d}{2q}}m^{-2}.
 \end{align}
 The latter result can be viewed as an oracle inequality for the estimation of squared $\bL_2$-norm.

 \paragraph{Proof of (\ref{eq1:proof-or-ineq}).} The proof is based an application of an abstract oracle inequality established in \cite{lepski2018}, Theorem 1, see also \cite{lepski2022}.

 $\mathbf{1^0.}$ To apply the aforementioned result we have first  to match our notations to those in \cite{lepski2018}.
 Thus, $n=m$, $A(f)=\|f\|_2^2$, $\bF=\bB_\infty\big(\alpha^{-1}_m\big)$, $\mH_n=\cH^d_m$,\; $\mh=\vec{h}$,\; $\eta=\vec{w}$,
  $(\mS_1,\ell)=(\bR,|\cdot|)$ and
 \begin{gather*}
 \widehat{A}_{\mh}=\widehat{N}_{\vec{h}}, \quad \widehat{A}_{\mh,\eta}=\widehat{N}_{\vec{h}\vee\vec{w}},\quad \Lambda_{\mh}(f)=\bE_f\big\{\widehat{N}_{\mh}\big\},\quad \Lambda_{\mh,\eta}(f)=\bE_f\big\{\widehat{N}_{\vec{h}\vee\vec{w}}\big\},
 \\*[3mm]
 \Psi_n(\mh)=9\widehat{\cU}_{\vec{h}},\quad\; \e_n=\big(2^{d+q-1}\mathbf{B\mathbf{}}[\ln(m)]^{d/2}m^{-2q}\big)^{\frac{1}{q}}=:\e_m,
 \end{gather*}
 where $\mathbf{B}=\Lambda^*_q(T,Q)+B_2$ and $\Lambda^*_q(T,Q)$ and $B_2$ are  the constants from Corollary \ref{cor1:th2} and
 Theorem \ref{th:random-upper-function} respectively.

 \par

 \noindent The oracle inequality in \cite{lepski2018} is proved under two hypotheses. The first one, called ${\bf A^{\text{permute}}}$, that is
 $\widehat{A}_{\mh,\eta}\equiv\widehat{A}_{\eta,\mh}$,  for any $\mh, \eta\in\mH_n$, is obviously fulfilled in our case because $\vec{h}\vee\vec{w}=\vec{w}\vee\vec{h}$.

 The second one, called  ${\bf A^{\text{upper}}}$, reduced to our consideration, consists in the following.
For any $m\geq 9$ one has to check that
\begin{align*}
(\mathbf{i}) &\;\;\sup_{f\in\bB_\infty(Q)}\bE_{f}\bigg(\max_{\vec{h}\in\cH^d_m}\Big[\big|\widehat{N}_{\vec{h}}-
\bE_f\big\{\widehat{N}_{\vec{h}}\big\}\big|-9\widehat{\cU}_{\vec{h}}\Big]_+^{q}\bigg)\leq \e_m^{q};
\\*[2mm]
(\mathbf{ii})&\;\;\sup_{f\in\bB_\infty(Q)}\bE_{f}\bigg(\max_{\vec{h},\vec{w}\in\cH^d_m}\Big[\big|\widehat{N}_{\vec{h}\vee\vec{w}}-
\bE_f\big\{\widehat{N}_{\vec{h}\vee\vec{w}}\big\}\big|-9\big(\widehat{\cU}_{\vec{h}}\wedge\widehat{\cU}_{\vec{w}}\big)\Big]_+^{q}\bigg)\leq \e_m^{q}.
\end{align*}

$\mathbf{2^0.}\;$ Let us check $(\mathbf{i})$ and $(\mathbf{ii})$. Set
$$
U^{**}_{\vec{h}}(f)=\max_{\vec{w}\in\cH_m^d}U^{*}_{\vec{h}\vee\vec{w}}(f),
$$
where $U^*_{\vec{h}}(f)$ is defined in (\ref{eq:definition-of-U}) with $U=T$,
and let
\begin{align*}
\Delta_f&=\bE_{f}\bigg(\sup_{\vec{h}\in\cH^d_m}\Big[\big|\widehat{N}_{\vec{h}}-
\bE_f\big\{\widehat{N}_{\vec{h}}\big\}\big|-U^{**}_{\vec{h}}(f)\Big]_+^{q}\bigg);
\\*[2mm]
\delta_f&=\bE_{f}\bigg(\sup_{\vec{h}\in\cH^d_m}\Big[U^{**}_{\vec{h}}(f)-9\widehat{\cU}_{\vec{h}}\Big]_+^{q}\bigg).
\end{align*}
Since $U^{*}_{\vec{h}}(f)\leq U^{**}_{\vec{h}}(f)$, we have in view of  Corollary \ref{cor1:th2}
\begin{equation}
\label{eq1:proof-th:oracle-inequality}
\sup_{f\in\bB_\infty(Q)}\Delta_f\leq \Lambda^*_q(T)[2\ln(m)]^dm^{-4q}.
\end{equation}
Using the trivial inequality
\begin{equation}
\label{eq:triv-ineq1}
\big(\sup_{\kappa}a_\kappa-\sup_{\kappa}b_\kappa\big)_{+}\leq \sup_{\kappa}\big(a_\kappa-b_\kappa\big)_+,
\end{equation}
verifying for arbitrary collections of real numbers, we get
$$
\max_{\vec{h}\in\cH^d_m}\Big[U^{**}_{\vec{h}}(f)-9\widehat{\cU}_{\vec{h}}\Big]_+\leq \max_{\vec{h}\in\cH^d_m}\sup_{\vec{w}\in\cH^d_m}
\Big[U^{*}_{\vec{h}\vee\vec{w}}(f)-9\cU_{\vec{h}\vee\vec{w}}\Big]_+.
$$
Taking into account that $\vec{h}\vee\vec{w}\in\cH_m^d$ for any $\vec{h},\vec{w}\in\cH_m^d$ we obtain
\begin{equation}
\label{eq2:proof-th:oracle-inequality-bis}
\max_{\vec{h}\in\cH^d_m}\Big[U^{**}_{\vec{h}}(f)-9\widehat{\cU}_{\vec{h}}\Big]_+\leq \max_{\vec{h}\in\cH^d_m}
\Big[U^{*}_{\vec{h}}(f)-9\cU_{\vec{h}}\Big]_+.
\end{equation}
Applying the  second statement of Theorem \ref{th:random-upper-function} with $U=T$, we assert that
\begin{equation}
\label{eq2:proof-th:oracle-inequality}
\sup_{f\in\bB_\infty(Q)}\delta_f\leq B_2[2\ln(m)]^{d/2}m^{-2q}.
\end{equation}

$\mathbf{2^0a.}\;$ We deduce from (\ref{eq1:proof-th:oracle-inequality}) and (\ref{eq2:proof-th:oracle-inequality}), using $(a+b)^q\leq 2^{q-1}(a^q+b^q), a,b\geq 0$,
\begin{align*}
\bE_{f}\bigg(\max_{\vec{h}\in\cH^d_m}\Big[\big|\widehat{N}_{\vec{h}}&-
\bE_f\big\{\widehat{N}_{\vec{h}}\big\}\big|-9\widehat{\cU}_{\vec{h}}\Big]_+^{q}\bigg)\leq 2^{q-1}\Big[\Delta_f+\delta_f\Big]
\\*[2mm]
&\leq 2^{d+q-1}\big[\Lambda^*_q(T,Q)+B_2\big][\ln(m)]^dm^{-2q}=\e^{q}_m
\end{align*}
and $(\mathbf{i})$ is verified since the right hand side of the obtained inequality is independent of $f$.

$\mathbf{2^0b.}\;$ Since $\vec{h}\vee\vec{w}\in\cH_m^d$ for any $\vec{h},\vec{w}\in\cH_m^d$ we first have
\begin{align*}
\max_{\vec{h},\vec{w}\in\cH^d_m}&\Big[\big|\widehat{N}_{\vec{h}\vee\vec{w}}-
\bE_f\big\{\widehat{N}_{\vec{h}\vee\vec{w}}\big\}\big|-U^*_{\vec{h}\vee w}(f)\Big]_+
\\*[2mm]
&=\max_{\vec{h}\in\cH^d_m}\Big[\big|\widehat{N}_{\vec{h}}-
\bE_f\big\{\widehat{N}_{\vec{h}}\big\}\big|-U^*_{\vec{h}}(f)\Big]_+.
\end{align*}
On the other hand $U^*_{\vec{h}\vee \vec{w}}(f)\leq U^{**}_{\vec{h}}(f)\wedge U^{**}_{ \vec{w}}(f)$ for any $\vec{h},\vec{w}\in\cH_m^d$ and, therefore,
\begin{equation}
\label{eq3:proof-th:oracle-inequality}
\bE_{f}\bigg(\max_{\vec{h},\vec{w}\in\cH^d_m}\Big[\big|\widehat{N}_{\vec{h}\vee\vec{w}}-
\bE_f\big\{\widehat{N}_{\vec{h}\vee\vec{w}}\big\}\big|-U^{**}_{\vec{h}}(f)\wedge U^{**}_{\vec{w}}(f)\Big]_+^{q}\bigg)\leq \Delta_f.
\end{equation}
Next, in view of the trivial inequality
$$
(a\wedge b-c\wedge d)_+\leq (a-c)_+ \vee (b-d)_+,\quad\;\forall a,b,c,d\in\bR,
$$
we get, using (\ref{eq2:proof-th:oracle-inequality-bis}),
\begin{align*}
\Big[U^{**}_{\vec{h}}(f)\wedge U^{**}_{\vec{w}}(f)-&9\widehat{\cU}_{\vec{h}}\wedge \widehat{\cU}_{\vec{w}}\Big]_+
\leq \Big[U^{**}_{\vec{h}}(f)-9\widehat{\cU}_{\vec{h}}\Big]_+\vee \Big[U^{**}_{\vec{w}}(f)-9\widehat{\cU}_{\vec{w}}\Big]_+
\\*[2mm]
&\leq \max_{\vec{h}\in\cH^d_m}\Big[U^{**}_{\vec{h}}(f)-9\widehat{\cU}_{\vec{h}}\Big]_+
\leq \sup_{\vec{h}\in\cH^d_m}
\Big[U^{*}_{\vec{h}}(f)-9\cU_{\vec{h}}\Big]_+.
\end{align*}
Hence
$$
\max_{\vec{h},\vec{w}\in\cH^d_m}\Big[U^{**}_{\vec{h}}(f)\wedge U^{**}_{\vec{w}}(f)-9\widehat{\cU}_{\vec{h}}\wedge \widehat{\cU}_{\vec{w}}\Big]_+
\leq \sup_{\vec{h}\in\cH^d_m}
\Big[U^{*}_{\vec{h}}(f)-9\cU_{\vec{h}}\Big]_+
$$
that yields together with (\ref{eq3:proof-th:oracle-inequality})
\begin{align*}
\bE_{f}\bigg(\sup_{\vec{h},\vec{w}\in\cH^d_m}&\Big[\big|\widehat{N}_{\vec{h}\vee\vec{w}}-
\bE_f\big\{\widehat{N}_{\vec{h}\vee\vec{w}}\big\}\big|-9\big(\widehat{\cU}_{\vec{h}}\wedge\widehat{\cU}_{\vec{w}}\big)\Big]_+^{q}\bigg)
\leq 2^{q-1}\Big[\Delta_f+\delta_f\Big]
\\*[2mm]
&\leq 2^{d+q-1}\big[\Lambda^*_q(T)+B_2\big][\ln(m)]^dm^{-2q}=\e^q_m.
\end{align*}
The assumption $(\mathbf{ii})$ is checked.

$\mathbf{3^0.}\;$ Applying Theorem 1 from \cite{lepski2018} we get for any $f\in\bB_\infty\big(\alpha_m^{-1}\big)$
\begin{equation}
\label{eq4:proof-th:oracle-inequality}
\cR(f)\leq \inf_{\vec{h}\in\cH_m^d}\left\{\cB_{\vec{h}}(f)+45\psi_{\vec{h}}(f)\right\}+6\e_m,
\end{equation}
where
\begin{align*}
\cB_{\vec{h}}(f)&=\big|\bE_f\big\{\widehat{N}_{\vec{h}}\big\}-\|f\|_2^2\big|+2\max_{\vec{w}\in \cH_m^d}
\big|\bE_f\big\{\widehat{N}_{\vec{h}\vee\vec{w}}\big\}-\bE_f\big\{\widehat{N}_{\vec{w}}\big\}\big|,
\end{align*}
and $\psi_{\vec{h}}(f)=\big(\bE_{f}\big\{\widehat{\cU}^q_{\vec{h}}\big\}\big)^{1/q}$.

$\mathbf{3^0a.}\;$ In view of Lemma \ref{lem:representation}
\begin{align*}
\big|\bE_f\big\{\widehat{N}_{\vec{h}}\big\}-\|f\|_2^2\big|&=\big\|B_{\vec{h}}\big\|_2^2;
\\*[2mm]
\max_{\vec{w}\in \cH_m^d}
\big|\bE_f\big\{\widehat{N}_{\vec{h}\vee\vec{w}}\big\}-\bE_f\big\{\widehat{N}_{\vec{w}}\big\}\big|
&=\max_{\vec{w}\in \cH_m^d}\Big|\big\|B_{\vec{h}\vee\vec{w}}\big\|^2_2-\big\|B_{\vec{w}}\big\|^2_2\Big|.
\end{align*}
and, therefore,
\begin{equation}
\label{eq5:proof-th:oracle-inequality}
\cB_{\vec{h}}(f)=\big\|B_{\vec{h}}\big\|_2^2+2\max_{\vec{w}\in \cH_m^d}\Big|\big\|B_{\vec{h}\vee\vec{w}}\big\|^2_2-\big\|B_{\vec{w}}\big\|^2_2\Big|.
\end{equation}

$\mathbf{3^0b.}\;$  Applying (\ref{eq:triv-ineq1}) we obtain similarly to (\ref{eq2:proof-th:oracle-inequality-bis})
\begin{equation}
\label{eq5000:proof-th:oracle-inequality}
\max_{\vec{h}\in \cH_m^d}\big[\widehat{\cU}_{\vec{h}}-2U^{**}_{\vec{h}}(f)\big]_+\leq
\max_{\vec{h}\in \cH_m^d}\big[\cU_{\vec{h}}-2U^{*}_{\vec{h}}(f)\big]_+.
\end{equation}
Hence, in view of the first assertion of Theorem \ref{th:random-upper-function}  we have
\begin{equation}
\label{eq6:proof-th:oracle-inequality}
\psi_{\vec{h}}(f)\leq 2\big[U^{**}_{\vec{h}}(f)\big] +\big(B_1[2\ln(m)]^{d/2}\big)^{\frac{1}{q}}m^{-2}.
\end{equation}
Since $U^{*}_{\vec{h}}(f)\leq \Omega_q(T)\mU_{\vec{h}}(f)$  for any $\vec{h}\in \cH_m^d$ (see the proof of Corollary \ref{cor2:th2}) and $f\in\bB_\infty\big(\alpha_m^{-1}\big)$ and because
$\mU_{\vec{h}}(f)\geq \mU_{\vec{h}\vee\vec{w}}(f)$ whatever  $\vec{h},\vec{w}\in \cH_m^d$,  we derive from (\ref{eq6:proof-th:oracle-inequality})
\begin{equation}
\label{eq7:proof-th:oracle-inequality}
\psi_{\vec{h}}(f)\leq 2\Omega_q(T)\mU_{\vec{h}}(f) +\big(B_1[2\ln(m)]^{d/2}\big)^{\frac{1}{q}}m^{-2}.
\end{equation}
The required inequality  (\ref{eq1:proof-or-ineq}) follows now
from (\ref{eq4:proof-th:oracle-inequality}), (\ref{eq5:proof-th:oracle-inequality}) and (\ref{eq7:proof-th:oracle-inequality}).

\paragraph{Completion of the proof of the theorem.} We obviously have for any density $f$
$$
\Big|\sqrt{|\widehat{N}_{\vec{h}^*}|}-\|f\|_2\Big|\leq \tfrac{|\widehat{N}_{\vec{h}^*}-\|f\|^2_2|}{\|f\|_2}\bigwedge
\sqrt{|\widehat{N}_{\vec{h}^*}-\|f\|_2|}.
$$
Hence
$$
 \cR^{(q)}_n\big[\widehat{N}_{\vec{h}^*},\|f\|_2\big]\leq \frac{\cR(f)}{\|f\|^2_2}\bigwedge \Big[\bE_{f}\Big\{\big|\widehat{N}_{\vec{h}^*}-\|f\|_2\big|^{q/2}\Big\}\Big]^{\frac{1}{q}}.
$$
Applying Cauchy-Schwartz inequality to the second term in the right hand side of the obtained inequality, we get
\begin{equation}
\label{eq8:proof-th:oracle-inequality}
\cR^{(q)}_n\big[\widehat{N}_{\vec{h}^*},\|f\|_2\big]\leq \frac{\cR(f)}{\|f\|_2}\bigwedge \sqrt{\cR(f)}.
\end{equation}
The assertion of the theorem follows from (\ref{eq1:proof-or-ineq}) and (\ref{eq8:proof-th:oracle-inequality}).
\epr

\subsection{Proof of Corollary \ref{cor2:th-oracle}.}
\label{sec:Proof-of-cor}
Let us recall the general form of the Young inequality, see \cite{Folland}, Theorem 8.9.
\begin{lemma}[Young inequality]
 Let $1\leq s, t, p\leq \infty$ and
$1+1/p=1/s+1/t$. If $w\in \bL_s\big(\bR^d\big)$ and $g\in \bL_t\big(\bR^d\big)$ then $w\star g\in \bL_p\big(\bR^d\big)$ and
\[
 \|w\star g\|_p \leq \|w\|_s \|g\|_t,
\]
where "$\star$" stands for convolution operator.
\end{lemma}

\paragraph{Proof of Corollary \ref{cor2:th-oracle}.} For any $\vec{h}\in\cH_m^d$ introduce
$$
\cH_m^d\big(\vec{h}\big)=\big\{\vec{w}\in\cH_m^d:\; \exists j\in\{1,\ldots,d\}:\; w_j<h_j\big\},\quad \bar{\cH}_m^d\big(\vec{h}\big)=\cH_m^d\setminus
\cH_m^d\big(\vec{h}\big).
$$
It is obvious that for any $\vec{h}\in\cH_m^d$
\begin{equation}
\label{eq100:proof-cor2-th-oracle}
\sup_{\vec{w}\in\bar{\cH}_m^d(\vec{h})}\Big|\big\|B_{\vec{h}\vee\vec{w}}\big\|^2_2-\big\|B_{\vec{w}}\big\|^2_2\Big|=0.
\end{equation}

\noindent We will need some additional notations and results. For any $J\subseteq\{1,\ldots d\}$ and  $y\in\bR^d$ set\footnote{Furthermore $|J|$ is used for the cardinality of $J$.} $y_J=\{y_j,\;j\in J\}\in\bR^{|J|}$ and we will write $y=\big(y_J,y_{\bar{J}}\big)$, where as usual $\bar{J}=\{1,\ldots d\}\setminus J$.

\noindent For any  $\vec{h}\in\bR_{+}^d$ and $J\neq\emptyset$ set
$
K_{\vec{h},J}(u_J)=\prod_{j\in J}h^{-1}_j\cK\big(u_j/h_j\big)
$
and define for any locally-integrable function $g:\bR^d\to\bR$
$$
\big[K_{\vec{h}}\star g\big]_{J}(y)=\int_{\bR^{|J|}}K_{\vec{h},J}(u_J-y_{J})g\big(u_J,y_{\bar{J}}\big)
\rd u_{J},\quad y\in\bR^d.
$$
By convention, $\big[K_{\vec{h}}\star g\big]_{\emptyset}\equiv g$.
The following result is a trivial consequence of the Young inequality, applying with $p=t$, $s=1$, and the Fubini theorem. For any $J\subseteq\{1,\ldots d\}$ and $p\in[1,\infty]$.
\begin{equation}
\label{eq101:proof-cor2-th-oracle}
\Big\|\big[K_{\vec{h}}\star g\big]_{J}\Big\|_{p}\leq \|\cK\|^{|J|}_{1}\|g\|_{p},\quad\forall\vec{h}\in\bR^d_+.
\end{equation}
For any $\vec{h},\vec{w}\in\cH_m^d$ denote
\begin{equation*}
B_{\vec{h},\vec{w}}(x)=B_{\vec{h}\vee \vec{w}}(x)-B_{ \vec{w}}(x),\quad x\in\bR^d.
\end{equation*}
and   let
$
\cJ=\{j\in\{1,\ldots,d\}:\; h_j=w_j\vee h_j\}.
$
Assume that  $\vec{w}\in\cH_m^d\big(\vec{h}\big)$. Then, $\cJ\neq\emptyset$ and let $\cJ=\{j_1<\cdots<j_k\}$, where $k=|\cJ|$.
Set at last $\cJ_l=\{j_1,\ldots, j_l\}$, $l=1,\ldots,k$.

\begin{lemma}
\label{lem:techniclal-proof-cor2}
 For any $\vec{h}\in\cH_m^d$,  $\vec{w}\in\cH_m^d\big(\vec{h}\big)$,    $x\in\bR^d$ and   $f:\bR^d\to\bR$
\begin{eqnarray*}
B_{\vec{h},\vec{w}}(x)&=& \sum_{l=1}^k\big\{\big[K_{\vec{w}}\star b_{h_l,l}\big]_{\bar{\cJ_l}}(x)+\big[K_{\vec{w}}\star b_{w_l,l}\big]_{\bar{\cJ_l}}(x)\big\};
\\
B_{\vec{h}}(x)&=& \sum_{l=1}^d
\left[K_{\vec{h}}\star b_{h_l,l} \right]_{\bar{J}_l}(x),\qquad J_l=\{1,\ldots, l\}.
\end{eqnarray*}
\end{lemma}
\noindent Proof of the lemma is postponed to Section \ref{sec:Proof of technical lemmas}.

\vskip0.1cm

$\mathbf{1^0.}\;$ Using (\ref{eq101:proof-cor2-th-oracle}) we derive from the first assertion of Lemma \ref{lem:techniclal-proof-cor2} that
for any $\vec{h}\in\cH_m^d$ and   $\vec{w}\in\cH_m^d\big(\vec{h}\big)$
$$
\big\|B_{\vec{h},\vec{w}}\big\|_2\leq \sum_{l=1}^k\|\cK\|_1^{d-|\cJ_l|}\big\{\big\|b_{h_l,l}\big\|_2+\big\|b_{w_l,l}\big\|_2\big\}.
$$
Noting that $\|\cK\|_1\geq 1$ since $\int \cK=1$ and that $w_l\leq h_l$ in view of the definition of $\cJ$ we can assert that
$$
\big\|B_{\vec{h},\vec{w}}\big\|_2\leq 2\sum_{l=1}^k\|\cK\|_1^{d-l}\sup_{a\leq h_l}\big\|b_{a,l}\big\|_2\leq 2A(\cK)\max_{l=1,\ldots,d}\sup_{a\leq h_l}\big\|b_{a,l}\big\|_2,
$$
where we have put $A(\cK)=\sum_{l=1}^d\|\cK\|_1^{l}$.
Noting that the right hand side of the latter inequality is independent of $\vec{w}$ and taking into account (\ref{eq100:proof-cor2-th-oracle})
we obtain
\begin{equation}
\label{eq102:proof-cor2-th-oracle}
\sup_{\vec{w}\in\cH_m^d}\big\|B_{\vec{h},\vec{w}}\big\|_2\leq 2A(\cK)\mb_{\vec{h}}(f),\quad \forall \vec{h}\in\cH_m^d.
\end{equation}
Similarly we derived from the second assertion of the lemma that
\begin{equation}
\label{eq103:proof-cor2-th-oracle}
\big\|B_{\vec{h}}\big\|_2\leq A(\cK)\mb_{\vec{h}}(f),\quad \forall \vec{h}\in\cH_m^d.
\end{equation}

$\mathbf{2^0.}\;$
We obviously have in view of Cauchy-Schwarz inequality
\begin{equation}
\label{eq1:proof-cor2-th-oracle}
\cB_{\vec{h},\vec{w}}:=\Big|\big\|B_{\vec{h}\vee\vec{w}}\big\|^2_2-\big\|B_{\vec{w}}\big\|^2_2\Big|\leq\big\|B_{\vec{h},\vec{w}}\big\|_2^2+
2\big\|B_{\vec{w}}\big\|_2\big\|B_{\vec{h},\vec{w}}\big\|_2.
\end{equation}
Applying triangle and the Young inequalities we get
\begin{align*}
\big\|B_{\vec{w}}\big\|_2^2&=\Big\|\int_{\bR^d} K_{\vec{w}}(\cdot-u)f(u)\rd u -f\Big\|_2^2\leq (\|K\|_1+1)^2\|f\|_2^2
\\*[2mm]
&\leq 4\|\cK\|_1^{2d}\|f\|_2^2.
\end{align*}
To get the last equality we have taken into account the product structure of the function $K$ and that $\|\cK\|_1\geq 1$.
We deduce from (\ref{eq1:proof-cor2-th-oracle}) that
\begin{equation}
\label{eq2:proof-cor2-th-oracle}
\Big|\big\|B_{\vec{h}\vee\vec{w}}\big\|^2_2-\big\|B_{\vec{w}}\big\|^2_2\Big|\leq\big\|B_{\vec{h},\vec{w}}\big\|_2^2+
4\|\cK\|_1^{d}\|f\|_2\big\|B_{\vec{h},\vec{w}}\big\|_2.
\end{equation}
The assertion of the corollary follows now from
 (\ref{eq102:proof-cor2-th-oracle}),  (\ref{eq103:proof-cor2-th-oracle}),  (\ref{eq2:proof-cor2-th-oracle}) and  Theorem \ref{th:oracle-inequality}.
\epr

 \subsection{Proof of Theorem \ref{th:adaptive-anisotr}.} From now on $c_1,c_2,\ldots,$ denote the constants independent of $n=2m$.
Also we will assume that $m$ is large enough in order to guarantee that $Q\leq \alpha^{-1}_m$ and, therefore, $\bB_\infty(Q)\subseteq\bB_\infty\big(\alpha^{-1}_m\big)$.

Introduce the following notation. For  $j=1,\ldots,d$ let
\begin{equation*}
\label{eq:2-kappa}
p_j:=\left\{
\begin{array}{ll}
\frac{r_j}{r_j-1}, &r_j\geq 2;
\\*[2mm]
r_j, &r_j<2,
\end{array}
\right.
\quad
z(\varsigma):=\left\{
\begin{array}{ll}
\;\;\frac{1}{\tau_\varsigma(1)},\quad & \vec{r}\in[2,\infty]^d;
\\*[2mm]
\frac{1}{2-\tau_\varsigma(\infty)},\; &\vec{r}\in[1,2]^d.
\end{array}
\right.
\end{equation*}
 The basic element of the proof of the theorem is the upper bound for the bias term formulated   below.
\begin{proposition}
\label{prop:bound-for-bias}
For any $\vec{\beta}\in (0,b)^d$, $\vec{r}\in [1,2]^d\cup(2,\infty]^d$, $\vec{L}\in (0,\infty)^d$,  $Q>0$,
$m\geq 9$ and any $f\in \bN_{\vec{r},d}(\vec{\beta},\vec{L})\cap\bB_\infty(Q)$
$$
\sup_{a\leq h_j} \big\|b_{a,j}\big\|_2\leq Ch^{\frac{\beta_jp_j}{2}}_{j},\; j=1,\ldots,d,\quad \forall \vec{h}\in\bR^d_+,
$$
where $C=C_1(1+\|\cK\|_1)(1\vee L_j)\sqrt{1\vee Q}$ and $C_1$ depends on $\cK$ and $b$ only.

\end{proposition}

\paragraph*{Proof of the theorem.}
Define $\boldsymbol{\vec{h}=(h_1,\ldots,h_d)}$ by
\begin{equation}
\label{eq:opt-bandwidth}
\boldsymbol{h_j}:=(\mu_m)^{\frac{2}{\beta_j p_j[1+1/\upsilon]}},\quad \tfrac{1}{\upsilon}:=\sum_{j=1}^d \tfrac{1}{p_j\beta_j}.
\end{equation}
Furthermore, without loss of generality, we will assume that $\boldsymbol{\vec{h}}\in\cH_m^d$. Otherwise, one can project $\boldsymbol{\vec{h}}$ coordinately
onto $\cH_m^d$. We begin the proof with several technical remarks.

\vskip0.1cm

$\mathbf{1^0.}\;$ Note that for any $ \varsigma\in(0,\infty)^d\times\big\{[1,2]^d\cup[2,\infty]^d\big\}$
\begin{equation}
\label{eq1:proof-th:adaptive-anisotr}
\mz^*(\varsigma)=z(\varsigma)\;\; \Leftrightarrow\;\; z(\varsigma)<\tfrac{1}{2};\qquad  \mz^*(\varsigma)=\tfrac{1}{2}\; \;\Leftrightarrow \;\; z(\varsigma)\geq\tfrac{1}{2}.
\end{equation}
Indeed, $\vec{r}\in[2,\infty]^d$ implies $\tau_\varsigma(2)\geq 1$ and, therefore, $\mz^*(\varsigma)=\frac{1}{\tau_\varsigma(1)}\wedge \frac{1}{2}$ while $z(\varsigma)=\frac{1}{\tau_\varsigma(1)}$ and (\ref{eq1:proof-th:adaptive-anisotr}) follows.
Note also that $\vec{r}\in[1,2]^d$ implies that $\tau_\varsigma(2)\leq 1$ and, therefore, $\mz^*(\varsigma)=z(\varsigma)<1/2$ if $\tau_\varsigma(\infty)<0$. On the other hand if $\tau_\varsigma(\infty)\geq 0$ then
$
z(\varsigma)\geq\tfrac{1}{2},
$
while $\mz^*(\varsigma)=\tfrac{1}{2}$
This completes the proof of (\ref{eq1:proof-th:adaptive-anisotr}).

$\mathbf{2^0.}\;$ Let us prove now that
\begin{equation}
\label{eq2:proof-th:adaptive-anisotr}
1+\tfrac{1}{\upsilon}=\tfrac{1}{z(\varsigma)}.
\end{equation}
Indeed if $\vec{r}\in[2,\infty]^d$ we have
$$
1+\tfrac{1}{\upsilon}=1+\sum_{j=1}^d\tfrac{r_j-1}{\beta_jr_j}=1+\tfrac{1}{\beta}-\tfrac{1}{\omega_\varsigma}=\tau_\varsigma(1).
$$
If $\vec{r}\in[1,2]^d$ then
$$
1+\tfrac{1}{\upsilon}=1+\sum_{j=1}^d\tfrac{1}{\beta_jr_j}=1+\tfrac{1}{\omega}=2-\tau(\infty),
$$
that completes the proof of (\ref{eq2:proof-th:adaptive-anisotr}).
 We derive from Proposition \ref{prop:bound-for-bias} and (\ref{eq2:proof-th:adaptive-anisotr})
\begin{equation}
\label{eq0:proof-th:adaptive-anisotr}
\sup_{f\in\bN_{\vec{r},d}(\vec{\beta},\vec{L})\cap\bB_\infty(Q)}\mb_{\vec{h}}(f)\leq C  (\mu_m)^{z(\varsigma)}.
\end{equation}

$\mathbf{3^0.}\;$ Note that
$
V_{\bsh}=(\mu_m)^{\frac{2}{1+\upsilon}}
$
and, therefore, in view of (\ref{eq2:proof-th:adaptive-anisotr})
\begin{equation}
\label{eq3:proof-th:adaptive-anisotr}
\mu^2_mV^{-1}_{\bsh}=(\mu_m)^{\frac{2}{1+1/\upsilon}}=(\mu_m)^{2z(\varsigma)}.
\end{equation}
 We derive from (\ref{eq3:proof-th:adaptive-anisotr}) that
$$
mV_{\bsh}=\big\{[\ln(m)]^{1-z(\varsigma)}\mathrm{1}\big(\mz^*(\varsigma)< 1/2\big)+[\ln(m)]^{2-2z(\varsigma)}\mathrm{1}\big(\mz^*(\varsigma)=1/2\big)\big\}
m^{2z(\varsigma)-1}.
$$
Thus, we can assert that for all $m$ large enough
\begin{equation}
\label{eq4:proof-th:adaptive-anisotr}
\bsh\in\left\{
\begin{array}{lll}
\mH^d_*,\;\;  &z(\varsigma)<\frac{1}{2};
\\*[2mm]
 \bar{H}^d_m,\; \; &z(\varsigma)\geq\frac{1}{2},
\end{array}
\right.
\end{equation}
In particular,
$$
\lim_{m\to\infty}\big|\tfrac{\ln(m)}{\ln(mV_{\bsh})}\big|=(1-2z(\varsigma))^{-1},\quad z(\varsigma)<1/2,
$$
and, therefore, for all  $m$ large enough
\begin{equation}
\label{eq5:proof-th:adaptive-anisotr}
\Upsilon_{\bsh}(m)\leq c_1\ln(m)\mathrm{1}\big(z(\varsigma)\geq 1/2\big)
+c_2\mathrm{1}\big(z(\varsigma)<1/2\big).
\end{equation}
Taking into account that $\mH^d_*\subset  H^d_m$  we also derive from (\ref{eq4:proof-th:adaptive-anisotr}) that for any $f\in\bB_\infty(Q)$ and all $m$ large enough
\begin{equation}
\label{eq6:proof-th:adaptive-anisotr}
\kappa_{\bsh}(f)\leq c_3Q\mathrm{1}\big(z(\varsigma)\geq 1/2\big)
+c_4\alpha^{-1}_m\mathrm{1}\big(z(\varsigma)< 1/2\big).
\end{equation}

$\mathbf{3^0.}\;$ Set for brevity $r_m=[\ln(m)]^{\frac{d}{2}}m^{-2}$. We obtain from (\ref{eq3:proof-th:adaptive-anisotr}), (\ref{eq5:proof-th:adaptive-anisotr}) and (\ref{eq6:proof-th:adaptive-anisotr}) that
\begin{align*}
\mU_{\bsh}(f)&\leq c_5\Big\{\|f\|_2 \sqrt{\tfrac{\ln(m)}{\alpha_m m}}+\tfrac{\ln(m)}{\alpha_m m}+\|f\|_2\mu_m^{z(\varsigma)}
+\mu_m^{2z(\varsigma)}+r_m\Big\},\quad z(\varsigma)< 1/2;
\\*[2mm]
\mU_{\bsh}(f)&\leq c_5\Big\{\|f\|_2 \sqrt{\mu_m}+\mu_m+\|f\|_2\mu_m^{z(\varsigma)}
+\mu_m^{2z(\varsigma)}+r_m\Big\},\quad\hskip1cm z(\varsigma)\geq 1/2.
\end{align*}
Remembering  that $\alpha^{-1}_m\leq\ln(m)$ and using that  $\mu_m\to 0$ we finally get
\begin{align}
\mU_{\bsh}(f)&\leq c_6\Big\{\|f\|_2\mu_m^{z(\varsigma)}
+\mu_m^{2z(\varsigma)}\Big\},\quad z(\varsigma)< 1/2;
\label{eq7:proof-th:adaptive-anisotr}
\\*[2mm]
\mU_{\bsh}(f)&\leq c_6\Big\{\|f\|_2 \sqrt{\mu_m}+\mu_m\Big\},\hskip0.7cm z(\varsigma)\geq 1/2.
\label{eq8:proof-th:adaptive-anisotr}
\end{align}
We derive from (\ref{eq0:proof-th:adaptive-anisotr}), (\ref{eq2:proof-th:adaptive-anisotr}), (\ref{eq7:proof-th:adaptive-anisotr}) and (\ref{eq8:proof-th:adaptive-anisotr}) that for any
$\varsigma\in(0,\infty)^d\times\big\{[1,2]^d\cup[2,\infty]^d\big\}$ and any $f\in\bN_{\vec{r},d}\big(\vec{\beta},\vec{L}\big)\cap\bB_\infty(Q)$
$$
\cO^*_f\leq c_7\Big\{\|f\|_2\mu_m^{z(\varsigma)}
+\mu_m^{2z(\varsigma)}\Big\},\quad z(\varsigma)< 1/2;
\\*[2mm]
$$
\begin{align*}
\cO^*_f&\leq c_7\Big\{\|f\|_2 \sqrt{\mu_m}+\mu_m+\|f\|_2\mu_m^{z(\varsigma)}
+\mu_m^{2z(\varsigma)}\Big\}
\\*[2mm]
&\leq 2c_7\Big\{\|f\|_2 \sqrt{\mu_m}+\mu_m\Big\},\quad z(\varsigma)\geq 1/2.
\end{align*}
In view of (\ref{eq1:proof-th:adaptive-anisotr})
$$
\mu_m^{z(\varsigma)}=\psi_m(\varsigma),\; z(\varsigma)< 1/2; \quad\; \mu_m=\psi_m(\varsigma),\; z(\varsigma)\geq 1/2
$$
and we obtain for any
$\varsigma\in(0,\infty)^d\times\big\{[1,2]^d\cup[2,\infty]^d\big\}$ and any $f\in\bN_{\vec{r},d}\big(\vec{\beta},\vec{L}\big)\cap\bB_\infty(Q)$
\begin{equation}
\label{eq9:proof-th:adaptive-anisotr}
\cO^*_f\leq c_8 \Big\{\|f\|_2 \psi_m(\varsigma)+\psi^2_m(\varsigma)\Big\}.
\end{equation}

$\mathbf{4^0.}\;$ Set $\bF=\bN_{\vec{r},d}\big(\vec{\beta},\vec{L}\big)\cap\bB_\infty(Q)\cap\bB_2\big(\psi_m(\varsigma)\big)$ and let $\bar{\bF}=\big[\bN_{\vec{r},d}\big(\vec{\beta},\vec{L}\big)\cap\bB_\infty(Q)\big]\setminus \bF$.
We have in view of (\ref{eq9:proof-th:adaptive-anisotr})
\begin{align}
\label{eq10:proof-th:adaptive-anisotr}
\sup_{f\in\bF}\sqrt{\cO^*_f}\leq c_9 \psi_m(\varsigma),\quad\; \sup_{f\in\bar{\bF}}\tfrac{\cO^*_f}{\|f\|_2}\leq c_9 \psi_m(\varsigma)
\end{align}
and, therefore, for any
$\varsigma\in(0,\infty)^d\times\big\{[1,2]^d\cup[2,\infty]^d\big\}$
$$
\sup_{f\in\bN_{\vec{r},d}\big(\vec{\beta},\vec{L}\big)\cap\bB_\infty(Q)}\Big[\tfrac{\cO^{*}_f}{\|f\|_2}\bigwedge \sqrt{\cO^{*}_f}\Big]\leq c_9 \psi_m(\varsigma).
$$
The assertion of the theorem follows now from Corollary \ref{cor2:th-oracle}.
\epr

\subsection{Proof of Theorem \ref{th:optim-adap-par-zone}.}
\label{sec:proof-th:optim-adap-par-zone} The proof is elementary and  furthermore, $c_1,c_2,\ldots,$ denote the constants independent of $m$.

 We derive from Propositions 3 and 4 in \cite{GL20b}, applied with $p=2$, that for any $f\in\bN_\vartheta$
$$
\bE_f\big\{\big(\widehat{N}_{\vec{\mathrm{h}}_m}-\|f\|_2^2\big)^2\big\}\leq c_1\Big[\sum_{j=1}^d\mathrm{h}_j^{2\blb\boldsymbol{p}}+\|f\|_2^2\big(\tfrac{1}{m}+\tfrac{1}{mV_{\vec{\mathrm{h}}}}\big)\Big]=c_2
\Big[\big(\tfrac{1}{m}\big)^{\frac{2\blb\boldsymbol{p}}{d}}+\tfrac{\|f\|_2^2}{m}\Big],
$$
where $\boldsymbol{p}=\frac{\bbr}{\bbr-1}\mathrm{1}_{[2,\infty]}(\bbr)+\bbr\mathrm{1}_{[1,2]}(\bbr)$. Noting that $\blb\boldsymbol{p}\geq d$ if $\theta\in\Theta^b_{(r),\text{param}}$ we finally get
\begin{equation}
\label{eq1:proof-th:optim-adap-par-zone}
R^2_m(f):=\bE_f\big\{\big(\widehat{N}_{\vec{\mathrm{h}}_m}-\|f\|_2^2\big)^2\big\}\leq c_3\Big[\tfrac{1}{m^2}+\tfrac{\|f\|_2^2}{m}\Big],\quad \forall f\in\bN_\vartheta,\;\forall \theta\in \Theta^b_{(r),\text{param}}.
\end{equation}
We have  for any density $f$
$$
\cR_{2m}\big[|\widehat{N}_{\vec{\mathrm{h}}_m}|^{\frac{1}{2}}; \|f\|_2\big]\leq \tfrac{R_m(f)}{\|f\|_2}\bigwedge \sqrt{R_m(f)}.
$$
Similarly to (\ref{eq10:proof-th:adaptive-anisotr}) we have
$$
\sqrt{R_m(f)}\leq c_4m^{-\tfrac{1}{2}},\;\; \forall f\in\bN_\vartheta\cap\bB_2\big(m^{-1/2}\big),
$$
and
$$\tfrac{R_m(f)}{\|f\|_2}\leq c_5m^{-\tfrac{1}{2}},\;\; \forall f\in\bN_\vartheta\cap\bar{\bB}_2\big(m^{-1/2}\big),
$$
and the assertion of the theorem follows.
\epr

\subsection{Proof of Theorem \ref{th:main-adap-isotrop}.} Later on $c_1,c_2,\ldots,D_1,D_2\ldots,$ denote the constants independent of $m$.
The initial step in our proof is to establish a crucial result presented in
(\ref{eq4-proof-th:main-adap-isotrop}) below.

\paragraph{Oracle inequality for losses.}
Setting
\begin{align*}
& \xi_1=\max_{\vec{h}\in\cH^d_m}\Big[\big|\widehat{N}_{\vec{h}}-
\bE_f\big\{\widehat{N}_{\vec{h}}\big\}\big|-9\widehat{\cU}_{\vec{h}}\Big]_+,
\\
& \xi_2=\max_{\vec{h},\vec{w}\in\cH^d_m}\Big[\big|\widehat{N}_{\vec{h}\vee\vec{w}}-
\bE_f\big\{\widehat{N}_{\vec{h}\vee\vec{w}}\big\}\big|-9\big(\widehat{\cU}_{\vec{h}}\wedge\widehat{\cU}_{\vec{w}}\big)\Big]_+,
\end{align*}
we derive from  \textit{the inequality (2.7) established in the the proof of Theorem 1 in \cite{lepski2018}}, applied with $\e_n^\prime=0$, that for any density $f\in\bB_\infty\big(\alpha^{-1}_m\big)$
\begin{equation}
\label{eq1-proof-th:main-adap-isotrop}
\big|\widehat{N}_{\vec{h}^*}-\|f\|^2_2\big|\leq \cB_{\vec{h}}(f)+45\widehat{\cU}_{\vec{h}}+3\xi_1+2\xi_2,\quad\forall\vec{h}\in\cH_m^d.
\end{equation}
Here, remind, in view of (\ref{eq5:proof-th:oracle-inequality})
\begin{equation*}
\cB_{\vec{h}}(f)=\big\|B_{\vec{h}}\big\|_2^2+2\max_{\vec{w}\in \cH_m^d}\Big|\big\|B_{\vec{h}\vee\vec{w}}\big\|^2_2-\big\|B_{\vec{w}}\big\|^2_2\Big|.
\end{equation*}
Setting also
$
\xi_3=\max_{\vec{h}\in \cH_m^d}\big[\cU_{\vec{h}}-2U^{*}_{\vec{h}}(f)\big]_+
$
we derive from (\ref{eq5000:proof-th:oracle-inequality}) that
\begin{equation}
\label{eq2-proof-th:main-adap-isotrop}
\widehat{\cU}_{\vec{h}}\leq 2U^{**}_{\vec{h}}(f) +\xi_3.
\end{equation}
Since $U^{*}_{\vec{h}}(f)\leq \Omega_4(T)\mU_{\vec{h}}(f)$  for any $\vec{h}\in \cH_m^d$ (see the proof of Corollary \ref{cor2:th2}) and $f\in\bB_\infty\big(\alpha_m^{-1}\big)$ and because
$\mU_{\vec{h}}(f)\geq \mU_{\vec{h}\vee\vec{w}}(f)$ whatever  $\vec{h},\vec{w}\in \cH_m^d$, we obtain from (\ref{eq2-proof-th:main-adap-isotrop})
$$
\widehat{\cU}_{\vec{h}}\leq 2\Omega_4(T)\mU_{\vec{h}}(f) +\xi_3.
$$
It yields together with (\ref{eq1-proof-th:main-adap-isotrop})
\begin{equation}
\label{eq3-proof-th:main-adap-isotrop}
\big|\widehat{N}_{\vec{h}^*}-\|f\|^2_2\big|\leq \cB_{\vec{h}}(f)+90\Omega_4(T)\mU_{\vec{h}}(f)+\eta,\quad\forall\vec{h}\in\cH_m^d,
\end{equation}
where we have denoted $\eta=3\xi_1+2\xi_2+45\xi_3$.

At last, in the proof of Corollary \ref{cor2:th-oracle} we have shown  that there exists $C(\cK)$ (explicitly known constant completely determined by the kernel $\cK$) such that
$$
\cB_{\vec{h}}(f)\leq C(\cK) \big[\mb^2_{\vec{h}}(f)+\|f\|_2\mb_{\vec{h}}(f)\big].
$$
Note that $T$ is completely determined by $\cK$, while  this function in its turn is completely determined by the function $D$ and the number $b$.
Hence, putting $C_b(D)=C(\cK)\vee [90\Omega_4(T)]$, we obtain from (\ref{eq3-proof-th:main-adap-isotrop}) for any $f\in\bB_\infty\big(\alpha_m^{-1}\big)$
\begin{equation}
\label{eq400-proof-th:main-adap-isotrop}
\big|\widehat{N}_{\vec{h}^*}-\|f\|^2_2\big|\leq C_b(D)\mathrm{O}^*_f+\eta,\quad\; \mathrm{O}^*_f=\min_{\vec{h}\in\cH^d_m}\Big\{\mb^2_{\vec{h}}(f)+\|f\|_2\mb_{\vec{h}}(f) +\mU_{\vec{h}}(f)\Big\}.
\end{equation}
The inequality (\ref{eq400-proof-th:main-adap-isotrop}) can be viewed as an \textit{oracle inequality for losses}.

\noindent Moreover, we have  shown (checking hypothesis ${\bf A^{\text{upper}}}$ within the proof of Theorem \ref{th:oracle-inequality}) that
\begin{equation*}
\sup_{f\in\bB_\infty(\alpha_m^{-1})}\bE_f\{\xi^q_i\}\leq 2^{d+q-1}\mathbf{B\mathbf{}}[\ln(m)]^{d/2}m^{-2q},\quad\;  i=1,2.
\end{equation*}
Also, we have in view of the first assertion of Theorem \ref{th:random-upper-function}
\begin{equation*}
\sup_{f\in\bB_\infty(\alpha_m^{-1})}\bE_f\big\{\xi_3^q\big\}\leq B_1[2\ln(m)]^{d/2}m^{-2q}.
\end{equation*}
Choosing $q=4$ and using Markov inequality, we derive from (\ref{eq400-proof-th:main-adap-isotrop}) that
\begin{equation}
\label{eq4-proof-th:main-adap-isotrop}
\sup_{f\in\bB_\infty(\alpha_m^{-1})}\bP_f\Big\{\big|\widehat{N}_{\vec{h}^*}-\|f\|^2_2\big|\geq 2C_b(D)\mathrm{O}^*_f+m^{-1}\Big\}\leq D_1[\ln(m)]^{d/2}m^{-4}.
\end{equation}

\paragraph{Proof of the theorem.} Introduce the random event
$$
\cA=\Big\{\big||\widehat{N}_{\vec{h}^*}|^{\frac{1}{2}}-|\widehat{N}_{\vec{\mathrm{h}}_m}|^{\frac{1}{2}}\big|\leq 2\ln(m)m^{-\frac{1}{2}}\Big\}
$$
and let $\bar{\cA}$ denote the complimentary event.  Consider separately two cases.

$\mathbf{1^0.}\;$ Assume first that $\theta\in\Theta^b_{(r),\text{nparam}}$. In view of the definition of the selection rule (\ref{eq:sel-rule-isotr})
\begin{align*}
\big||\widehat{N}_{\vec{\mathbf{h}}}|^{\frac{1}{2}}-\|f\|_2\big|&=\big||\widehat{N}_{\vec{\mathrm{h}}_m}|^{\frac{1}{2}}-\|f\|_2\big|\mathrm{1}_{\cA}+
\big||\widehat{N}_{\vec{h}^*}|^{\frac{1}{2}}-\|f\|_2\big|\mathrm{1}_{\bar{\cA}}\leq \big||\widehat{N}_{\vec{h}^*}|^{\frac{1}{2}}-\|f\|_2\big|
\\*[2mm]
&\quad+\big||\widehat{N}_{\vec{h}^*}|^{\frac{1}{2}}-
|\widehat{N}_{\vec{\mathrm{h}}_m}|^{\frac{1}{2}}\big|\mathrm{1}_{\cA}\leq \big||\widehat{N}_{\vec{h}^*}|^{\frac{1}{2}}-\|f\|_2\big|+2\ln(m)m^{-\frac{1}{2}}.
\end{align*}
Taking into account that $\psi_m(\varsigma)=\boldsymbol{\psi}_m(\theta)$ for any $\theta\in\Theta^b_{(r),\text{nparam}}$ and that
$$
\lim_{m\to\infty}\boldsymbol{\psi}^{-1}_m(\theta)\ln(m)m^{-\frac{1}{2}}=0,
$$
we derive from Theorem \ref{th:adaptive-anisotr}, applied with $q=4$ that for any $\theta\in\Theta^b_{(r),\text{nparam}}$
\begin{align}
\label{eq5-proof-th:main-adap-isotrop}
\limsup_{m\to\infty}\boldsymbol{\psi}^{-1}_m(\theta)\;\cR_{2m}\Big[|\widehat{N}_{\vec{\mathbf{h}}}|^{\frac{1}{2}}, \bN_\vartheta\Big]&\leq
\limsup_{m\to\infty}\boldsymbol{\psi}^{-1}_m(\theta)\;\cR_{2m}\Big[|\widehat{N}_{\vec{h}^*}|^{\frac{1}{2}}, \bN_\vartheta\Big]
\nonumber\\*[2mm]
&\leq \limsup_{m\to\infty}\boldsymbol{\psi}^{-1}_m(\theta)\;\cR^{(4)}_{2m}\Big[|\widehat{N}_{\vec{h}^*}|^{\frac{1}{2}}, \bN_\vartheta\Big]<\infty.
\end{align}

$\mathbf{2^0.}\;$ Assume now that $\theta\in\Theta^b_{(r),\text{param}}$. We have for any $f\in\bB_\infty(\alpha_m^{-1})$
\begin{align}
\label{eq6-proof-th:main-adap-isotrop}
\cR^2_{2m}\Big[|\widehat{N}_{\vec{\mathbf{h}}}|^{\frac{1}{2}}, \|f\|_2\Big]&\leq \cR^2_{2m}\Big[|\widehat{N}_{\vec{\mathrm{h}}_m}|^{\frac{1}{2}}, \|f\|_2\Big]+\bE_f\Big\{\big(|\widehat{N}_{\vec{h}^*}|^{\frac{1}{2}}-\|f\|_2\big)^2\mathrm{1}_{\bar{\cA}}\Big\}
\nonumber\\*[2mm]
&\leq \cR^2_{2m}\Big[|\widehat{N}_{\vec{\mathrm{h}}_m}|^{\frac{1}{2}}, \|f\|_2\Big]+\bigg(\cR^{(4)}_{2m}\Big[|\widehat{N}_{\vec{h}^*}|^{\frac{1}{2}}, \|f\|_2\Big]\bigg)^{2}\sqrt{\bP_f(\bar{\cA})}.
\end{align}
Here we also used  Cauchy-Schwartz inequality.
We have for any $\theta\in\Theta^b_{(r),\text{param}}$ in view of Theorems \ref{th:optim-adap-par-zone} and \ref{th:adaptive-anisotr} respectively
\begin{align*}
&\limsup_{m\to\infty}\sqrt{m}\;\cR_{2m}\Big[|\widehat{N}_{\vec{\mathrm{h}}_m}|^{\frac{1}{2}}, \bN_\vartheta\Big]<\infty,\qquad
\limsup_{m\to\infty}\sqrt{\tfrac{m}{\ln(m)}}\;\cR^{(4)}_{2m}\Big[|\widehat{N}_{\vec{h}^*}|^{\frac{1}{2}}, \bN_\vartheta\Big]<\infty.
\end{align*}
Hence, if we show that for any $\theta\in\Theta^b_{(r),\text{param}}$
\begin{align}
\label{eq7-proof-th:main-adap-isotrop}
\limsup_{m\to\infty}\ln^2(m)\sup_{f\in\bN_\vartheta}\bP_f(\bar{\cA})<\infty,
\end{align}
then we can derive from (\ref{eq6-proof-th:main-adap-isotrop}) that for any $\theta\in\Theta^b_{(r),\text{param}}$
\begin{align*}
\limsup_{m\to\infty}\sqrt{m}\;\cR_{2m}\Big[|\widehat{N}_{\vec{\mathbf{h}}}|^{\frac{1}{2}}, \bN_\vartheta\Big]<\infty.
\end{align*}
This together with (\ref{eq5-proof-th:main-adap-isotrop}) will complete the proof of the theorem.

\vskip0.1cm

$\mathbf{3^0.}\;$ Thus, let us prove (\ref{eq7-proof-th:main-adap-isotrop}). Introduce the following events
$$
\cB=\Big\{\big||\widehat{N}_{\vec{h}^*}|^{\frac{1}{2}}-\|f\|_2\big|> \ln(m)m^{-\frac{1}{2}}\Big\},\quad
\cC=\Big\{\big||\widehat{N}_{\vec{\mathrm{h}}_m}|^{\frac{1}{2}}\big|-\|f\|_2\big| > \ln(m)m^{-\frac{1}{2}}\Big\}
$$
and note that obviously $\bar{\cA}\subseteq\cB\cup\cC$.

$\mathbf{3^0a.}\;$ Applying Markov inequality we obtain
$$
\sup_{f\in\bN_\vartheta}\bP_f(\cC)\leq \tfrac{m \cR^2_{2m}[|\widehat{N}_{\vec{\mathrm{h}}_m}|^{\frac{1}{2}}, \bN_\vartheta]}{\ln^2(m)}
$$
and, therefore, in view of Theorem \ref{th:optim-adap-par-zone}
\begin{align}
\label{eq8-proof-th:main-adap-isotrop}
\limsup_{m\to\infty}\ln^2(m)\sup_{f\in\bN_\vartheta}\bP_f(\cC)<\infty,\quad\forall \theta\in\Theta^b_{\text{param}}.
\end{align}

$\mathbf{3^0b.}\;$ We have similarly to (\ref{eq9:proof-th:adaptive-anisotr}) that for any $\theta\in\Theta^b_{(r),\text{param}}$ and any
$f\in\bN_\vartheta$
$$
\mathrm{O}^*_f\leq c_1\big(\|f\|_2\sqrt{\mu_m}+\mu_m\big)\quad \Rightarrow\quad 2C_b(D)\mathrm{O}^*_f+m^{-1}\leq c_2\big(\|f\|_2\sqrt{\mu_m}+\mu_m\big),
$$
where, remind,  $\mu_m=m^{-1}\ln(m)$. Here $c_2$ is a constant depending only on $L, Q, b$ and the function $D$.
From now on we will assume that $m$ is sufficiently large in order to guarantee $2c_2<\ln(m)$.
We note that for any $f\in\bN_\vartheta\cap\bB_2(\sqrt{\mu_m})$
$$
2C_b(D)\mathrm{O}^*_f+m^{-1}\leq 2c_2\mu_m=\tfrac{2c_2\ln(m)}{m}\leq \tfrac{\ln^2(m)}{m}.
$$
We have in view of (\ref{eq4-proof-th:main-adap-isotrop}) for any $\theta\in\Theta^b_{(r),\text{param}}$ and $f\in\bN_\vartheta\cap\bB_2(\sqrt{\mu_m})$
\begin{align}
\label{eq9-proof-th:main-adap-isotrop}
\bP_f(\cB)&\leq \bP_f\Big\{\big|\widehat{N}_{\vec{h}^*}-\|f\|^2_2\big|\geq m^{-1}\ln^2(m)\Big\}
\nonumber\\*[2mm]
&\leq \bP_f\Big\{\big|\widehat{N}_{\vec{h}^*}-\|f\|^2_2\big|\geq 2C_b(D)\mathrm{O}^*_f+m^{-1}\Big\}\leq D_1[\ln(m)]^{d/2}m^{-4}.
\end{align}
To get the first inequality we have used that $|\sqrt{|a|}-\sqrt{|b|}|\leq \sqrt{|a-b|}$.
It remains to note that for any $\theta\in\Theta^b_{(r),\text{param}}$ and $f\in\bN_\vartheta\cap\bar{\bB}_2(\sqrt{\mu_m})$
$$
2C_b(D)\mathrm{O}^*_f+m^{-1}\leq 2c_2\|f\|_2\sqrt{\mu_m}\leq \|f\|_2\ln^{\frac{3}{2}}(m)m^{-\frac{1}{2}}
$$
and, therefore, for any $\theta\in\Theta^b_{(r),\text{param}}$ and $f\in\bN_\vartheta\cap\bar{\bB}_2(\sqrt{\mu_m})$
\begin{align}
\label{eq10-proof-th:main-adap-isotrop}
\bP_f(\cB)&\leq \bP_f\Big\{\big|\widehat{N}_{\vec{h}^*}-\|f\|^2_2\big|\geq \|f\|_2m^{-1}\ln^2(m)\Big\}
\nonumber\\*[2mm]
&\leq \bP_f\Big\{\big|\widehat{N}_{\vec{h}^*}-\|f\|^2_2\big|\geq 2C_b(D)\mathrm{O}^*_f+m^{-1}\Big\}\leq D_1[\ln(m)]^{d/2}m^{-4}.
\end{align}
To get the first inequality we have used that $|\sqrt{|a|}-\sqrt{|b|}|\leq |a-b|(|a|\vee |b|)^{-\frac{1}{2}}$.

We derive from (\ref{eq9-proof-th:main-adap-isotrop}) and (\ref{eq10-proof-th:main-adap-isotrop})
\begin{align}
\label{eq11-proof-th:main-adap-isotrop}
\limsup_{m\to\infty}\ln^2(m)\sup_{f\in\bN_\vartheta}\bP_f(\cB)=0,\quad\forall \theta\in\Theta^b_{(r),\text{param}}.
\end{align}
The required inequality (\ref{eq7-proof-th:main-adap-isotrop}) follows from (\ref{eq8-proof-th:main-adap-isotrop}) and (\ref{eq11-proof-th:main-adap-isotrop}). Theorem is proved.
\epr

\subsection{Proof of Proposition \ref{prop:bound-for-bias}.}
\label{sec:Proof of auxiliary results}
Let us recall the interpolation inequality.
Let $1\leq s_0<s<s_1\leq\infty$.
  If $f\in \bL_{s_0}\big(\bR^d\big)\cap \bL_{s_1}\big(\bR^d\big)$ then $f\in \bL_s\big(\bR^d\big)$ and
\[
 \big\|f\big\|_s\leq \big(\big\|f\big\|_{s_0}\big)^{\frac{(s_1-s)s_0}{(s_1-s_0)s}} \big(\big\|f\big\|_{s_1}\big)^{\frac{(s-s_0)s_1}{(s_1-s_0)s}}.
\]
We will also  need the following lemma.
\begin{lemma}
\label{lem:bias-via-emb}
Let $\vec{\beta}\in (0,b)^d$, $\vec{r}\in [1,\infty]^d$, $\vec{L}\in (0,\infty)^d$ and $Q>0$ be fixed.

\vskip0.1cm

1) There exists $1\leq C_1<\infty$ depending on $\cK$ and $b$ only such that for any $f\in\bN_{\vec{r},d}\big(\vec{\beta},\vec{L}\big)$
$$
\big\|b_{a,j}\big\|_{r_j}\leq C_1 L_ja^{\beta_j},\quad\forall a>0 .
$$

\vskip0.1cm

2) For any $f\in\bB_\infty(Q)$ and  $a>0$
$$
\big\|b_{a,j}\big\|_{p}\leq (1+\|\cK\|_1)Q^{1-1/p},\quad \forall p\in[1,\infty].
$$
\end{lemma}
Proof of the lemma is postponed to Section \ref{sec:Proof of technical lemmas}.

\paragraph{Proof of the proposition.}
Consider separately two cases.

\vskip0.1cm

$\mathbf{1^0.}\;$
Let $\vec{r}\in[2,\infty]^d$. Then, for any $j=1,\ldots,d,$ applying the interpolation inequality with $s_0=1$, $s=p$ and $s_1=r_j$,
and the both assertions of Lemma \ref{lem:bias-via-emb} for any $a>0$
 we obtain
\begin{align*}
 \big\|b_{a,j}\big\|_2&\leq (1+\|\cK\|_1)\big\|b_{a,j}\big\|_{r_j}^{\frac{r_j}{2(r_j-1)}}\leq (1+\|\cK\|_1)\big(C_1L_ja^{\beta_j}\big)^{\frac{p_j}{2}}
 \\
 &=(1+\|\cK\|_1)(C_1L_j)^{\frac{p_j}{2}}a^{\frac{\beta_jp_j}{2}},
\end{align*}
for any $f\in\bN_{\vec{r},d}\big(\vec{\beta},\vec{L}\big)$.
It yields for any $j=1,\ldots,d,$
$$
\sup_{a\leq h_j} \big\|b_{a,j}\big\|_2\leq C_1(1+\|\cK\|_1)(1\vee L_j)h_j^{\frac{\beta_jp_j}{2}}.
$$
and the assertion of the proposition follows in the considered case.

\vskip0.1cm

$\mathbf{2^0.}\;$  Let now  $\vec{r}\in[1,2]^d$. Applying the interpolation inequality with $s_0=r_j, s=2$ and $s_1=\infty$ and the both assertions of Lemma \ref{lem:bias-via-emb}
we get for  any $j=1,\ldots,d$ and  any  $f\in\bN_{\vec{r},d}\big(\vec{\beta},\vec{L}\big)\cap\bB_q(Q)$
$$
\big\|b_{a,j}\big\|_2\leq \big\|b_{a,j}\big\|_{r_j}^{\frac{r_j}{2}}\;Q^{\frac{2-r_j}{2}}
\leq (C_1L_j)^{\frac{p_j}{2}}a^{\frac{\beta_jr_j}{2}}=(C_1L_j)^{\frac{p_j}{2}}Q^{\frac{2-r_j}{2}}\;a^{\frac{\beta_jp_j}{2}}.
$$
Hence, taking into account that $p_j\leq 2$ since $r_j\leq 2$ we get
$$
\sup_{a\leq h_j} \big\|b_{a,j}\big\|_2\leq C_1\sqrt{Q\vee 1}(1\vee L_j)h_j^{\frac{\beta_jp_j}{2}},
$$
that completes the proof of the of proposition.
\epr

\subsection{Proof of technical lemmas.}
\label{sec:Proof of technical lemmas}

 \paragraph{Proof of Lemma \ref{lem:representation}.} The following result
 is the direct consequence of Lemma 1 in \cite{GL20b} in the case $p=2$.
 \begin{eqnarray}
  && \|f\|_2^2 =  2\int_{\bR^d} S_{\vec{h}}(x) f(x)\rd x-\int_{\bR^d} S_{\vec{h}}^2(x)\rd x
 + \int_{\bR^d}  B^2_{\vec{h}}(x)\rd x.
\label{eq:N-representation}
 \end{eqnarray}
We have obviously
\begin{align*}
 \int_{\bR^d} S_{_{\vec{h}}}(x) f(x)\rd x&=\int_{\bR^d}\int_{\bR^d}K_{\vec{h}}(x-y)f(x)f(y)\rd x\rd y
\\*[2mm]
\int_{\bR^d} S_{\vec{h}}^2(x)\rd x&=\int_{\bR^d}\int_{\bR^d}\bigg[\int_{\bR^d}K_{\vec{h}}(x-y)K_{\vec{h}}(x-z)\rd x\bigg]f(z)f(y)\rd z\rd y.
 \end{align*}
It remains to note that for any $y,z\in\bR^d$
$$
\int_{\bR^d}K_{\vec{h}}(x-y)K_{\vec{h}}(x-z)\rd x=\int_{\bR^d}K(v)V^{-1}_{\vec{h}}K\Big(v+\tfrac{y-z}{\vec{h}}\Big)\rd v
$$
and, therefore,
$$
 2\int_{\bR^d} S_{\vec{h}}(x) f(x)\rd x-\int_{\bR^d} S_{\vec{h}}^2(x)\rd x=\int_{\bR^d}\int_{\bR^d} T_{\vec{h}}(y-z)f(y)f(z)\rd y\rd z.
$$
The assertion of the lemma follows now from (\ref{eq:N-representation}).
\epr

\paragraph{Proof of Lemma \ref{lem:techniclal-proof-cor2}.}

Note that for any $x\in\bR^d$
\begin{align*}
B_{\vec{h},\vec{w}}(x)&=\big|S_{\vec{h}\vee\vec{w}}(x)-S_{\vec{w}}(x)\big|
\\*[2mm]
&=\bigg|\int_{\bR^d}K_{\vec{h}\vee \vec{w}}(u-x)\big[f(u)-f\big(x_\cJ,u_{\bar{\cJ}}\big)\big]\rd u
\\
&\hspace{2cm}-\int_{\bR^d}K_{\vec{w}}(u-x)\big[f(u)-f\big(x_\cJ,u_{\bar{\cJ}}\big)\big]\rd u\bigg|.
\end{align*}
Here we have used that $K_{\vec{h}\vee \vec{\eta},\bar{J}}\equiv K_{\vec{\eta},\bar{J}}$ in view of the definition of $J$ and $\int\cK=1$. Indeed, it yields
\begin{align*}
\int_{\bR^d}K_{\vec{h}\vee \vec{w}}(u-x)f\big(x_\cJ,u_{\bar{\cJ}}\big)\rd u&=\int_{\bR^{|\bar{\cJ}|}}K_{\vec{w},\bar{\cJ}}(u_{\bar{\cJ}}-x_{\bar{\cJ}})f\big(x_\cJ,u_{\bar{\cJ}}\big)\rd u_{\bar{\cJ}}
\\*[2mm]
&=\int_{\bR^d}K_{ \vec{w}}(u-x)f\big(x_\cJ,u_{\bar{\cJ}}\big)\rd u.
\end{align*}
Putting formally $f\big(x_\emptyset,u_{\bar{\emptyset}}\big)=f(u)$ we obtain
\begin{gather}
\label{eq1:proof-lem-tech}
f(u)-f\big(x_\cJ,u_{\bar{\cJ}}\big)=\sum_{l=1}^k f\big(x_{\cJ_{l-1}},u_{\bar{\cJ}_{l-1}}\big)-f\big(x_{\cJ_{l}},u_{\bar{\cJ}_{l}}\big),
\end{gather}
where $\cJ_0=\emptyset$. Note that for any $l=1,\ldots,k$
\begin{eqnarray*}
\int_{\bR}\cK_{h_l}(u_l-x_l)\big[f\big(x_{\cJ_{l-1}},u_{\bar{\cJ}_{l-1}}\big)-f\big(x_{\cJ_{l}},u_{\bar{\cJ}_{l}}\big)\big]\rd u_l=b_{h_l,l}\big(x_{\cJ_{l}},u_{\bar{\cJ}_{l}}\big)
\\*[2mm]
\int_{\bR}\cK_{w_l}(u_l-x_l)\big[f\big(x_{\cJ_{l-1}},u_{\bar{\cJ}_{l-1}}\big)-f\big(x_{\cJ_{l}},u_{\bar{\cJ}_{l}}\big)\big]\rd u_l=b_{w_l,l}\big(x_{\cJ_{l}},u_{\bar{\cJ}_{l}}\big).
\end{eqnarray*}
This yields together with (\ref{eq1:proof-lem-tech}) that for any $x\in\bR^d$
\begin{eqnarray}
\label{eq2:proof-lem-tech}
\int_{\bR^d}K_{\vec{h}\vee \vec{w}}(u-x)\big[f(u)-f\big(x_\cJ,u_{\bar{\cJ}}\big)\big]\rd u&=&\sum_{l=1}^k\big[K_{\vec{w}}\star b_{h_l,l}\big]_{\bar{\cJ_l}}(x);
\\*[2mm]
\int_{\bR^d}K_{\vec{w}}(u-x)\big[f(u)-f\big(x_\cJ,u_{\bar{\cJ}}\big)\big]\rd u&=&\sum_{l=1}^k\big[K_{\vec{w}}\star b_{w_l,l}\big]_{\bar{\cJ_l}}(x);
\nonumber
\end{eqnarray}
and the first assertion of the lemma is established.
It remains to note that the  second assertion follows from (\ref{eq2:proof-lem-tech}) if one chooses $\vec{w}=\vec{h}$.
\epr

\paragraph{Proof of Lemma \ref{lem:bias-via-emb}.}
 We obviously have for any $a>0$
\begin{eqnarray*}
\big|b_{a,j}(x)\big|=
\bigg|\int_{\bR} \cK(u)\big[f\big(x+ ua\mathbf{e}_j\big)-f(x)\big]\rd u \bigg|
=
\bigg|\int_{\bR} \cK(u)\big[\Delta_{ua,j}f(x)\big]\rd u  \bigg|.
\end{eqnarray*}
We have in view of the definition of the function $\cK$
\begin{eqnarray*}
\int_{\bR} \cK(u)\big[\Delta_{ua,j}f(x)\big]\rd u &=&
\int_{\bR} \sum_{j=1}^b (-1)^{j+1}\tbinom{b}{j}\tfrac{1}{j} D\big(\tfrac{u}{j}\big)\big[\Delta_{ua, j}f(x)\big]
\rd u
\nonumber
\\*[2mm]
&=& (-1)^{-b+1}\int_{\bR} D(z) \sum_{j=1}^b \tbinom{b}{j} (-1)^{j+b}\big[\Delta_{jza, j}f(x)\big]
\rd z
\nonumber\\*[2mm]
&=& (-1)^{-b+1} \int_{\bR} D(z) \big[[\Delta^b_{jza, j}\, f(x)\big]\rd z.
\label{eq:int-representation}
\end{eqnarray*}
The last equality follows from the definition of $b$-th order difference operator
(\ref{eq:Delta}).

Thus, for  any $ x\in\bR^d$ and $a>0$
\begin{equation}
\label{eq:B-h-j}
\big|b_{a,j}(x)\big|=\bigg| \int_{\bR} D(z) \big[[\Delta^b_{jza, j}\, f(x)\big]\rd z\bigg|
\end{equation}
Therefore, by the Minkowski inequality for integrals [see, e.g., \cite[Theorem~6.19]{Folland}]
 we obtain
\begin{eqnarray*}
&& \big\|b_{a,j}(\cdot,f)\big\|_{r_j}
\leq
  \int_{\bR} |D(z)|\big\| \Delta^b_{za, j}\, f\big\|_{r_j} \rd z .
\end{eqnarray*}
Since $f\in \bN_{\vec{r},d}\big(\vec{\beta},\vec{L}\big)$ we have for any $a>0$
$$
\big\| \Delta^b_{za, j}\, f\big\|_{r_j}\leq L_j |z|^{\beta_j}a^{\beta_j}.
$$
and, therefore,
$$
 \big\|b_{a,j}\big\|_{r_j}
\leq L_ja^{\beta_j}\int_{-b/2}^{b/2}|D(z)|\; |z|^{\beta_j}\rd z.
$$
Here we have used that $D$ is supported in $[-\frac{b}{2},\frac{b}{2}]$. It remains to note that
$$
\int_{-b/2}^{b/2}|D(z)|\; |z|^{\beta_j}\rd z\leq 2\|D\|_\infty (b/2\vee 1)^{b/2+1}
$$
and the first assertion of the lemma follows.
To prove the second assertion let us first note that for any $a>0$
\begin{align*}
 \big\|b_{a,j}\big\|_{1}&\leq 1+\int_{\bR^d}\int_{\bR}|\cK(u)|f\big(x+ua\mathbf{e}_j\big)\rd u\rd x=
 1+\int_{\bR}|\cK(u)|\int_{\bR^d}f\big(x+ua\mathbf{e}_j\big)\rd x\rd u
 \\
& =1+\|\cK\|_1,
\end{align*}
since $f$ is  a probability density. Next, for any $f\in\bB_\infty(Q)$
$$
 \big\|b_{a,j}\big\|_{\infty}\leq (1+\|\cK\|_1)Q.
$$
Applying the interpolation inequality with $s_0=1$, $s=p$ and $s_1=\infty$ we come to the third assertion. Lemma is proved.
\epr

\section{Proof of the results from Section \ref{sec:subsec-U-stat}.}
\label{sec:proof-U-stat}

\subsection{Basic result.}
Let $\vec{h}\in\mH^d$ be fixed. Introduce the following notations.
\begin{gather*}
I_{\vec{h}}(u)=\int_{\bR^d} U_{\vec{h}}(z-u)f(z)\rd z,\quad\; \cE_{\vec{h}}(f)=\int_{\bR^d} I_{\vec{h}}(u)f(u)\rd u;
\\*[2mm]
\xi_{\vec{h}}(x)=m^{-1}\sum_{j=1}^m \big(I_{\vec{h}}(x_j)-\cE_{\vec{h}}(f)\big),\quad x\in\bR^{dm}.
\end{gather*}
Set also for any  $x\in\bR^{dm}$ and $y\in\bR^d$
\begin{equation}
\label{eq:def-of-G}
G(y,x)=\tfrac{1}{m}\sum_{j=1}^m\big[U_{\vec{h}}\big(y-x_j\big)-I_{\vec{h}}(y)\big].
\end{equation}
We remark that
\begin{align}
\label{eq:U-stat-representation}
\widehat{U}_{\vec{h}}-\bE_f[\widehat{U}_{\vec{h}}]=\widehat{U}_{\vec{h}}-\cE_{\vec{h}}(f)
=\tfrac{1}{m}\sum_{i=1}^m G(Y_i,X^{(m)})+\xi_{\vec{h}}(Y^{(m)}).
\end{align}

Recall that $s=12q+2$ and set
\begin{align*}
G_{\vec{h}}(f)&=\Big[\tfrac{4s W_{\vec{h}}(f) J_{\vec{h}}(f)\ln(m)}{m}\Big]^{\frac{1}{2}}+\Big[\tfrac{16s\varpi_UJ_{\vec{h}}(f)\ln(m) }{m^2V_{\vec{h}}}\Big]^{\frac{1}{2}}+\tfrac{20sW_{\vec{h}}(f)\ln(m)}{m}
\\*[2mm]
&\quad + \Big[\tfrac{\Upsilon_{\vec{h}}(m)s^2\varpi_U \ln(m)}{m^2V_{\vec{h}}}\Big]
\\*[2mm]
\Lambda_q(U)&=2^{q+1}\varpi_U^q\Big(1 +2^{2q-1}\big[1+2^q\Gamma(q+1)\big]
\\*[2mm]
&\quad+2^{2q-2}\Gamma(q+1)\Big[2^{-1}\big(3\sqrt{12q+2}\big)^{q}
+\big(28(6q+1)\big)^{q}\Big]\Big).
\end{align*}

\begin{theorem}
\label{th:bound-for-func-G-final}

For any $m\geq 9$, $q\geq 1$, $f\in\bB_\infty\big(\alpha_m^{-1}\big)$  and $\vec{h}\in\mH_d$
\begin{equation*}
\bE_f\bigg\{\Big(\Big|m^{-1}\sum_{i=1}^m G\big(Y_i,X^{(m)}\big)\Big|-G_{\vec{h}}(f)\Big)_+^q\bigg\}\leq \Lambda_q(U)m^{-4q}.
\end{equation*}

\end{theorem}

\subsection{Proofs of Theorems \ref{th:bound-for-func-G-final}.}

\subsubsection{Auxiliary propositions.}

Recall that
$\mH^d_*=\big\{\vec{h}\in\mH^d:\; mV_{\vec{h}}\leq \alpha_m\big\}$,  $\bar{\mH}^d_*=\mH^d\setminus\mH^d_*$.
Set for any $s> 0$
\begin{eqnarray*}
\rho_{\vec{h}}(f)=\Big[\Big(\tfrac{2s\varpi_U W_{\vec{h}}(f)\ln(m)}{mV_{\vec{h}}}\Big)^{\frac{1}{2}}+\tfrac{4s\varpi_U \ln(m)}{3mV_{\vec{h}}}\Big]
\mathrm{1}_{\bar{\mH}^d_*}\big(\vec{h}\big)
+ \Big[\tfrac{2s\varpi_U \ln(m)}{mV_{\vec{h}}|\ln(mV_{\vec{h}})|}\Big]
\mathrm{1}_{\mH^d_*}\big(\vec{h}\big).
\end{eqnarray*}

\begin{proposition}
\label{prop:bound-for-func-G}

For any $m\geq 9$, $q,s\geq 1$, $f\in\bB_\infty\big(\alpha_m^{-1}\big)$ and $\vec{h}\in\mH^d$   one has
\begin{equation*}
\bP_f\Big(\cup_{i=1}^m\big\{|G(Y_i,X^{(m)})|> \rho_{\vec{h}}(f)\big\}\Big)\leq 2m^{-\frac{s}{2}+1}.
\end{equation*}

\end{proposition}

\noindent Set for any $x\in\bR^{dm}$
$$
\upsilon^2(x)=\bE_f\big\{|G(Y_1,x)|^2\big\}.
$$
and introduce the following quantity:
$$
c_{\vec{h}}(f)=\Big[\tfrac{4\varpi_U J_{\vec{h}}(f)}{mV_{\vec{h}}}\Big]^{\frac{1}{2}}+\Big[\tfrac{8s\ln(m)}{m}\Big]^{\frac{1}{2}}W_{\vec{h}}(f)
+\Big[\tfrac{11s^2\varpi_U W_{\vec{h}}(f)\ln^2(m)}{m^2V_{\vec{h}}}\Big]^{\frac{1}{2}}.
$$

\begin{proposition}
\label{prop:bound-for-second-term2}
For any $m\geq 9$, $q,s\geq 1$, $f\in\bB_\infty\big(\alpha_m^{-1}\big)$ and $\vec{h}\in\mH^d$   one has
\begin{equation*}
\label{eq1:random-bound-for-second-term2}
\bE_f\big\{\big(\upsilon\big(X^{(m)}\big)-c_{\vec{h}}(f)\big)_+^q\big\}\leq \Gamma(q+1)\big[6\varpi_U\big]^{q}m^{\frac{3q}{2}-s}.
\end{equation*}

\end{proposition}
\noindent The proofs of both propositions are given in Section \ref{sec:Proof of auxiliary results}.

\subsubsection{Probabilistic background.}

\begin{lemma}
\label{lem:lb-deviation-new}
Let $\e_j,j=1,\ldots, k$ be i.i.d. \textbf{nonnegative} random variables defined on a measurable space equipped with a probability measure $\mathrm{P}$. Assume that $\mathrm{P}\{\e_1\leq T\}=1$ for some $T>0$.

Set $\cE_k=\sum_{j=1}^k\e_j$ and let $\alpha=\mathrm{E}(\e_1)$, where $\mathrm{E}$ is the mathematical expectation w.r.t. $\mathrm{P}$.
Then for any $z\geq  k\alpha$ one has
\[
 \mathrm{P}\big(\cE_k \geq z\big) \leq \big(\tfrac{k\alpha}{z}\big)^{\frac{z}{T}} e^{\frac{z-k\alpha}{T}}.
\]
\end{lemma}
\noindent The authors are sure that this result is known but they cannot find the exact reference. We put the elementary
proof of the lemma  into Section \ref{sec:Proof of technical lemmas}. In the case when $\cE_m$ is a binomial r.v. with parameters $m,p$, this result with $T=1$ and  $\alpha=p$ can be found in \cite{Bouch-Lugosi-Massart}.

\begin{lemma}
\label{lem:bernstein-new}
Let $\e_j,j=1,\ldots, m$ be i.i.d. \textbf{centered} random variables defined on measurable space equipped with a probability measure $\mathrm{P}$. Assume that $\mathrm{P}\{|\e_1|\leq T\}=1$ for some $T>0$.
Set $\cE_m=m^{-1}\sum_{k=1}^m\e_k$ and let $v^2=\mathrm{E}(\e^2_1)$.
Then for any $m\geq 3$ and $s>0$ one has
\[
 \mathrm{P}\Big\{|\cE_m| \geq \big[\tfrac{2sv^2\ln(m)}{m}\big]^{\frac{1}{2}}+\tfrac{2sT\ln(m)}{3m}\Big\} \leq 2m^{-s}.
\]
Additionally for any $p\geq 1$, $s\geq 1$ and $m\geq 3$
$$
 \mathrm{E} \Big\{\Big(|\cE_m|-\big[\tfrac{4sv^2\ln(m)}{m}\big]^{\frac{1}{2}}-\tfrac{4sT\ln(m)}{3m}\Big)_+^p\Big\} \leq 2^{p+1}\Gamma(p+1)T^p m^{-s-\frac{p}{2}},
$$
where $\Gamma$ is gamma-function.
\end{lemma}
\noindent Both assertions of the lemma are direct consequences of the Bernstein inequality and their proofs are omitted.

\smallskip

\noindent Let $\cS$ be a countable set of functions $S:\cN\to\bR$ defined on some set $\cN$. Let $\eta_i, i=1,\ldots m,$ be i.i.d $\cN$-valued random variables defined on some measurable space endowed with probability measure $\mathrm{P}$.
Suppose
that $\mathrm{E} \{S(\eta_1)\}=0$, $\|S\|_\infty\leq T<\infty$,  $\forall S\in\cS$, and  put
\begin{eqnarray*}
Z=\sup_{S\in\cS}\sum_{i=1}^{m}S\big(\eta_i\big),\;\quad
\varsigma^2=\sup_{S\in\cS}\mathrm{E} [S(\eta_1)]^2,
\end{eqnarray*}
where as previously $\mathrm{E}$ is the mathematical expectation w.r.t. $\mathrm{P}$.

\begin{lemma}[\cite{bousquet}]
\label{bousquet}
For any $x\geq 0$
$$
\mathrm{P}\left\{Z-\mathrm{E}\{Z\}\geq x\right\}\leq
e^{-\frac{x^2}{2m\varsigma^2 + 4T\bE\{Z\} + \frac{2}{3}T x}}.
$$
\end{lemma}

\noindent The following useful bound can be easily obtained by the integration of Bousquet inequality. For any $t\geq 1$, $p\geq 1$ and $m\in\bN^*$
\begin{eqnarray}
\label{eq0:Bousquet-expectation}
&&\mathrm{E}\Big\{Z- \mathrm{E}\{Z\}-\sqrt{(2m\varsigma^2+4T\bE\{Z\}) t}-(2/3)Tt\Big\}_+^{p}
\nonumber\\
&&\qquad\leq \Gamma(p+1)\Big(\sqrt{2m\varsigma^2 + 4T\bE\{Z\}}+(2/3)T\Big)^{p}e^{-t}.
\end{eqnarray}
The latter inequality can be simplified.  First we note that
$$
\sqrt{(2m\varsigma^2+4T\bE\{Z\}) t}\leq \varsigma\sqrt{2mt}+Tt +\bE\{Z\}.
$$
Next, for any $m\geq 3$
$$
\sqrt{2m\varsigma^2 + 4T\bE\{Z\}}+(2/3)T\leq T\sqrt{m}\Big[\sqrt{6}+\tfrac{2}{3\sqrt{3}}\Big]\leq 3T\sqrt{m}.
$$
Thus, we deduce from (\ref{eq0:Bousquet-expectation}) for any $t\geq 1$
\begin{eqnarray}
\label{eq:Bousquet-expectation}
\mathrm{E}\Big\{Z- 2\mathrm{E}\{Z\}-\varsigma\sqrt{2mt}-(5/3)Tt\Big\}_+^{p}
\leq 3^{p}\Gamma(p+1)T^{p}m^{\frac{p}{2}}e^{-t}.
\end{eqnarray}

\subsubsection{Proof of Theorem \ref{th:bound-for-func-G-final}.}

\quad $\mathbf{1^0.}\;$ Recall that $s=12q+2$ and
set  $C_{\vec{h}}(f)=C^*_{\vec{h}}(f)+d_{\vec{h}}(f)$, where
\begin{align*}
C^*_{\vec{h}}(f)&=\big[\tfrac{4s\ln(m)}{m}\big]^{\frac{1}{2}}c_{\vec{h}}(f)+\tfrac{8s\rho_{\vec{h}}(f)\ln(m)}{3m};
\nonumber\\*[2mm]
 d_{\vec{h}}(f)&=\Big[\tfrac{4s W_{h}(f) J_{h}(f)\ln(m)}{m}\Big]^{1/2}+\tfrac{8sW_{\vec{h}}(f)\ln(m)}{3m}.
\end{align*}
Our first goal is to show that
\begin{equation}
\label{eq1:proof-th1}
C_{\vec{h}}(f)\leq G_{\vec{h}}(f),\quad\forall \vec{h}\in\mH^d,\; f\in\bB_\infty(Q).
\end{equation}
Indeed, the simple algebra shows that
\begin{align*}
C_{\vec{h}}(f)&\leq \Big[\tfrac{4s W_{h}(f) J_{h}(f)\ln(m)}{m}\Big]^{1/2}+\Big[\tfrac{16s\varpi_UJ_{\vec{h}}(f)\ln(m) }{m^2V_{\vec{h}}}\Big]^{\frac{1}{2}}+
\\*[2mm]
&\quad +\Big[\tfrac{121 s^3\varpi_U W_{\vec{h}}(f)\ln^3(m)}{m^3V_{\vec{h}}}\Big]^{\frac{1}{2}}+\tfrac{14sW_{\vec{h}}(f)\ln(m)}{m}
\\*[2mm]
&\quad+\Big[\tfrac{32s^2\varpi_U \ln^2(m)}{9m^2V_{\vec{h}}}\Big]
\mathrm{1}_{\bar{\mH}^d_*}\big(\vec{h}\big) +\Big[\tfrac{16s^2\varpi_U \ln^2(m)}{3m^2V_{\vec{h}}|\ln(mV_{\vec{h}})|}\Big]
\mathrm{1}_{\mH^d_*}\big(\vec{h}\big)
\end{align*}
Consider separately two cases.
First let $\vec{h}\in\bar{\mH}^d_*$. Then using $\sqrt{ab}\leq (a+b)/2, a,b\geq 0$, we get
$$
\Big[\tfrac{121 s^3\varpi_U W_{\vec{h}}(f)\ln^3(m)}{m^3V_{\vec{h}}}\Big]^{\frac{1}{2}}\leq \tfrac{5.5 sW_{\vec{h}}(f)\ln(m)}{m}+\tfrac{5.5 s^2\varpi_U \ln^2(m)}{m^2V_{\vec{h}}}
$$
and, therefore,
\begin{align}
\label{eq2:proof-th1}
C_{\vec{h}}(f)&\leq \Big[\tfrac{4s W_{h}(f) J_{h}(f)\ln(m)}{m}\Big]^{1/2}+\Big[\tfrac{16s\varpi_UJ_{\vec{h}}(f)\ln(m) }{m^2V_{\vec{h}}}\Big]^{\frac{1}{2}}+\tfrac{20sW_{\vec{h}}(f)\ln(m)}{m}
\nonumber\\*[2mm]
&\quad+\Big[\tfrac{10s^2\varpi_U \ln^2(m)}{m^2V_{\vec{h}}}\Big]
\mathrm{1}_{\bar{\mH}^d_*}\big(\vec{h}\big).
\end{align}
Since $W_{\vec{h}}(f)\leq \varpi_U \alpha_m^{-1}$ for any $\vec{h}\in\mH^d$ and $f\in\bB_\infty\big(\alpha^{-1}_m\big)$ we have for any $\vec{h}\in\mH^d_*$
\begin{align*}
&\Big[\tfrac{121 s^3\varpi_U W_{\vec{h}}(f)\ln^3(m)}{m^3V_{\vec{h}}}\Big]^{\frac{1}{2}}\leq [11s^2\varpi_U\ln^2(m)]\big[\tfrac{1}{\alpha_m\ln(m)}\big]^{\frac{1}{2}}m^{-1}[mV_{\vec{h}}]^{-\frac{1}{2}}
\\*[2mm]
&\leq [11s^2\varpi_U\ln^2(m)]m^{-1}[mV_{\vec{h}}]^{-\frac{1}{2}}\leq \tfrac{11s^2\varpi_U \ln^2(m)}{m^2V_{\vec{h}}|\ln(mV_{\vec{h}})|}
\big[(mV_{\vec{h}})\ln^2(mV_{\vec{h}})\big]^{1/2}
\\*[2mm]
&\leq \tfrac{11s^2\varpi_U \ln^2(m)}{m^2V_{\vec{h}}|\ln(mV_{\vec{h}})|}.
\end{align*}
Here we used that $Q\leq \ln(m)$,
$\alpha_m\ln(m)\geq 1$ and  that $x\mapsto x\ln^2(x)$ is monotonically increasing on $[0,e^{-2}]$
and that $mV_{\vec{h}}\leq \alpha^2_m\leq 1/9$.
\begin{align}
\label{eq3:proof-th1}
C_{\vec{h}}(f)&\leq \Big[\tfrac{4s W_{h}(f) J_{h}(f)\ln(m)}{m}\Big]^{1/2}+\Big[\tfrac{16s\varpi_UJ_{\vec{h}}(f)\ln(m) }{m^2V_{\vec{h}}}\Big]^{\frac{1}{2}}+\tfrac{14sW_{\vec{h}}(f)\ln(m)}{m}
\nonumber\\*[2mm]
&\quad+\Big[\tfrac{17s^2\varpi_U \ln^2(m)}{m^2V_{\vec{h}}|\ln(mV_{\vec{h}})|}\Big]
\mathrm{1}_{\mH^d_*}\big(\vec{h}\big).
\end{align}
Thus, (\ref{eq1:proof-th1}) follows from (\ref{eq2:proof-th1}) and (\ref{eq3:proof-th1}).

$\mathbf{2^0.}\;$ Recall that
$$
\mathfrak{p}_f(x)=\prod_{i=1}^m f(x_i), x=(x_1,\ldots,x_m)\in\bR^{dm}.
$$
and
introduce the  event
$$
\cA=\cap_{i=1}^m\big\{|G(Y_i,X^{(m)})|\leq \varpi_U\;\rho_{\vec{h}}(f)\big\}.
$$
As usual the complimentary event will be denoted by $\bar{\cA}$. Note that
\begin{align}
\label{eq11-proof-prop1-new00}
\sup_{x\in\bR^{dm},  y\in\bR^d}|G(y,x)|\leq 2\varpi_U V^{-1}_{\vec{h}}\leq 2\varpi_Um^2, \quad\forall \vec{h}\in\mH^d.
\end{align}
We obviously have using (\ref{eq1:proof-th1}) and (\ref{eq11-proof-prop1-new00})
\begin{align}
\label{eq8-proof-prop1-new00}
\bE_f&\bigg\{\Big(\Big|m^{-1}\sum_{i=1}^m G(Y_i,X^{(m)})\Big|-G_{\vec{h}}(f)\Big)_+^q\bigg\}
\nonumber\\*[2mm]
&\leq \bE_f\bigg\{\Big(\Big|m^{-1}\sum_{i=1}^m G(Y_i,X^{(m)})\Big|-C_{\vec{h}}(f)\Big)_+^q\mathrm{1}_{\cA}\bigg\}+\big(2\varpi_Um^2\big)^{q}\bP_f\big(\bar{\cA}\big)
\nonumber\\*[2mm]
&=\bE_f\bigg\{\Big(\Big|m^{-1}\sum_{i=1}^m G(Y_i,X^{(m)})\Big|-C_{\vec{h}}(f)\Big)_+^q\mathrm{1}_{\cA}\bigg\}
\nonumber\\*[2mm]
&\quad+2^{q+1}\varpi_U^q m^{-\frac{s}{2}+1+2q}
\end{align}
in view of Proposition  \ref{prop:bound-for-func-G}. We obviously have
\begin{align*}
\Big\{&\Big|m^{-1}\sum_{i=1}^m G(Y_i,X^{(m)})\Big|- C_{\vec{h}}(f)\Big\}_+\mathrm{1}_{\cA}
\nonumber\\*[2mm]
&\subseteq \Big\{\Big|m^{-1}\sum_{i=1}^m G(Y_i,X^{(m)})\mathrm{1}\big(|G(Y_i,X^{(m)})|\leq \rho_{\vec{h}}(f)\big)\Big|- C_{\vec{h}}(f) \Big\}_+.
\end{align*}
Setting for any $i=1,\ldots,m$, $y\in\bR^d$ and $x\in\bR^{dm}$
$$
\cG(y,x)=G(y,x)\mathrm{1}\big(|G(y,x)|\leq \rho_{\vec{h}}(f)\big).
$$
we have in view of  (\ref{eq8-proof-prop1-new00})
\begin{align}
\label{eq9-proof-prop1-new00}
\bE&_f\bigg\{\Big(\Big|m^{-1}\sum_{i=1}^m G(Y_i,X^{(m)})\Big|-C_{\vec{h}}(f)\Big)_+^q\bigg\}
\nonumber\\*[2mm]
&\leq
\bE_f\bigg\{\Big(\Big|m^{-1}\sum_{i=1}^m \cG(Y_i,X^{(m)})\Big|-C_{\vec{h}}(f)\Big)_+^q\bigg\}
\nonumber\\*[2mm]
&\quad+2^{q+1}\varpi_U^q m^{-\frac{s}{2}+1+2q}
=:E_m+2^{q+1}\varpi_U^q m^{-\frac{s}{2}+1+2q}.
\end{align}

$\mathbf{3^0.}\;$ We have in view of the independence of $X^{(m)}$ and $Y^{(m)}$
\begin{align}
\label{eq10-proof-prop1-new00}
E_m=\int_{\bR^{dm}}\bE_f\bigg\{\Big(\Big|m^{-1}\sum_{i=1}^m \cG(Y_i,x)\Big|-C_{\vec{h}}(f)\Big)_+^q\bigg\}\mathfrak{p}_f(x)\rd x.
\end{align}
Set for any $x\in\bR^{dm}$
$$
Z(x)=\bE_f\big\{\cG(Y_1,x)\big\},\quad Z_*(x)=\bE_f\big\{G(Y_1,x)\big\}-Z(x)=\xi_{\vec{h}}(x)-Z(x).
$$
Using the trivial inequalitis $(a+b)_+\leq (a)_++(b)_+$, $a,b\in\bR$ and $(c+d)^q\leq 2^{q-1}[c^q+d^q]$, $c,d\geq 0$, we obtain that
\begin{align}
\label{eq13-proof-prop1-new00}
E_m&\leq2^{q-1}\int_{\bR^{dm}}\bE_f\bigg\{\Big(\Big|m^{-1}\sum_{i=1}^m \cG(Y_i,x)-Z(x)\Big|-C^*_{\vec{h}}(f)\Big)_+^q\bigg\}\mathfrak{p}_f(x)\rd x&
\nonumber\\*[2mm]
+&2^{2q-1}\int_{\bR^{dm}}|Z_*(x)|^{q}\mathfrak{p}_f(x)\rd x+2^{2q-1}\int_{\bR^{dm}}\big\{\big|\xi_{\vec{h}}(x)\big|-d_{\vec{h}}(f)\big\}_+^q\mathfrak{p}_f(x)\rd x
\nonumber\\*[2mm]
&=:2^{q-1}E_m^*+2^{2q-1}\int_{\bR^{dm}}|Z_*(x)|^{q}\mathfrak{p}_f(x)\rd x
\nonumber\\*[2mm]
&\quad+2^{2q-1}\bE_f\big\{\big(\big|\xi_{\vec{h}}\big(X^{(m)}\big)\big|-d_{\vec{h}}(f)\big)\big\}_+^q.
\end{align}
Some remarks are in order. First we note that for any $x\in\bR^{dm}$
\begin{align*}
\big|Z_*(x)\big|&=\big|\bE_f\big\{G(Y_1,x)\mathrm{1}\big(|G(Y_1,x)|> \rho_{\vec{h}}(f)\big)\big\}\big|.
\end{align*}
By Jensen inequality we have for any $q\geq 1$ and $x\in\bR^{dm}$
\begin{align*}
\big|Z_*(x)\big|^q&\leq\bE_f\big\{\big|G(Y_1,x)\big|^q\mathrm{1}\big(|G(Y_1,x)|> \rho_{\vec{h}}(f)\big)\big\}.
\end{align*}
Thus, we obtain using  first Fubini theorem and next the bound (\ref{eq11-proof-prop1-new00})
\begin{align*}
\int_{\bR^{dm}}&|Z_*(x)|^{q}\mathfrak{p}_f(x)\rd x
\nonumber\\*[2mm]
&\leq \int_{\bR^{dm}}\bE_f\big\{\big|G(Y_1,x)\big|^q\mathrm{1}\big(|G(Y_1,x)|> \rho_{\vec{h}}(f)\big)\big\}\mathfrak{p}_f(x)\rd x
\nonumber\\*[2mm]
&=\int_{\bR^{dm}}\Big[\int_{\bR^d}\big|G(y,x)\big|^q\mathrm{1}\big(|G(y,x)|> \rho_{\vec{h}}(f)\big)f(y)\rd y\Big]\mathfrak{p}_f(x)\rd x
\nonumber\\*[2mm]
&=\bE_f\Big\{\big|G(Y_1,X^{(m)})\big|^q\mathrm{1}\big(|G(Y_1,X^{(m)})|> \rho_{\vec{h}}(f)\big)\Big\}
\nonumber\\*[2mm]
&\leq \big(2\varpi_Um^2\big)^q\bP_f\big\{|G(Y_1,X^{(m)})|> \rho_{\vec{h}}(f)\big\}.
\end{align*}
Using Proposition \ref{prop:bound-for-func-G} we finally get
$$
\int_{\bR^{dm}}|Z_*(x)|^{q}\mathfrak{p}_f(x)\rd x\leq 2^{q+1}\varpi_U^q m^{-\frac{s}{2}+2q+1}
$$
and, therefore, we deduce from (\ref{eq13-proof-prop1-new00}) that
\begin{align*}
E_m\leq 2^{q-1}E_m^*+2^{3q}\varpi_U^q m^{-\frac{s}{2}+2q+1}+2^{2q-1}\bE_f\big\{\big(\big|\xi_{\vec{h}}\big(X^{(m)}\big)\big|-d_{\vec{h}}(f)\big)\big\}_+^q.
\end{align*}
We note that $\xi_{\vec{h}}\big(X^{(m)}\big)$ is the arithmetic mean of the independent, centered and bounded by $2W_{\vec{h}}\leq 2\varpi_U m^2$ random variables.
Hence, applying second assertion of Lemma \ref{lem:bernstein-new} with $T=2W_{\vec{h}}$, $p=q$ and
\begin{align}
\label{eq1:new-new-new}
v^2=\int_{\bR^d}I^2_{\vec{h}}(z)f(z)\rd z\leq  W_{h}(f)J_{h}(f),
\end{align}
we obtain
$$
\bE_f\big\{\big(\big|\xi_{\vec{h}}\big(X^{(m)}\big)\big|-d_{\vec{h}}(f)\big)\big\}_+^q\leq 2^{2q+1}\Gamma(q+1)\varpi^q_U m^{-s+\frac{3q}{2}}.
$$
Thus, we get
\begin{align}
\label{eq12-proof-prop1-new00}
E_m\leq 2^{q-1}E_m^*+2^{3q}\varpi_U^q\big[1+2^q\Gamma(q+1)\big]m^{-\frac{s}{2}+2q+1}.
\end{align}

$\mathbf{4^0.}\;$ Set now for any $x\in\bR^{dm}$
$$
\zeta(x)=m^{-1}\sum_{i=1}^m \cG(Y_i,x)-Z(x).
$$
We note that $\zeta(x)$ is the arithmetic mean of the independent, centered and bounded by $2\rho_{\vec{h}}(f)$ random variables.
Hence, applying second assertion of Lemma \ref{lem:bernstein-new} with $T=2\rho_{\vec{h}}(f)$, $p=q$ and $v^2=\varsigma^2(x):=\bE_f\{\cG^2(Y_1,x)\}$
we obtain for any $x\in\bR^{dm}$
 \begin{equation}
 \label{eq14-proof-prop1-new00}
 \bE_f \Big\{\Big(|\zeta(x)|-\big[\tfrac{4s\varsigma^2(x)\ln(m)}{m}\big]^{\frac{1}{2}}-\tfrac{8s\rho_{\vec{h}}(f)\ln(m)}{3m}\Big)_+^q\Big\} \leq \lambda_q\big[2\rho_{\vec{h}}(f)\big]^{q}m^{-s-\frac{q}{2}},
\end{equation}
where we put for brevity $\lambda_q=2^{q+1}\Gamma(q+1)$.
We have
\begin{align*}
\label{eq15-proof-prop1-new00}
\varsigma^2(x)=\bE_f\{\cG^2(Y_1,x)\}\leq \bE_f\{G^2(Y_1,x)\}=:\upsilon^2(x).
\end{align*}
It yields together with  (\ref{eq14-proof-prop1-new00})
\begin{equation}
 \label{eq16-proof-prop1-new00}
 \bE_f \big\{\big(|\zeta(x)|-\Delta_{\vec{h}}(x)\big)_+^q\big\} \leq \lambda_q\big[2\rho_{\vec{h}}(f)\big]^{q}m^{-s-\frac{q}{2}},\;\quad\forall x\in\bR^d.
\end{equation}
where we have put
$$
\Delta_{\vec{h}}(x)=\big[\tfrac{4s\upsilon^2(x)\ln(m)}{m}\big]^{\frac{1}{2}}+\tfrac{8s\rho_{\vec{h}}(f)\ln(m)}{3m}.
$$
We have in view of (\ref{eq16-proof-prop1-new00}) for any $x\in\bR^{dm}$
\begin{align}
\label{eq17-proof-prop1-new00}
&\bE_f\bigg\{\Big(\Big|m^{-1}\sum_{i=1}^m \cG(Y_i,x)-Z(x)\Big|-C_{\vec{h}}(f)\Big)_+^q\bigg\}
\nonumber\\*[2mm]
&\leq 2^{q-1}\bE_f\big\{\big|\zeta(x)\big|-\Delta_{\vec{h}}(x)\big)_+^q\big\}+2^{q-1}\bE_f\big\{\big(\Delta_{\vec{h}}(x)-C_{\vec{h}}(f)\big)_+^q\big\}
\nonumber\\*[2mm]
&\leq 2^{q-1}\bE_f\big\{\big(\Delta_{\vec{h}}(x)-C_{\vec{h}}(f)\big)_+^q\big\}+2^{2q-1}\lambda_q\big[\rho_{\vec{h}}(f)\big]^{q}m^{-s-\frac{q}{2}}
\nonumber\\*[2mm]
&=2^{q-1}\big[\tfrac{4s\ln(m)}{m}\big]^{\frac{q}{2}}\big(\upsilon(x)-c_{\vec{h}}(f)\big)_+^q
+2^{2q-1}\lambda_q\big[\rho_{\vec{h}}(f)\big]^{q}m^{-s-\frac{q}{2}}.
\end{align}
The last equality follows from the fact that $\Delta_{\vec{h}}(x)$ does not depend on $Y^{(m)}$.
Thus, we obtain
\begin{align*}
E^*_m&\leq 2^{2q-1}\big[\tfrac{s\ln(m)}{m}\big]^{\frac{q}{2}}\int_{\bR^{dm}}\big(\upsilon(x)-c_{\vec{h}}(f)\big)_+^{q}\mathfrak{p}_f(x)\rd x
\nonumber\\*[2mm]
&\quad+2^{2q-1}\lambda_q\big[\rho_{\vec{h}}(f)\big]^{q}m^{-s-\frac{q}{2}}
\nonumber\\*[2mm]
&=2^{2q-1}\big[\tfrac{s\ln(m)}{m}\big]^{\frac{q}{2}}\bE_f\big\{\big(\upsilon\big(X^{(m)}\big)-c_{\vec{h}}(f)\big)_+^{q}\big\}
\nonumber\\*[2mm]
&\quad
+2^{2q-1}\lambda_q\big[\rho_{\vec{h}}(f)\big]^{q}m^{-s-\frac{q}{2}}.
\end{align*}
Applying  Proposition \ref{prop:bound-for-second-term2} we finally obtain
\begin{align*}
E^*_m&\leq 2^{-1}\Gamma(q+1)\big[12\sqrt{s}\varpi_U\big]^{q}m^{-s+\frac{3q}{2}}
\nonumber\\*[2mm]
&\quad+2^{3q}\Gamma(q+1)\big[\rho_{\vec{h}}(f)\big]^{q}m^{-s-\frac{q}{2}}.
\end{align*}
Here we have also used that $m^{-1}\ln(m)\leq 1/4$ for any $m\geq 9$.

It remains to note that
$
\rho_{\vec{h}}(f)\leq 7s\varpi_Um
$
for any $\vec{h}\in\mH^d$ and we get
\begin{align}
\label{eq1800-proof-prop1-new00}
E^*_m\leq (4\varpi_U)^q\Gamma(q+1)\Big[2^{-1}\big(3\sqrt{s}\big)^{q}
+\big(14s\big)^{q}\Big]m^{-s+\frac{3q}{2}}.
\end{align}
The assertion of the theorem follows from (\ref{eq9-proof-prop1-new00}), (\ref{eq12-proof-prop1-new00}) and (\ref{eq1800-proof-prop1-new00}).
\epr

\subsection{Proofs of Theorems \ref{th:non-random-upper-function}, \ref{th:random-upper-function} and Corollaries \ref{cor1:th2}, \ref{cor2:th2}.}

\subsubsection{Proof of Theorem \ref{th:non-random-upper-function}.} First, we remark that $U_{\vec{h}}(f)> G_{\vec{h}}(f)+d_{\vec{h}}(f)$. Next
in view of the representation of $\widehat{U}_{\vec{h}}$ given in (\ref{eq:U-stat-representation}) we obviously have
\begin{align}
\label{eq1:proof-th2}
\bE_f&\Big\{\Big(\big|\widehat{U}_{\vec{h}}-\bE_{f}[\widehat{U}_{\vec{h}}]\big|-U_{\vec{h}}(f)\Big)_+^q\Big\}
\nonumber\\*[2mm]
&\quad \leq 2^{q-1}\bE_f\bigg\{\Big(\Big|m^{-1}\sum_{i=1}^m G(Y_i,X^{(m)})\Big|-G_{\vec{h}}(f)\Big)_+^q\bigg\}
\nonumber\\*[2mm]
&\qquad + 2^{q-1}\bE_f\Big\{\Big(\big|\xi_{\vec{h}}\big(Y^{(m)}\big)\big|-d_{\vec{h}}(f)\Big)_+^q\Big\}.
\end{align}
Note that  $\xi_{\vec{h}}(Y^{(m)})$ is the average of $\varepsilon_{i}:=I_{\vec{h}}(Y_i)-\cE_{\vec{h}}(f)$, which are i.i.d. centered random variables.
Moreover  $\|I\|_\infty\leq W_{\vec{h}}(f)$ and, see (\ref{eq1:new-new-new}),
$$
\bE_{f}(\varepsilon_1^2)\leq \bE_{f}\big\{I^2_{\vec{h}}(Y_1)\big\}\leq  W_{h}(f)J_{h}(f).
$$
Note also that $W_{\vec{h}}(f)\le \varpi_{U}V_{\vec{h}}^{-1}\leq \varpi_{U}m^2$, verifying for any  $\vec{h}\in\mH^d$.
Applying the second assertion of Lemma \ref{lem:bernstein-new} with $T=2W_{\vec{h}}$ and $s=12q+2$
\begin{align}
\label{eq2:proof-th2}
\bE_f\big\{\big(|\xi_{\vec{h}}(Y^{(m)})|-d_{\vec{h}}(f)\big)_+^q\big\}
\leq 2\Gamma(q+1)\big[4\varpi_{U}\big]^q m^{-10q-2-\frac{q}{2}}.
\end{align}
Thus, we deduce from
(\ref{eq1:proof-th2}), (\ref{eq2:proof-th2}) and the statement of Theorem \ref{th:bound-for-func-G-final} that
\begin{align}
\label{eq4:proof-th2}
\bE_f\bigg\{&\Big(\big|\widehat{U}_{\vec{h}}-\bE_{f}[\widehat{U}_{\vec{h}}]\big|-U_{\vec{h}}(f)\Big)_+^q\bigg\}\leq 2^{q-1}\Lambda_q(U)m^{-4q}
\nonumber\\*[2mm]
&\quad+\Gamma(q+1)\big[8\varpi_{U}\big]^q m^{-12q-1-\frac{q}{2}}\leq \Lambda^*_q(U)m^{-4q}
\end{align}
for any $\vec{h}\in \mH^d$. Noting that
\begin{align*}
\bE_f\bigg\{&\sup_{\vec{h}\in\cH^d_m}\Big(\big|\widehat{U}_{\vec{h}}-\bE_{f}[\widehat{U}_{\vec{h}}]\big|-U_{\vec{h}}(f)\Big)_+^q\bigg\}
\nonumber\\*[2mm]
&\leq \sum_{\vec{h}\in\cH^d_m}\bE_f\bigg\{\sup_{\vec{h}\in\mH^d}\Big(\big|\widehat{U}_{\vec{h}}-\bE_{f}[\widehat{U}_{\vec{h}}]\big|-U_{\vec{h}}(f)\Big)_+^q\bigg\}
\end{align*}
and  that the cardinality of $\cH^d_m$ is less than $[2\ln(m)]^d$, we deduce the assertion of the theorem from  (\ref{eq4:proof-th2}) since $\cH^d_m\subset\mH^d$.
\epr

\subsubsection{Proof of Theorem \ref{th:random-upper-function}.}
Throughout this section $A,A_1,\ldots, B,B_1,\ldots$ denote the constants which may depend on $q$ and $U$ only.
We well need the following result. Set
$$
\widetilde{\cW}_{\vec{h}}(f)=\sup_{\vec{k}\in\bZ^d}V_{\vec{h}}^{-1}\int_{\bR^d}\mathrm{1}_{\Pi^*_{\vec{k}}}\big(\tfrac{y}{\vec{h}}\big)f(y)\rd y,\quad
\widetilde{\cW}_{\vec{h}}=\sup_{\vec{k}\in\bZ^d}\tfrac{1}{mV_{\vec{h}}}\sum_{j=1}^m\mathrm{1}_{\Pi^*_{\vec{k}}}\big(\tfrac{X_j}{\vec{h}}\big).
$$

\begin{proposition}
\label{prop:bound-for-W}
For any $q\geq 1$, $m\geq 3$ and any probability density $f$
$$
\bE_f\bigg\{\sup_{\vec{h}\in \cH_m^d}\Big(\big|\widetilde{\cW}_{\vec{h}}-\widetilde{\cW}_{\vec{h}}(f)\big|-\tfrac{25}{192}\cW_{\vec{h}}(f)\Big)_+^q\bigg\}\leq C_{q}
[\ln(m)]^d m^{-4q},
$$
where $C_{q}=6^d2^{q+1}\Gamma(q+1)$.

\end{proposition}

\noindent The proof of the proposition is given in Section \ref{sec:Proof of auxiliary results}.

\paragraph{Proof of the theorem.}  We break the proof into several steps.

$\mathbf{1^0}.\;$ Since
$$
\widehat{\cW}_{\vec{h}}=\widetilde{\cW}_{\vec{h}}\vee \tfrac{256s\ln(m)}{mV_{\vec{h}}},\quad\; \cW_{\vec{h}}(f)=\widetilde{\cW}_{\vec{h}}(f) \vee \tfrac{256s\ln(m}{mV_{\vec{h}}},
$$
applying the trivial inequality
$$
\big|(u\vee w)-(v\vee w)\big|\leq |u-v|,\quad \forall u,v,w\in\bR,
$$
we deduce from Proposition  \ref{prop:bound-for-W} that
\begin{align}
\label{eq000:proof-th3}
\bE_f\bigg\{\sup_{\vec{h}\in H_m^d}\Big(\big|\widehat{\cW}_{\vec{h}}-\cW_{\vec{h}}(f)\big|-\tfrac{25}{192}\cW_{\vec{h}}(f)\Big)_+^q\bigg\}\leq C_{q}
[\ln(m)]^d m^{-4q}.
\end{align}

$\mathbf{2^0}.\;$ We note that $\widehat{J}_{\vec{h}}$ is a decoupled $U$-statistics with the kernel $|U|$. Hence Theorem \ref{th:non-random-upper-function} is applicable to $\widehat{J}_{\vec{h}}$. Moreover we note  that all quantities and constants  appeared in the assertion of the latter theorem remain the same if one replaces $U$ by $|U|$.
Thus, recalling that $J_{\vec{h}}(f)=\bE_f\{\widehat{J}_{\vec{h}}\}$, we have in view of Corollary \ref{cor1:th2}
\begin{align}
\label{eq1:proof-th3}
\bE_f\bigg\{\sup_{\vec{h}\in\cH^d_m}\Big(\big|\hat{J}_{\vec{h}}-J_{\vec{h}}(f)\big|-U^*_{\vec{h}}(f)\Big)_+^q\bigg\}
\leq B_3[2\ln(m)]^dm^{-4q},
\end{align}
where  $B_3=\Lambda^*_q(|U|)=\Lambda^*_q(U)$.


Set as before $s=12q+2$ and introduce the following notations.
\begin{align*}
\eta&=\sup_{\vec{h}\in\cH^d_m}\big(\big|\widehat{J}_{\vec{h}}-J_{\vec{h}}(f)\big|- U^*_{\vec{h}}(f)\big)_+;
\nonumber\\*[2mm]
\chi&=\sup_{\vec{h}\in H_m^d}\Big(\big|\widehat{\cW}_{\vec{h}}-\cW_{\vec{h}}(f)\big|-\tfrac{25}{192}\cW_{\vec{h}}(f)\Big)_+
\end{align*}
Set also
$$
a=\tfrac{16s \varpi_U \ln(m)}{m},\;\;\; b=\tfrac{16s\varpi_U\ln(m) }{m^2V_{\vec{h}}},\;\;\; c=\tfrac{147s \varpi_U\ln(m)}{m}.
$$
We have
\begin{align}
\label{eq:proof-th3-100new}
\big|\cU_{\vec{h}}-U^*_{\vec{h}}(f)\big|&\leq \sqrt{a}\Big|\sqrt{\widehat{W}_{\vec{h}} \widehat{J}_{\vec{h}}}-\sqrt{W^*_{\vec{h}}(f) J_{\vec{h}}(f)}\Big|+\sqrt{b}\Big|\sqrt{\widehat{J}_{\vec{h}}}-\sqrt{ J_{\vec{h}}(f)}\Big|
\nonumber\\*[2mm]
&\quad  +c\big|\widehat{W}_{\vec{h}}-W^*_{\vec{h}}(f) \big|.
\end{align}

$\mathbf{3^0}.\;$ We obviously have
\begin{align*}
\sqrt{b}\Big|\sqrt{\widehat{J}_{\vec{h}}}-\sqrt{ J_{\vec{h}}(f)}\Big|&\leq \sqrt{bU^*_{\vec{h}}(f)}+\sqrt{b\eta}
\nonumber\\*[2mm]
&\leq (b/2)+2^{-1}U^*_{\vec{h}}(f)+\sqrt{b\eta}.
\end{align*}
Noting that $\Upsilon_{\vec{h}}(m)\geq 17$ for any $\vec{h}\in\mH^d$ and $m\geq m^*$ we assert that
$$
b\leq \cU_{\vec{h}}\wedge U^{*}_{\vec{h}}(f)
$$
and therefore
\begin{align}
\label{eq:proof-th3-200new}
\sqrt{b}\Big|\sqrt{\widehat{J}_{\vec{h}}}-\sqrt{ J_{\vec{h}}(f)}\Big|\leq 2^{-1}\big[\cU_{\vec{h}}\wedge U^{*}_{\vec{h}}(f)+U^{*}_{\vec{h}}(f)\big]
+\sqrt{b\eta}.
\end{align}
Note also another bound which can be useful in the sequel.

Since $\sqrt{bU^*_{\vec{h}}(f)}\leq b+4^{-1}U^*_{\vec{h}}(f)$ we obtain
\begin{align}
\label{eq:proof-th3-2000new}
\sqrt{b}\Big|\sqrt{\widehat{J}_{\vec{h}}}-\sqrt{ J_{\vec{h}}(f)}\Big|\leq \tfrac{16s\varpi_U\ln(m) }{m^2V_{\vec{h}}}+4^{-1}U^*_{\vec{h}}(f)
+\sqrt{b\eta}.
\end{align}

Consider now separately two cases.

\vskip0.1cm

$\mathbf{4^0}.\;$ Assume that $\vec{h}\in H^d_m$. In this case  we have
\begin{align*}
\sqrt{a}\Big|\sqrt{\widehat{W}_{\vec{h}} \widehat{J}_{\vec{h}}}&-\sqrt{W^*_{\vec{h}}(f) J_{\vec{h}}(f)}\Big|=(a/\alpha_m)^{\frac{1}{2}}
\Big|\sqrt{\widehat{J}_{\vec{h}}}-\sqrt{ J_{\vec{h}}(f)}\Big|
\nonumber\\*[2mm]
&\leq\sqrt{(a/\alpha_m)U^*_{\vec{h}}(f)}+\sqrt{(a/\alpha_m)\eta}
\nonumber\\*[2mm]
&\leq \tfrac{a}{2\alpha_m}+2^{-1}U^*_{\vec{h}}(f)+\sqrt{(a/\alpha_m)\eta}.
\end{align*}
Noting that
$$
\tfrac{a}{\alpha_m}=\tfrac{16s \varpi_U \ln(m)}{\alpha_m m}\leq\tfrac{147s \varpi_U\ln(m)}{\alpha_m m}\leq \cU_{\vec{h}}\wedge U^{*}_{\vec{h}}(f)
$$
we get
\begin{align}
\label{eq:proof-th3-300new}
\sqrt{a}\Big|\sqrt{\widehat{W}_{\vec{h}} \widehat{J}_{\vec{h}}}&-\sqrt{W^*_{\vec{h}}(f) J_{\vec{h}}(f)}\Big|
\nonumber\\*[2mm]
&\leq 2^{-1}\big[\cU_{\vec{h}}\wedge U^{*}_{\vec{h}}(f)+U^{*}_{\vec{h}}(f)\big]+\sqrt{(a/\alpha_m)\eta}.
\end{align}
It remains to note that $\widehat{W}_{\vec{h}}=W^*_{\vec{h}}(f)$ in the considered case and we deduce from (\ref{eq:proof-th3-100new}), (\ref{eq:proof-th3-200new}) and (\ref{eq:proof-th3-300new}) that
\begin{align*}
\big|\cU_{\vec{h}}-U^*_{\vec{h}}(f)\big|
\leq 2^{-1}\big[\cU_{\vec{h}}\wedge U^{*}_{\vec{h}}(f)+U^{*}_{\vec{h}}(f)\big]+\sqrt{b\eta}+\sqrt{(a/\alpha_m)\eta}.
\end{align*}
It yields
\begin{align*}
\sup_{\vec{h}\in H_m^d}\big(\cU_{\vec{h}}-2U^*_{\vec{h}}(f)\big)_+
\leq A_1\sqrt{\eta},\quad\;\sup_{\vec{h}\in H_m^d}\big(U^*_{\vec{h}}-3\cU_{\vec{h}}(f)\big)_+
\leq 2A_1\sqrt{\eta}.
\end{align*}
It remains to note that
$$
\bE_f\big\{\eta^{q/2}\big\}\leq \sqrt{\bE_f\{\eta^{q}\}}\leq  \sqrt{B_2}\;[2\ln(m)]^{d/2}m^{-2q}
$$
in view of (\ref{eq1:proof-th3}) and we get
\begin{eqnarray}
\label{eq:proof-th3-400new}
&&\bE_f\bigg\{\sup_{\vec{h}\in H_m^d}\big(\cU_{\vec{h}}-2U^*_{\vec{h}}(f)\big)^q_+\bigg\}\leq B_4[2\ln(m)]^{d/2}m^{-2q};
\\*[2mm]
&&\bE_f\bigg\{\sup_{\vec{h}\in H_m^d}\big(U^*_{\vec{h}}-3\cU_{\vec{h}}(f)\big)^q_+\bigg\}\leq B_5[2\ln(m)]^{d/2}m^{-2q}.
\label{eq:proof-th3-401new}
\end{eqnarray}

$\mathbf{5^0}.\;$  It is easily seen that there exists $A_0$ such that for any $m\geq 9$, $f\in\bB_\infty\big(\alpha_m^{-1}\big)$ and  $\vec{h}\in \bar{H}^d_m$
\begin{equation}
\label{eq:proof-th3-40005new}
\max\big\{U^*_{\vec{h}}(f),\cW_{\vec{h}}(f),J_{\vec{h}}(f)\big\}\leq A_0\alpha^{-1}_m\leq  A_0\ln(m).
\end{equation}
Assume that $\vec{h}\in \bar{H}^d_m$. We have in this case
\begin{align*}
&\big|\widehat{W}_{\vec{h}}  \widehat{J}_{\vec{h}}-W^*_{\vec{h}}(f) J_{\vec{h}}(f)\big|=\big|\widehat{\cW}_{\vec{h}}  \widehat{J}_{\vec{h}}-\cW_{\vec{h}}(f) J_{\vec{h}}(f)\big|
\nonumber\\*[3mm]
&\leq \big|\widehat{\cW}_{\vec{h}}-\cW_{\vec{h}}(f)\big| \big|\widehat{J}_{\vec{h}}-J_{\vec{h}}(f)\big|+\big|\widehat{J}_{\vec{h}}-J_{\vec{h}}(f)\big|\cW_{\vec{h}}(f)+\big|\widehat{\cW}_{\vec{h}}
-\cW_{\vec{h}}(f)\big|J_{\vec{h}}(f)
\nonumber\\*[3mm]
&\leq\big(U^*_{\vec{h}}(f)+\eta\big)\big(\tfrac{25}{192}\cW_{\vec{h}}(f)+\chi\big)+U^*_{\vec{h}}(f)\cW_{\vec{h}}(f)+\cW_{\vec{h}}(f)\eta
\nonumber\\*[3mm]
&\hskip2cm +\tfrac{25}{192}\cW_{\vec{h}}(f)J_{\vec{h}}(f)+J_{\vec{h}}(f)\chi
\nonumber\\*[3mm]
&\leq \tfrac{217}{192}U^*_{\vec{h}}(f)\cW_{\vec{h}}(f)+\tfrac{25}{192}\cW_{\vec{h}}(f)J_{\vec{h}}(f)
 +A_1\eta+A_2\chi+\eta\chi.
\end{align*}
Thus, we have
\begin{align}
\label{eq:proof-th3-402new}
\sqrt{a}\Big|&\sqrt{\widehat{W}_{\vec{h}} \widehat{J}_{\vec{h}}}-\sqrt{W^*_{\vec{h}}(f) J_{\vec{h}}(f)}\Big|\leq 4^{-1}U^*_{\vec{h}}(f)
+\tfrac{217}{882}\Big[\tfrac{147s \cW_{\vec{h}}(f)\varpi_U \ln(m)}{m}\Big]
\nonumber\\*[2mm]
&+\tfrac{5}{13}\Big[\tfrac{16s \varpi_U \cW_{\vec{h}}(f)J_{\vec{h}}(f)\ln(m)}{m}\Big]^{\frac{1}{2}}
 +A_3\sqrt{\eta}+A_4\sqrt{\chi}+2^{-1}(\eta+\chi).
\end{align}
Here we have also used (\ref{eq:proof-th3-40005new}).
At last
$$
c\big|\widehat{W}_{\vec{h}}-W^*_{\vec{h}}(f) \big|=c\big|\widehat{\cW}_{\vec{h}}-\cW_{\vec{h}}(f) \big|\leq \tfrac{25}{192}\Big[\tfrac{147s \varpi_U \cW_{\vec{h}}(f)\ln(m)}{m}\Big]+c\chi.
$$
This yields together with (\ref{eq:proof-th3-100new}), (\ref{eq:proof-th3-2000new}) and (\ref{eq:proof-th3-402new})
\begin{align*}
\big|\cU_{\vec{h}}-U^*_{\vec{h}}(f)\big|&\leq  2^{-1}U^*_{\vec{h}}(f)
+\big(\tfrac{217}{882}+\tfrac{25}{192}\big)\Big[\tfrac{147s \cW_{\vec{h}}(f)\varpi_U \ln(m)}{m}\Big]
\nonumber\\*[2mm]
&\;+\tfrac{5}{13}\Big[\tfrac{16s \varpi_U \cW_{\vec{h}}(f)J_{\vec{h}}(f)\ln(m)}{m}\Big]^{\frac{1}{2}}
 +\tfrac{16s\varpi_U\ln(m) }{m^2V_{\vec{h}}}
\nonumber\\*[2mm]
&
 \;+A_5\sqrt{\eta}+A_4\sqrt{\chi}+2^{-1}(\eta+\chi).
\end{align*}
Noting that $\tfrac{217}{882}+\tfrac{25}{192}\leq \tfrac{5}{13}$ and that
$$
\tfrac{16s\varpi_U\ln(m) }{m^2V_{\vec{h}}}\leq 1.6s^{-1}\Big[\tfrac{10s\varpi_U\ln^2(m) }{m^2V_{\vec{h}}}\Big]\leq \tfrac{5}{13}\Big[\tfrac{10s\varpi_U\ln^2(m) }{m^2V_{\vec{h}}}\Big]
$$
because $s=12q+2\geq 10$, we obtain finally that
\begin{align}
\label{eq:proof-th3-403new}
\big|\cU_{\vec{h}}-U^*_{\vec{h}}(f)\big|&\leq  2^{-1}U^*_{\vec{h}}(f)+\tfrac{5}{13}U^*_{\vec{h}}(f)
+A_5\sqrt{\eta}+A_4\sqrt{\chi}+2^{-1}(\eta+\chi).
\nonumber\\*[2mm]
&=\tfrac{23}{26}U^*_{\vec{h}}(f)+A_6\xi,
\end{align}
where we put for brevity $\xi=\sqrt{\eta}+\sqrt{\chi}+\eta+\chi$. We deduce from (\ref{eq:proof-th3-403new})
\begin{align*}
\sup_{\vec{h}\in \bar{H}_m^d}\big(\cU_{\vec{h}}-2U^*_{\vec{h}}(f)\big)_+
\leq A_6\xi,\quad\;\sup_{\vec{h}\in \bar{H}_m^d}\big(U^*_{\vec{h}}-9\cU_{\vec{h}}(f)\big)_+
\leq 9A_6\xi.
\end{align*}
Applying (\ref{eq1:proof-th3}) and Proposition \ref{prop:bound-for-W} we obtain similarly to (\ref{eq:proof-th3-400new}) and  (\ref{eq:proof-th3-401new})
\begin{eqnarray}
\label{eq:proof-th3-404new}
&&\bE_f\bigg\{\sup_{\vec{h}\in \bar{H}_m^d}\big(\cU_{\vec{h}}-2U^*_{\vec{h}}(f)\big)^q_+\bigg\}\leq B_6[2\ln(m)]^{d/2}m^{-2q};
\\*[2mm]
&&\bE_f\bigg\{\sup_{\vec{h}\in \bar{H}_m^d}\big(U^*_{\vec{h}}-9\cU_{\vec{h}}(f)\big)^q_+\bigg\}\leq B_7[2\ln(m)]^{d/2}m^{-2q}.
\label{eq:proof-th3-405new}
\end{eqnarray}
The assertion of the theorem with $B_1=B_4\vee B_6$ and $B_2=B_5\vee B_7$ follows now from (\ref{eq:proof-th3-400new}), (\ref{eq:proof-th3-401new}), (\ref{eq:proof-th3-404new}) and
(\ref{eq:proof-th3-405new}).
\epr

\subsubsection{Proof of Corollary \ref{cor1:th2}.} Set
$
\Pi_{\vec{k}}=\otimes_{i=1}^d \big(k_ih_it,(k_i+1)h_it\big], \vec{k}\in\bZ^d.
$
Since $\{\Pi_{\vec{k}}, \vec{k}\in\bZ^d\}$ forms the partition of $\bR^d$ on obtain
$$
W_{\vec{h}}(f)=\sup_{\vec{k}\in\bZ^d}\sup_{x\in\Pi_{\vec{k}}}\int_{\bR^d}\big|U_{\vec{h}}(x-y)\big|f(y)\rd y.
$$
Noting that for any $x\in\Pi_{\vec{k}}$ and $y\in\bR^d$
$$
\big|U_{\vec{h}}(x-y)\big|\leq \varpi_UV_{\vec{h}}^{-1}\mathrm{1}_{\Pi^*_{\vec{k}}}\big(\tfrac{y}{\vec{h}}\big)
$$
we can state that
$
W_{\vec{h}}(f)\leq \cW_{\vec{h}}(f)
$
for all $\vec{h}\in\bR^d_+$ and density $f$.

On the other hand $W_{\vec{h}}(f)\leq \varpi_U \alpha_m^{-1}$  for any $\vec{h}\in\bR^d_+$,  $f\in\bB_\infty\big(\alpha_m^{-1}\big)$. Both bounds yield together  $W_{\vec{h}}(f)\leq W^*_{\vec{h}}(f)$. Additionally we replace multiplicative constant $23$ in the third term of $U^*_{\vec{h}}(f)$ by larger constant $147$. All of this yields  $U_{\vec{h}}(f)<U^*_{\vec{h}}(f)$ and the assertion of the corollary follows from Theorem \ref{th:non-random-upper-function}.
\epr

\subsubsection{Proof of Corollary \ref{cor2:th2}.} First we remark that that
$$
\sup_{\vec{k}\in\bZ^d}V_{\vec{h}}^{-1}\int_{\bR^d}\mathrm{1}_{\Pi^*_{\vec{k}}}\big(\tfrac{y}{\vec{h}}\big)f(y)\rd y\leq 3^{d}\|f\|_\infty
$$
for any $\vec{h}\in\bR^d_+$ and any density $f\in\bB_\infty\big(\alpha_m^{-1}\big)$. It yields
\begin{align}
\label{eq1:proof-cor1}
W^{*}_{\vec{h}}(f)\leq 3^{d}\kappa_{\vec{h}}(f),\quad \vec{h}\in\cH_m^d.
\end{align}
Next, in view of the Cauchy-Schwartz inequality
\begin{align*}
J_{\vec{h}}(f)&:=\int_{\bR^d}\int_{\bR^d}|U_{\vec{h}}(y-z)|f(y)f(z)\rd y\rd z
\\*[2mm]
&\leq \|f\|_2 \bigg[\int_{\bR^d}\bigg(\int_{\bR^d}|U_{\vec{h}}(y-z)|f(y)\rd y\bigg)^2\rd z\bigg]^{\frac{1}{2}}=
\|f\|_2\big\||U_{\vec{h}}|\star f\big\|_2.
\end{align*}
Applying the Young inequality  we obtain
\begin{align}
\label{eq2:proof-cor1}
J_{\vec{h}}(f)\leq \|U\|_1\|f\|_2^2\leq \varpi_U\|f\|_2^2.
\end{align}
We deduce from (\ref{eq1:proof-cor1}) and  (\ref{eq2:proof-cor1}) that $U_{\vec{h}}(f)\leq \Omega_q(U)\mU_{\vec{h}}(f)$ and, therefore,
the assertion of the corollary follows from   Theorem \ref{th:non-random-upper-function}.
\epr

\section{Proofs of auxiliary results from Section \ref{sec:proof-U-stat}.}

\paragraph{Proof of Proposition \ref{prop:bound-for-func-G}.}

We obviously have
\begin{align}
\label{eq0-proof-prop1-new00}
\bP_f\Big(\cup_{i=1}^m\big\{|G(Y_i,X^{(m)})&|> \rho_{\vec{h}}(f)\big\}\Big)\leq \sum_{i=1}^m\bP_f\big\{|G(Y_i,X^{(m)})|> \rho_{\vec{h}}(f)\big\}
\nonumber\\*[2mm]
&=m\int_{\bR^d}\bP_f\big\{|G(y,X^{(m)})|> \rho_{\vec{h}}(f)\big\}f(y)\rd y.
\end{align}
To get  last equality we used that that $X^{(m)}$ and $Y_i$
 are independent.

$\mathbf{1^0}.\;$ Noting that $G(y,X^{(m)})$ is the arithmetic mean of the independent, centered and bounded by $2\varpi_U V^{-1}_{\vec{h}}$ random variables we deduce from
the first assertion of Lemma \ref{lem:bernstein-new}
 that
\begin{eqnarray*}
 &&\bP_f\Big\{|G(y,X^{(m)})|> \Big(\tfrac{2s \sigma^2(y)\ln(m)}{m}\Big)^{\frac{1}{2}}+\tfrac{4s\varpi_U \ln(m)}{3mV_{\vec{h}}}\Big\}\leq 2m^{-s}.
\end{eqnarray*}
where
\begin{eqnarray*}
\sigma^2(y)=\bE_f\big(U^2_{\vec{h}}(y-X_1)\big)=\int_{\bR^d}U^2_{\vec{h}}(y-z)f(z)\rd z
\leq\varpi_U V^{-1}_{\vec{h}}\;W_{\vec{h}}(f).
\end{eqnarray*}
Thus, we obtain  for any $\vec{h}\in\mH^d$
 \begin{equation}
\label{eq1-proof-prop1-new00}
\bP_f\Big\{|G(y,X^{(m)})|> \Big(\tfrac{2s\varpi_U W_{\vec{h}}(f)\ln(m)}{mV_{\vec{h}}}\Big)^{\frac{1}{2}}+\tfrac{4s\varpi_U \ln(m)}{3mV_{\vec{h}}}\Big\}\leq 2m^{-s}.
\end{equation}
We deduce the assertion of the proposition in the case $\vec{h}\in\bar{\mH}^d_*$ from (\ref{eq0-proof-prop1-new00}) and (\ref{eq1-proof-prop1-new00}).

\smallskip

$\mathbf{2^0}.\;$
Denote by $U^*_{\vec{h}}(x)=|U(\frac{x_1}{h_1},\ldots,\frac{x_d}{h_d})|$ and remark that for any $y\in\bR^d$
\begin{equation*}
|G(y,X^{(m)})|\leq  \big\|I_{\vec{h}}\big\|_\infty+\tfrac{\pi(y,X^{(m)})}{mV_{\vec{h}}},\;\quad \pi(y,X^{(m)})=\sum_{j=1}^m\;U^*_{\vec{h}}(y-X_j).
\end{equation*}
Noting that $\|I_{\vec{h}}\|_\infty\leq \varpi_U\alpha^{-1}_m$ for any $f\in\bB_\infty\big(\alpha^{-1}_m\big)$
and setting
$$
\mu_{\vec{h}}(f)=\varpi_U\alpha^{-1}_m+\tfrac{Z}{mV_{\vec{h}}},\quad\; Z=\tfrac{\varpi_U s\ln(m)}{|\ln(mV_{\vec{h}})|}
$$
 we obviously have
\begin{align}
\label{eq4-proof-prop1-new00}
\bP_f\Big(|G(y,X^{(m)})|> \mu_{\vec{h}}(f)\big\}\Big)\leq \bP_f\big\{\pi(y,X^{(m)})\geq Z\big\}.
\end{align}
We remark that  Lemma \ref{lem:lb-deviation-new} is applicable to $\pi(y,X^{(m)})$ with $k=m$ and
$$
\e_j=U^*_{\vec{h}}(y-X_j),\;\; T=\varpi_U,\;\; \alpha=\int_{\bR^d}U^*_{\vec{h}}(y-x)f(x)\rd x,\;\; z=Z.
$$
Indeed, for any $\vec{h}\in\mH^d_*$ and  $$f\in\bB_\infty\big(\alpha^{-1/2}_m\big)$$
$$
m\alpha\leq \|U\|_1 mV_{\vec{h}}\alpha^{-1}_m\leq\varpi_U \alpha_m\leq \varpi_U.
$$
On the other hand, since $mV_{\vec{h}}\geq m^{-1}$ in view of the definition of $\mH^d$ we have
$Z\geq \varpi_U s\geq \varpi_U$.

\noindent Thus, we derive from  Lemma \ref{lem:lb-deviation-new}, putting for brevity  $z=Z/\varpi_U$,
that
for any $\vec{h}\in\mH^d_*$ and $y\in\bR^d$
\begin{align}
\label{eq5-proof-prop1-new00}
\bP_f\big\{\pi(y,X^{(m)})\geq Z\big\}
&\leq  \Big(\tfrac{mV_{\vec{h}}\varpi_U \alpha^{-1}_m}{Z}\Big)^{\frac{Z}{\varpi_U}} e^{\frac{Z}{\varpi_U}}
\leq \Big(\tfrac{\sqrt{mV_{\vec{h}}}}{z}\Big)^{z}e^{z}
\nonumber\\*[2mm]
&= e^{-(s/2)\ln(m)-z\ln(z)+z}
\leq  m^{-\frac{s}{2}}.
\end{align}
To get the last inequality we used $z\geq 1$ and that $z\mapsto z\ln(z)-z$ in increasing function on $[1,\infty)$.

Since $\alpha^{-2}_m\leq (mV_{\vec{h}})^{-1}$ in view of the definition of  $\mH_*^d$ we have,
$$
\mu_{\vec{h}}(f)\leq \tfrac{2s\varpi_U \ln(m)}{mV_{\vec{h}}|\ln(mV_{\vec{h}})|}.
$$
Thus, we deduce from (\ref{eq4-proof-prop1-new00}) and (\ref{eq5-proof-prop1-new00}) that for any $\vec{h}\in\mH_*^d$ and $y\in\bR^d$
\begin{align}
\label{eq6-proof-prop1-new00}
\bP_f\Big(|G(y,X^{(m)})|> \tfrac{2s\varpi_U \ln(m)}{mV_{\vec{h}}|\ln(mV_{\vec{h}})|}\Big)\leq  m^{-s/2}.
\end{align}
The assertion of the proposition in the case $\vec{h}\in\mH_*^d$ follows now from (\ref{eq0-proof-prop1-new00}) and (\ref{eq6-proof-prop1-new00}).
\epr

\paragraph{Proof of Proposition \ref{prop:bound-for-second-term2}.}
We begin the proof with presenting the version  of Shur lemma,  see \cite{Folland}, Theorem 6.18, adapted to our purposes.

\begin{lemma}[Shur lemma.]
Let $\bR^d$ be equipped with two $\sigma$-finite measures $\mu$ and $\nu$ and let  $p\geq 1$.
Assume that $\lambda:\bR^d\times\bR^d \to \bR$ and $g:\bR^d\to\bR$ be such that $g\in\bL_p(\bR^d,\nu)$ and
$$
\mathbf{C}_\lambda:=\Big[\sup_{x\in\bR^d}\int_{\bR^d}|\lambda(x,y)|\nu(\rd y)\Big]\bigvee\Big[\sup_{y\in\bR^d}\int_{\bR^d}|\lambda(x,y)|\mu(\rd x)\Big]<\infty.
$$
Then, $\int_{\bR^d}\lambda(\cdot,y)g(y)\nu(\rd y)\in \bL_p(\bR^d,\mu)$ and
$$
\Big\|\int_{\bR^d}\lambda(\cdot,y)g(y)\nu(\rd y)\Big\|_p\leq \mathbf{C}_\lambda\|g\|_p.
$$
\end{lemma}

\noindent Denote by $\bB_{2,f}$ the following set of functions.
$$
\bB_{2,f}=\bigg\{\ell:\bR^d\to\bR:\;\int_{\bR^d}\ell^2(y)f(y)\rd y\leq 1\bigg\}.
$$
We have in view of the duality argument
\begin{align*}
\upsilon\big(X^{(m)}\big)&=\bigg[\int_{\bR^d}\big[G\big(y,X^{(m)}\big)\big]^2 f(y)\rd y\bigg]^{\frac{1}{2}}
\nonumber\\*[2mm]
&=\sup_{\ell\in\bB_{2,f}}\int_{\bR^d}\ell(y)G\big(y,X^{(m)}\big) f(y)\rd y
\nonumber\\*[2mm]
&=\tfrac{1}{m}\sup_{\ell\in\bB_{2,f}}\sum_{j=1}^m\int_{\bR^d}\ell(y) \big[U_{\vec{h}}\big(y-X_j\big)-I_{\vec{h}}(y)\big] f(y)\rd y.
\end{align*}
Since $\bB_{2,f}$ is the unit ball in a Banach space, there exists a countable dense subset, say  $\bB^*_{2,f}$, of $\bB_{2,f}$.
Obviously
$$
\ell\mapsto \sum_{j=1}^m\int_{\bR^d}\ell(y) \big[U_{\vec{h}}\big(y-X_j\big)-I_{\vec{h}}(y)\big] f(y)\rd y
$$
is continuous functional and, therefore,
 \begin{align}
\label{eq181-proof-prop1-new00}
\upsilon\big(X^{(m)}\big)=\tfrac{1}{m}\sup_{\ell\in\bB^*_{2,f}}\sum_{j=1}^m\int_{\bR^d}\ell(y) \big[U_{\vec{h}}\big(y-X_j\big)-I_{\vec{h}}(y)\big] f(y)\rd y.
\end{align}
Setting for any $x\in\bR^d$
$$
S(x)=\tfrac{1}{m}\int_{\bR^d}\ell(y) \big[U_{\vec{h}}\big(y-x\big)-I_{\vec{h}}(y)\big] f(y)\rd y
$$
and denoting by $\cS=\big\{S:\;\ell\in\bB^*_{2,f}\big\}$ we obtain from (\ref{eq181-proof-prop1-new00})
 \begin{align}
\label{eq181-2-proof-prop1-new00}
\upsilon\big(X^{(m)}\big)=\sup_{S\in\cS}\sum_{j=1}^m S(X_j)=:Z.
\end{align}
We note that $\bE_f\{S(X_1)\}=0$ for any $S\in\cS$, since
\begin{align}
\label{eq181-1-proof-prop1-new00}
\bE_f\big\{U_{\vec{h}}\big(y-x\big)\big\}=I_{\vec{h}}(y),\quad \forall y\in\bR^d.
\end{align}
Moreover, by Cauchy-Schwartz inequality, for any $S\in\cS$
$$
\|S\|^2_\infty\leq \tfrac{1}{m^2}\sup_{x\in\bR^d}\int_{\bR^d}\big[U_{\vec{h}}\big(y-x\big)-I_{\vec{h}}(y)\big]^2 f(y)\rd y\leq \tfrac{4\varpi_U W_{\vec{h}}(f)}{m^2V_{\vec{h}}}=:\mathbf{T}^2.
$$
The idea is to apply the inequality (\ref{eq:Bousquet-expectation}) to the random variable $Z$.

First of all we remark that in view of (\ref{eq181-1-proof-prop1-new00})
\begin{align}
\label{eq186-proof-prop1-new00}
\bE_f\{Z\}= \bE_f\big\{\upsilon\big(X^{(m)}\big)\big\}\leq \big[\bE_f\big\{\upsilon^2\big(X^{(m)}\big)\big\}\big]^{\frac{1}{2}}.
\end{align}
We obviously have
\begin{align*}
\upsilon^2\big(X^{(m)}\big)&=
\tfrac{1}{m^2}\sum_{j=1}^m\int_{\bR^d}\big[U_{\vec{h}}\big(y-X_j\big)-I_{\vec{h}}(y)\big]^2f(y)\rd y
\nonumber\\*[2mm]
&\quad +\tfrac{1}{m^2}\sum_{i,j=1,j\neq i}^m\int_{\bR^d}\big[U_{\vec{h}}\big(y-X_i\big)-I_{\vec{h}}(y)\big]\big[U_{\vec{h}}\big(y-X_j\big)\big]f(y)\rd y
\end{align*}
and, therefore,
\begin{align*}
\bE_f\big\{\upsilon^2\big(X^{(m)}\big)\big\}&=
\tfrac{1}{m}\int_{\bR^d}\bE_f\big[U_{\vec{h}}\big(y-X_1\big)-I_{\vec{h}}(y)\big]^2f(y)\rd y
\nonumber\\*[2mm]
&\leq\tfrac{1}{m}\int_{\bR^d}\int_{\bR^d}U^{2}_{\vec{h}}\big(y-x\big)f(x)f(y)\rd x\rd y
\nonumber\\*[2mm]
&\leq \tfrac{\varpi_U}{mV_{\vec{h}}} \int_{\bR^d}\int_{\bR^d}\big|U_{\vec{h}}\big(y-x\big)\big|f(y)f(z)\rd y\rd z=\tfrac{\varpi_U J_{\vec{h}}(f)}{mV_{\vec{h}}}.
\end{align*}
It yields together with (\ref{eq186-proof-prop1-new00})
\begin{align}
\label{eq186-1-proof-prop1-new00}
\bE_f\{Z\} \leq \Big[\tfrac{\varpi_U J_{\vec{h}}(f)}{mV_{\vec{h}}}\Big]^{\frac{1}{2}}.
\end{align}
Note at last that
\begin{align*}
\bE_f\big\{S^2&(X_1)\big)\big\}=
\tfrac{1}{m^2}\int_{\bR^d}\bigg[\int_{\bR^d}\ell(y)\big[U_{\vec{h}}\big(y-x\big)-I_{\vec{h}}(y)\big]f(y)\rd y\bigg]^2f(x)\rd x.
\end{align*}
Applying Shur lemma with $p=2$, $\nu(dy)=f(y)\rd y$, $\mu(dx)=f(x)\rd x$ and
$$
\lambda(x,y)=U_{\vec{h}}\big(y-x\big)-I_{\vec{h}}(y),\quad g(y)=\ell(y),
$$
we obtain for any $S\in\cS$
\begin{align*}
\bE_f\big\{S^2&(X_1)\big)\big\}\leq \mathbf{C}^2_\lambda m^{-2},
\end{align*}
where, remind, $\mathbf{C}_\lambda$ is given by
\begin{align*}
\mathbf{C}_\lambda&:=\Big[\sup_{x\in\bR^d}\int_{\bR^d}\big|U_{\vec{h}}\big(y-x\big)-I_{\vec{h}}(y)\big|f(y)\rd y
\Big]
\nonumber\\*[2mm]
&\quad\;\vee\Big[\sup_{y\in\bR^d}\int_{\bR^d}\big|U_{\vec{h}}\big(y-x\big)-I_{\vec{h}}(y)\big|f(x)\rd x\Big]
\nonumber\\*[2mm]
&\leq 2\|I_{\vec{h}}\|_\infty\leq 2W_{\vec{h}}(f).
\end{align*}
Thus, applying (\ref{eq:Bousquet-expectation}) with $p=q$, $\varsigma^2=4m^{-2}W^2_{\vec{h}}(f)$, $T=\mathbf{T}$, $t=s\ln(m)$, noting that $\mathbf{T}\leq 2\varpi_Um$  and
$$
2\mathrm{E}\{Z\}+\varsigma\sqrt{2mt}+(5/3)Tt\leq c_{\vec{h}}(f),
$$
and taking into account (\ref{eq186-1-proof-prop1-new00}) we come to
$$
\bE_f\big\{\big(\upsilon\big(X^{(m)}\big)-c_{\vec{h}}(f)\big)_+^q\big\}\leq \Gamma(q+1)\big[6\varpi_U\big]^{q}m^{\frac{3q}{2}-s}.
$$
The proposition is proved.
\epr

\paragraph{Proof of Proposition \ref{prop:bound-for-W}.} For any  $\vec{k}\in\bZ^d$ set
$$
\mathbf{f}_{\vec{k}}=\int_{\bR^d}\mathrm{1}_{\Pi^*_{\vec{k}}}\big(\tfrac{y}{\vec{h}}\big)f(y)\rd y,\quad \zeta_{\vec{k}}=\tfrac{1}{mV_{\vec{h}}}\sum_{j=1}^m \big[\mathrm{1}_{\Pi^*_{\vec{k}}}\big(\tfrac{X_j}{\vec{h}}\big)-\mathbf{f}_{\vec{k}}\big].
$$
Obviously,
$
\bE_f\big\{\zeta_{\vec{k}}\big\}=0.
$
Introduce
$$
\bK^d=\big\{\vec{k}\in\bZ^d:\; \mathbf{f}_{\vec{k}}\geq m^{-6q-2}\big\},\quad\; \bar{\bK}^d=\bZ^d\setminus\bK^d.
$$

$\mathbf{1^0a}.\;$ Set as before $s=12q+2$.

Applying the second assertion of Lemma \ref{lem:bernstein-new} with  $\e_j=V^{-1}_{\vec{h}}\big[\mathrm{1}_{\Pi^*_{\vec{k}}}\big(\tfrac{X_j}{\vec{h}}\big)-\mathbf{f}_{\vec{k}}\big]$, $p=q$, $T=V^{-1}_{\vec{h}}$
and $v^2=\mathbf{f}_{\vec{k}}V^{-2}_{\vec{h}}$
to $\zeta_{\vec{k}}$, we obtain for any $\vec{k}\in\bK^d$
$$
\bE_f\Big\{\Big(\big|\zeta_{\vec{k}}\big|-\big[\tfrac{4s \mathbf{f}_{\vec{k}}\ln(m)}{mV^{2}_{\vec{h}}}\big]^{\frac{1}{2}}-\tfrac{4s\ln(m)}{3mV_{\vec{h}}}\Big)_+^q\Big\} \leq 2^{q+1}\Gamma(q+1) m^{-s+\frac{3q}{2}}.
$$
Here we also used that $V_{\vec{h}}\geq m^{-2}$ for any $\vec{h}\in\mH^d$. It yields for any $\vec{k}\in\bK^d$
$$
\bE_f\Big\{\Big(\big|\zeta_{\vec{k}}\big|-\big[\tfrac{4s \mathbf{f}_{\vec{k}}\ln(m)}{mV^{2}_{\vec{h}}}\big]^{\frac{1}{2}}-\tfrac{4s\ln(m)}{3mV_{\vec{h}}}\Big)_+^q\Big\} \leq 2^{q+1}\Gamma(q+1) m^{-\frac{9q}{2}}\mathbf{f}_{\vec{k}}.
$$
Putting for brevity $\lambda_m=2^{q+1}\Gamma(q+1) m^{-\frac{9q}{2}}$ and denoting by
$$
\varrho_{\vec{h}}\big(f,\vec{k}\big)=\big[\tfrac{4s \mathbf{f}_{\vec{k}}\ln(m)}{mV^{2}_{\vec{h}}}\big]^{\frac{1}{2}}+\tfrac{4s\ln(m)}{3mV_{\vec{h}}}
$$
we obtain
\begin{align}
\label{eq1:proof-lem-bound-forW}
\bE_f\bigg\{\sup_{\vec{k}\in\bK^d}\Big(\big|\zeta_{\vec{k}}\big|-\varrho_{\vec{h}}\big(f,\vec{k}\big)\Big)_+^q\bigg\} &\leq \lambda_m\sum_{\vec{k}\in\bK^d}\mathbf{f}_{\vec{k}}
\nonumber\\*[2mm]
&\leq\lambda_m\int_{\bR^d}\bigg(\sum_{\vec{k}\in\bZ^d}\mathrm{1}_{\Pi^*_{\vec{k}}}\big(\tfrac{y}{\vec{h}}\big)\bigg)f(y)\rd y.
\end{align}

$\mathbf{1^0b}.\;$ Let us prove that
\begin{align}
\label{eq2:proof-lem-bound-forW}
\sum_{\vec{k}\in\bZ^d}\mathrm{1}_{\Pi^*_{\vec{k}}}(u)\leq 3^{d},\quad\forall u\in\bR^d.
\end{align}
Let us fix $u\in\bR^d$ and let $\vec{\ell}\in\bZ^d$ be such that\footnote{Obviously $\{\pi_{\vec{l}},\vec{l}\in\bZ^d\}$ forms the partition of $\bR^d$.}
$u\in\otimes_{r=1}^d \big(\ell_rt,(\ell_r+1)t\big]=:\pi_{\vec{\ell}}$.
In view of the definition of  $\Pi^*_{\vec{k}}$
$$
\Pi^*_{\vec{k}}\cap\pi_{\vec{\ell}}=\emptyset,  \quad \forall \vec{k}\in\bZ^d:\; \max_{r=1,\ldots,d}|k_r-\ell_r|>1.
$$
It remains to note that $\text{card}\{\vec{k}\in\bZ^d:\; \max_{r=1,\ldots,d}|k_r-\ell_r|\leq1\}=3^d$ for any $\vec{\ell}\in\bZ^d$ and, therefore,
(\ref{eq2:proof-lem-bound-forW}) is established.

We deduce finally from (\ref{eq1:proof-lem-bound-forW}) and (\ref{eq2:proof-lem-bound-forW}) that
\begin{align}
\label{eq3:proof-lem-bound-forW}
\bE_f\bigg\{\sup_{\vec{k}\in\bK^d}\Big(\big|\zeta_{\vec{k}}\big|-\varrho_{\vec{h}}\big(f,\vec{k}\big)\Big)_+^q\bigg\} &\leq 3^d2^{q+1}\Gamma(q+1) m^{-\frac{9q}{2}}.
\end{align}

$\mathbf{2^0}.\;$ Fix $\vec{k}\in\bar{\bK}^d$ and introduce the following random event
$$
\cA=\bigg\{\sum_{j=1}^m \mathrm{1}_{\Pi^*_{\vec{k}}}\big(\tfrac{X_j}{\vec{h}}\big)\geq 2\bigg\}.
$$
Let us prove that
\begin{align}
\label{eq4:proof-lem-bound-forW}
\bP_f(\cA)\leq  2m^{-6q}\mathbf{f}_{\vec{k}}.
\end{align}
We will apply Lemma \ref{lem:lb-deviation-new} with $k=m$, $\e_j=\mathrm{1}_{\Pi^*_{\vec{k}}}\big(\tfrac{X_j}{\vec{h}}\big)$, $T=1$, $\alpha=\mathbf{f}_{\vec{k}}$
and $z=2$. We remark $z=2\geq m^{-6q-1}\geq m\mathbf{f}_k=m\alpha$ and, therefore, Lemma \ref{lem:lb-deviation-new} is applicable. It yields
\[
 \bP_f(\cA) \leq (e/2)^2m^{2}\mathbf{f}_{\vec{k}}^2\leq 2m^{-6q}\mathbf{f}_{\vec{k}},\quad \forall \vec{k}\in\bar{\bK}^d.
\]
and (\ref{eq4:proof-lem-bound-forW}) is established.

Note that
$
\big|\zeta_{\vec{k}}\big|\mathrm{1}_{\bar{\cA}}\leq \big(mV_{\vec{h}})^{-1},
$
for any $\vec{k}\in\bar{\bK}^d$
and, therefore,
$$
\big(\big|\zeta_{\vec{k}}\big|-\varrho_{\vec{h}}\big(f,\vec{k}\big)\big)_+\mathrm{1}_{\bar{\cA}}=0.
$$
On the other hand $\big|\zeta_{\vec{k}}\big|\leq V^{-1}_{\vec{h}}\leq m^{2}$ if $\vec{k}\in\bar{\bK}^d$ and we obtain for any $\vec{k}\in\bar{\bK}^d$
\begin{align*}
\bE_f\bigg\{\Big(\big|\zeta_{\vec{k}}\big|-\varrho_{\vec{h}}\big(f,\vec{k}\big)\Big)_+^q\bigg\}\leq 2m^{-4q}\mathbf{f}_{\vec{k}}.
\end{align*}
It yields together with (\ref{eq2:proof-lem-bound-forW})
\begin{align}
\label{eq5:proof-lem-bound-forW}
\bE_f\bigg\{\sup_{\vec{k}\in\bar{\bK}^{d}}\Big(\big|\zeta_{\vec{k}}\big|-\varrho_{\vec{h}}\big(f,\vec{k}\big)\Big)_+^q\bigg\}\leq 3^{d}2m^{-4q}.
\end{align}
We deduce from  (\ref{eq3:proof-lem-bound-forW}) and (\ref{eq5:proof-lem-bound-forW}) that for any $\vec{h}\in\mH^d$
\begin{align}
\label{eq6:proof-lem-bound-forW}
\bE_f\bigg\{\sup_{\vec{k}\in\bZ^d}\Big(\big|\zeta_{\vec{k}}\big|-\varrho_{\vec{h}}\big(f,\vec{k}\big)\Big)_+^q\bigg\} &\leq 3^d2^{q+1}\Gamma(q+1) m^{-4q}.
\end{align}

$\mathbf{3^0}.\;$ Obviously
$
\big|\widetilde{\cW}_{\vec{h}}-\widetilde{\cW}_{\vec{h}}(f)\big|\leq \sup_{\vec{k}\in\bZ^d}\big|\zeta_{\vec{k}}\big|
$
and, therefore,
\begin{align}
\label{eq7:proof-lem-bound-forW}
\Big(\big|\widetilde{\cW}_{\vec{h}}-\widetilde{\cW}_{\vec{h}}(f)\big|-\sup_{\vec{k}\in\bZ^d}\varrho_{\vec{h}}\big(f,\vec{k}\big)\Big)_+\leq \sup_{\vec{k}\in\bZ^d}\Big(\big|\zeta_{\vec{k}}\big|-\varrho_{\vec{h}}\big(f,\vec{k}\big)\Big)_+.
\end{align}
It remains to note that
\begin{align*}
\sup_{\vec{k}\in\bZ^d}\varrho_{\vec{h}}\big(f,\vec{k}\big)=\big[\tfrac{4s \widetilde{\cW}_{\vec{h}}(f)\ln(m)}{mV_{\vec{h}}}\big]^{\frac{1}{2}}+\tfrac{4s\ln(m)}{3mV_{\vec{h}}}\leq \tfrac{25}{192}\cW_{\vec{h}}(f)
\end{align*}
and we deduce from (\ref{eq6:proof-lem-bound-forW}) and (\ref{eq7:proof-lem-bound-forW}) that for any $\vec{h}\in\mH^d$
\begin{align*}
\bE_f\bigg\{\Big(\big|\widetilde{\cW}_{\vec{h}}-\widetilde{\cW}_{\vec{h}}(f)\big|-\tfrac{25}{192}\cW_{\vec{h}}(f)\Big)_+^q\bigg\} &\leq 3^d2^{q+1}\Gamma(q+1) m^{-4q}.
\end{align*}
It yields finally
\begin{align*}
\bE_f\bigg\{\sup_{\vec{h}\in \cH_m^d}\Big(\big|\widetilde{\cW}_{\vec{h}}-\widetilde{\cW}
_{\vec{h}}(f)\big|-\tfrac{25}{192}\cW_{\vec{h}}(f)\Big)_+^q\bigg\} &\leq 6^d2^{q+1}\Gamma(q+1)[\ln(m)]^d m^{-4q}.
\end{align*}
Proposition is proved.
\epr

\paragraph{Proof of Lemma \ref{lem:lb-deviation-new}.}
Applying exponential Markov inequality we have
$$
\mathrm{P} \big(\cE_m \geq z\big)\leq e^{-z\lambda}\big[\mathrm{E}\big\{e^{\lambda\e_1}\big\}\big]^k,\quad\forall \lambda\geq 0.
$$
Using Taylor expansion we get
$$
\mathrm{E}\big\{e^{\lambda\e_1}\big\}=1+\sum_{r=1}^\infty \tfrac{\lambda^r \mathrm{E}\{\e^r_1\}}{r!}
\leq 1+\alpha\sum_{r=1}^\infty \tfrac{\lambda^r T^{r-1}}{r!}=1+\alpha T^{-1}\big(e^{T\lambda}-1\big).
$$
It yields since $1+y\leq e^y$
$$
\mathrm{P} \big(\cE_m \geq z\big)\leq e^{\psi(z)},\quad\; \psi(z)= \inf_{\lambda\geq 0}\big[-z\lambda+k\alpha T^{-1}(e^{T\lambda}-1)\big].
$$
It is easily seen that the minimum is attained  at  $\lambda^*=T^{-1}\ln(z/(k\alpha))$ and $\lambda^*\geq 0$ since $z\geq k\alpha$.
It yields
$$
\psi(z)=-zT^{-1}\ln\big(\tfrac{z}{k\alpha}\big)+T^{-1}(z-k\alpha).
$$
Lemma is proved.
\epr

 \bibliographystyle{agsm}

\end{document}